\documentclass[11pt,reqno]{amsart}
\usepackage{amsmath}
\usepackage{amscd}
\usepackage{comment}
\usepackage{amssymb}
\usepackage{graphicx}
\usepackage{subcaption}
\usepackage{color}

\usepackage[cmtip,arrow]{xy}
\usepackage{pb-diagram,pb-xy}

\def\corcommstyle#1{\bf\small\tt\textcolor{blue}{#1}}
\def\corcommredstyle#1{\bf\small\tt\textcolor{red}{#1}}

\def\corrl #1<<#2||#3>>{
\if\visiblecomments y
  \begin{quote} {\corcommstyle{ $<<$COMMENT$>>$ #1\marginpar{!!}\\$<<$OLD$<<$}} \end{quote}
  #2
  \begin{quote} {\corcommstyle{ ==NEW== }} \end{quote}
  #3
  \begin{quote} {\corcommstyle{ $>>$END$>>$ }} \end{quote}
 \else
  #3
 \fi
}

\long\def\longcorrl #1<<#2||#3>>{
\if\visiblecomments y
  \begin{quote} {\corcommstyle{ $<<$COMMENT$>>$ #1\marginpar{!!}\\$<<$OLD$<<$}} \end{quote}
  #2
  \begin{quote} {\corcommstyle{ ==NEW== }} \end{quote}
  #3
  \begin{quote} {\corcommstyle{ $>>$END$>>$ }} \end{quote}
 \else
  #3
 \fi
}

\def\corrq #1<<#2>>{
\if\visiblecomments y
  \begin{quote} {\corcommstyle{ $<<$COMMENT$>>$ #1\marginpar{!!}\\$<<$BEG$<<$}} \end{quote}
  #2
  \begin{quote} {\corcommstyle{ $>>$END$>>$ }} \end{quote}
 \else
  #2
 \fi
}

\long\def\longcorrq #1<<#2>>{
\if\visiblecomments y
  \begin{quote} {\corcommstyle{ $<<$COMMENT$>>$ #1\marginpar{!!}\\$<<$BEG$<<$}} \end{quote}
  #2
  \begin{quote} {\corcommstyle{ $>>$END$>>$ }} \end{quote}
 \else
  #2
 \fi
}

\def\corrc #1<<>>{
\if\visiblecomments y
  \begin{quote} {\corcommstyle{ $<<$COMMENT$>>$ #1\marginpar{!!}}} \end{quote}
\fi
}

\def\corrse #1<<>>{
\if\visiblecomments y
  \begin{quote} {\corcommstyle{ $<<$SECOND EDITION$>>$ #1\marginpar{!!}}} \end{quote}
\fi
}

\def\corre #1<<#2||#3>>{
\if\visiblecomments y
  #3\marginpar{\corcommstyle{ #1}}
 \else
  #3
 \fi
}

\long\def\longcorre #1<<#2||#3>>{
\if\visiblecomments y
  #3\marginpar{\corcommstyle{ #1}}
 \else
  #3
 \fi
}

\def\corrs #1<<#2||#3>>{
\if\visiblecomments y
  \textcolor{red}{#3}\marginpar{\corcommstyle{ #2 $\rightarrow$ #3\\ #1}}
 \else
  #3
 \fi
}

\def\corrm #1<<#2>>{
\if\visiblecomments y
  \textcolor{blue}{#2}\marginpar{\corcommstyle{#1}}
 \else
  #2
 \fi
}

\def\corrmr #1<<#2>>{
\if\visiblecomments y
  \textcolor{blue}{#2}\marginpar{\corcommredstyle{#1}}
 \else
  #2
 \fi
}

\def\corro #1<<#2||#3>>{
#2}

\def\corrn #1<<#2||#3>>{
#3}

\long\def\longcorro #1<<#2||#3>>{
#2}

\long\def\longcorrn #1<<#2||#3>>{
#3}

\long\def\underconstruction #1<<<#2>>>{
\if\visiblecomments y
  \begin{quote} {\corcommstyle{ $<<$UNDER CONSTRUCTION - BEGIN$>>$ #1\marginpar{!!}}} \end{quote}
  #2
  \begin{quote} {\corcommstyle{ $>>$UNDER CONSTRUCTION - END$>>$ }} \end{quote}
 \else
 \fi
}

\def\showcomments{
  \let\visiblecomments y
}

\def\hidecomments{
  \let\visiblecomments n
}

\usepackage{amssymb}

\def\refeq#1{\if\workingver y(\ref{#1})-[[#1]]\else(\ref{#1})\fi}
\def\refth#1{\if\workingver y\ref{#1}-[[#1]]\else\ref{#1}\fi}
\def\mylabel#1{\if\workingver y\label{#1}{\bf\ \ [[#1]]\ \ }\else\label{#1}\fi}
\def\mybibitem#1{\if\workingver y\bibitem{#1}{\bf\ \ [[#1]]\ \
}\else\bibitem{#1}\fi}

\def\articletheorems{
\newtheorem{thm}{Theorem}[section]
\newtheorem{lem}[thm]{Lemma}

\newtheorem{defn}[thm]{Definition}
\newtheorem{cor}[thm]{Corollary}
\newtheorem{prop}[thm]{Proposition}

}

\newcommand{\mto}{\multimap}
\newcommand{\pto}{\nrightarrow}

\renewcommand{\emptyset}{\varnothing}
\renewcommand{\rho}{\varrho}
\renewcommand{\epsilon}{\varepsilon}

\def\cA{\text{$\mathcal A$}}
\def\cB{\text{$\mathcal B$}}

\def\cH{\text{$\mathcal H$}}

\def\cL{\text{$\mathcal L$}}

\def\cN{\text{$\mathcal N$}}

\def\cP{\text{$\mathcal P$}}

\def\cR{\text{$\mathcal R$}}
\def\cS{\text{$\mathcal S$}}

\def\cV{\text{$\mathcal V$}}

\def\cX{\text{$\mathcal X$}}
\def\cY{\text{$\mathcal Y$}}
\def\cZ{\text{$\mathcal Z$}}

\newcommand{\id}{\operatorname{id}}
\newcommand{\cl}{\operatorname{cl}}

\newcommand{\Cl}{\operatorname{cl}}
\newcommand{\Int}{\operatorname{int}}
\newcommand{\inte}{\operatorname{int}}

\newcommand{\bd}{\operatorname{bd}}
\newcommand{\Bd}{\operatorname{Bd}}

\newcommand{\dom}{\operatorname{dom}}

\newcommand{\im}{\operatorname{im}}

\renewcommand{\emptyset}{\varnothing}

\newcommand{\Inv}{\operatorname{Inv}}

\def\proof{{\bf Proof:\ }}

\def\mathobj#1{\mbox{$#1$}}

\def\NN{\mathobj{\mathbb{N}}}

\def\QQ{\mathobj{\mathbb{Q}}}
\def\RR{\mathobj{\mathbb{R}}}

\def\TT{\mathobj{\mathbb{T}}}

\def\ZZ{\mathobj{\mathbb{Z}}}

\def\implies{\;\Rightarrow\;}
\def\iff{\;\Leftrightarrow\;}

\def\setof#1{\mbox{$\{\,#1\,\}$}}

\let\visiblecomments y

\newcommand{\cXtop}{\cX_{\rm top}}
\newcommand{\cAtop}{\cA^{\rm top}}
\newcommand{\cPb}{\overline{\cP}}
\newcommand{\Fr}{\operatorname{Fr}}
\newcommand{\Star}{\operatorname{Star}}
\def\opc#1{\protect\mbox{$\stackrel{\circ}{#1}$}}
\def\cell#1{\sigma_{#1}}
\newcommand{\Cell}{\Sigma}
\def\cebd{\operatorname{\partial}}
\renewcommand{\Int}{\operatorname{Int}}
\renewcommand{\Cl}{\operatorname{Cl}}
\newcommand{\mo}{\operatorname{mo}}
\newcommand{\Mo}{\operatorname{Mo}}
\newcommand{\Opn}{\operatorname{Opn}}
\newcommand{\ift}{\operatorname{ift}}
\newcommand{\ipt}{\operatorname{ipt}}
\newcommand{\Sol}{\operatorname{Sol}}
\newcommand{\kmax}{k_{\max}}

\DeclareMathOperator{\clos}{cl}

\usepackage{ifthen}
\usepackage{xcolor}
\newboolean{showNotes}

\setboolean{showNotes}{true}

\newcommand{\noteM}[1]{\ifthenelse{\boolean{showNotes}}{\textbf{\textcolor{red}{M: #1}}}{}}
\newcommand{\noteR}[1]{\ifthenelse{\boolean{showNotes}}{\textbf{\textcolor{green}{R: #1}}}{}}
\newcommand{\noteT}[1]{\ifthenelse{\boolean{showNotes}}{\textbf{\textcolor{blue}{T: #1}}}{}}

\newcommand{\DS}{\displaystyle}
\renewcommand{\phi}{\varphi}

\articletheorems

\begin{document}

\author{Marian Mrozek}
\address{Marian Mrozek, Division of Computational Mathematics,
  Faculty of Mathematics and Computer Science,
  Jagiellonian University, ul.~St. \L{}ojasiewicza 6, 30-348~Krak\'ow, Poland
}
\email{Marian.Mrozek@uj.edu.pl}
\author{Roman Srzednicki}
\address{Roman Srzednicki, Institute of Mathematics,
  Faculty of Mathematics and Computer Science,
  Jagiellonian University, ul.~St. \L{}ojasiewicza 6, 30-348~Krak\'ow, Poland
}
\email{Roman.Srzednicki@im.uj.edu.pl}
\author{Justin Thorpe}
\address{Justin Thorpe, Department of Mathematical Sciences,
George Mason University, Fairfax, Virginia 22030, USA
}
\email{jthorpe3@masonlive.gmu.edu}
\author{Thomas Wanner}
\address{Thomas Wanner, Department of Mathematical Sciences,
George Mason University, Fairfax, Virginia 22030, USA
}
\email{twanner@gmu.edu}

\thanks{
  Research of  M.M.\ was partially supported by
  the Polish National Science Center under Ma\-estro Grant No.\ 2014/14/A/ST1/00453
  and Opus Grant No.\ 2019/35/B/ST1/00874.
  Research of  R.S.\ was partially supported by
  the Polish National Science Center under Ma\-estro Grant No.\ 2014/14/A/ST1/00453.
  T.W.\ was partially supported by NSF grant DMS-1407087 and by the Simons Foundation
  under Award~581334.
}
\subjclass[2010]{Primary: 37B30; Secondary: 37E15, 57M99, 57Q05, 57Q15.}
 \keywords{Combinatorial vector field, isolated invariant set, Conley index, periodic orbit.}

\title%[Combinatorial vs. classical dynamics: recurrence]
{Combinatorial vs. classical dynamics: recurrence}

\date{\today}

\begin{abstract}
Establishing the existence of periodic orbits is one of the crucial
and most intricate topics in the study of dynamical systems, and 
over the years, many methods have been developed to this end. On 
the other hand, finding closed orbits in discrete contexts, such as
graph theory or in the recently developed field of combinatorial
dynamics, is straightforward and computationally feasible. In this
paper, we present an approach to study classical dynamical systems
as given by semiflows or flows using techniques from combinatorial
topological dynamics. More precisely, we present a general existence
theorem for periodic orbits of semiflows which is based on suitable
phase space decompositions, and indicate how combinatorial techniques
can be used to satisfy the necessary assumptions. In this way, one can
obtain computer-assisted proofs for the existence of periodic orbits
and even certain chaotic behavior.
\end{abstract}

\maketitle

\tableofcontents

%%%%%%%%%%%%%%%%%%%%%%%%%%%%%%%%%%%%%%%%%%%%%%%%%%%%%%%%%%%%%%%%%
%%%%%%%%%%%%%%%%%%%%%%%%%%%%%%%%%%%%%%%%%%%%%%%%%%%%%%%%%%%%%%%%%
%%%%%%%%%%%%%%%%%%%%%%%%%%%%%%%%%%%%%%%%%%%%%%%%%%%%%%%%%%%%%%%%%
\section{Introduction}

Understanding the asymptotic dynamics of solutions of differential equations is vital in many areas of science.
Since in most cases no closed formulas for solutions are available, various other approaches have been invented
to get some insight into the behavior of solutions. Qualitative methods based on tools from analysis or topology
provide invariants and theorems detecting certain types of solutions. Unfortunately, analytic verification of the assumptions
of these theorems in many cases is very arduous or just infeasible. 
Computer simulations give quick answers, but are not very reliable because of the approximate nature of numerical methods.  
Computer-assisted proofs in dynamics emerged at the end of the 20th century as a remedy to this situation. 
They combine qualitative tools with rigorous numerics as a method replacing the arduous analytic computations. 

A typical computer-assisted proof requires an a-priori selection of the qualitative method for the particular problem and fine tuning 
of several parameters, for instance the step of the numerical method,  bounds for the area of interest, or the location of a Poincar\'e section. 
To properly choose the method and get the parameters right, many numerical experiments have to be performed first. 
This effort pays for itself in the outcome which often is not only the proof of the existence of the object of interest, 
but also sharp bounds for its location in space. 

Another approach to computer-assisted proofs in dynamics aims at building a global combinatorial model with some formal ties
between the model dynamics and the dynamics of interest.  This approach turned out to be very successful in the study of the gradient
structure of a dynamical system \cite{AraiEtal:2009} by providing a Morse decomposition of the system together with the associated Conley-Morse graph. 
In this approach, individual Morse sets in the decomposition are identified on the combinatorial level as strongly connected components
in the directed graph representation of the combinatorial dynamics. 
What we propose in this paper is an extension of the combinatorial approach to get insight into the Morse sets  in the case when a Morse set
exhibits non-trivial recurrent behavior. 

The extension combines a theorem characterizing the existence of periodic orbits for semiflows
in terms of Conley index \cite{mccord:etal:95a} with the recently developed Conley theory for combinatorial multivector fields \cite{lipinski:etal:p20a,mrozek:17a}.
The combinatorial multivector field is obtained by constructing a polygonal decomposition of space whose faces are crossed by the flow transversally, a method already
applied in \cite{boczko:etal:07a,eidenschink:95a}. By lifting the problem to combinatorial multivector fields we can not only identify the Morse sets, but also look 
inside the Morse sets searching for recurrent isolated invariant sets. Once an isolated invariant set with a combinatorial Poincar\'e section and the Conley index of a hyperbolic periodic orbit 
is located, we can use the theory presented in this paper to prove the existence of a periodic orbit in the original flow. 

The results of the present paper apply to flows and the existence of a periodic orbit.
But, they may be extended to semiflows and the existence of chaotic invariant sets.
Whenever possible, we have formulated our results for the most general situation.
In this way, we will be able to present some applications to semiflows and chaos already in
Section~\ref{sec:app-main-results}. However, the full details of the general extension
to semiflows, even though they are the same in spirit, are technically more complicated
and will be presented elsewhere.

The remainder of this paper is organized as follows. In Section~\ref{sec:app-main-results} we give a broad
overview over potential applications for our main result. For this, we demonstrate how one can establish 
the existence of periodic orbits in planar systems, show that for a semiflow on a branched manifold
it can be used to establish chaos, and then present a template for an analogous chaos-inducing 
situation in three dimensions. After we present necessary background definitions and results in
Section~\ref{sec:prelims}, the following Section~\ref{sec:periodicex} is devoted to our main result.
Based on the notion of brick decompositions, we can present a general existence result for periodic
orbits which is based on Conley index methods. In contrast to earlier results, we do not need the
existence of an explicit Poincar\'e section. Rather, our result relies on a much weaker notion of
flow directionality. Finally, Section~\ref{sec:ctdrecurrence} shows how methods from combinatorial
topological dynamics based on multivector fields can be used to establish the assumptions of our
main result.

%%%%%%%%%%%%%%%%%%%%%%%%%%%%%%%%%%%%%%%%%%%%%%%%%%%%%%%%%%%%%%%%%
%%%%%%%%%%%%%%%%%%%%%%%%%%%%%%%%%%%%%%%%%%%%%%%%%%%%%%%%%%%%%%%%%
%%%%%%%%%%%%%%%%%%%%%%%%%%%%%%%%%%%%%%%%%%%%%%%%%%%%%%%%%%%%%%%%%
\section{Applications of the main results}
\label{sec:app-main-results}

Before presenting the mathematical theory behind our approach to
recurrence, the current section is devoted to illustrating a few applications.
This is accomplished in three parts. We begin by demonstrating
the basic approach in a simple setting, namely flows generated by
ordinary differential equations in the plane. This is followed by
an example taken from the sequence of papers~\cite{batko:etal:20a,
kaczynski:etal:16a, mrozek:wanner:21a}, where we can use our results
to establish the existence of infinitely many different periodic 
orbits on a branched manifold and with respect to a large class of
semiflows. Finally, the third example is concerned with a Lorenz-type
situation in three dimensions and which also exhibits chaos.

%%%%%%%%%%%%%%%%%%%%%%%%%%%%%%%%%%%%%%%%%%%%%%%%%%%%%%%%%%%%%%%%%
%%%%%%%%%%%%%%%%%%%%%%%%%%%%%%%%%%%%%%%%%%%%%%%%%%%%%%%%%%%%%%%%%
\subsection{Periodic orbits in planar flows}

As mentioned in the introduction, our main result for proving the
existence of periodic orbits is based on a suitable decomposition
of phase space. In this section, we illustrate how this can be used
to study periodic orbits in planar flows generated by ordinary
differential equations. We would like to emphasize that our
examples were primarily chosen to illustrate the main ideas. In
these simple cases, the existence of the periodic solutions can
easily be obtained by different means. Yet, the approach outlined
below can also be used in a higher-dimensional setting.

Our main result, Theorem~\ref{thm:bricks} below, relies on a decomposition
of a region of interest in phase space into small pieces called bricks. In
the simplest case, this decomposition can take the form of a triangulation,
even though this will be generalized later on. Our approach for establishing
the existence of periodic orbits can then be summarized in the following
five steps:
\begin{itemize}
\item[(a)] {\bf Numerical determination of a flow transverse triangulation:}
After identifying a region of interest in phase space, use local flow or
vector field information to construct a triangulation of the region in the
following way. Across each edge of the triangulation, the considered flow
is transverse, and there are no equilibrium solutions of the flow in any
of the triangles.
\item[(b)] {\bf Rigorous validation of the triangulation properties:}
Once a numerical candidate for the triangulation has been found, rigorous
computations based on interval arithmetic have to be used to validate the
flow behavior across the triangulation edges and inside the triangles.
\item[(c)] {\bf Create a combinatorial multivector field:}
Based on the rigorous information from the previous step, one can
construct a combinatorial multivector field associated with the flow
as described in Proposition~\ref{prop:cmvf}.
\item[(d)] {\bf Graph-theoretic analysis of the multivector field:}
Finding recurrent sets in the combinatorial setting is now a
graph-theoretic problem, which can be solved efficiently using the
concept of strongly connected components in a digraph. If the triangulation
constructed in~(1) was appropriate, then we can obtain a candidate
region which can contain a periodic orbit.
\item[(e)] {\bf Existence proof for periodic solutions:}
As a last step, the usually large triangulation has to be coarsened
in such a way that Theorem~\ref{thm:bricks} can be applied to prove
the existence of the periodic orbit.
\end{itemize}
The crucial part of this procedure is clearly the first step, since
we need to find an ``appropriate triangulation'' to apply our results.
Nevertheless, a number of options exist, for example the methods described
in~\cite{boczko:etal:07a, eidenschink:95a}. For our examples below we
will follow a different approach. We would like to point out, however,
that the use of a triangulation is not necessary. From a mathematical
point of view, the bricks which form the phase space decomposition
are characterized by having well-defined entry and exit behavior
with respect to the underlying semiflow or flow, and they have to be
topologically simple. This will be discussed in more detail in
Section~\ref{sec:ctdrecurrence}.

To illustrate the above approach we consider two simple planar ordinary
differential equations. The first one is given by
\begin{equation} \label{eqn:excircles}
  \begin{array}{rcl}
    \DS \dot{x} & = & \DS -y + x \left( x^2 + y^2 - 4 \right)
      \left( x^2 + y^2 - 1 \right) \\[1ex]
    \DS \dot{y} & = & \DS \;\;\, x + y \left( x^2 + y^2 - 4 \right)
      \left( x^2 + y^2 - 1 \right)
  \end{array}
\end{equation}
and one can easily see that this system has an unstable equilibrium at 
the origin, as well as two invariant circles which are centered at~$(0,0)$
and have radii~$1$ and~$2$, respectively. For our second example, we
consider the celebrated van der Pol equation
\begin{equation} \label{eqn:exvdpol}
  \begin{array}{rcl}
    \DS \dot{x} & = & \DS y \\[0.5ex]
    \DS \dot{y} & = & \DS \mu y \left( 1 - x^2 \right) - x
  \end{array}
\end{equation}
at the parameter value $\mu = 1$.

\begin{figure}[tb]
  \begin{center}
  \includegraphics[width=0.48\textwidth]{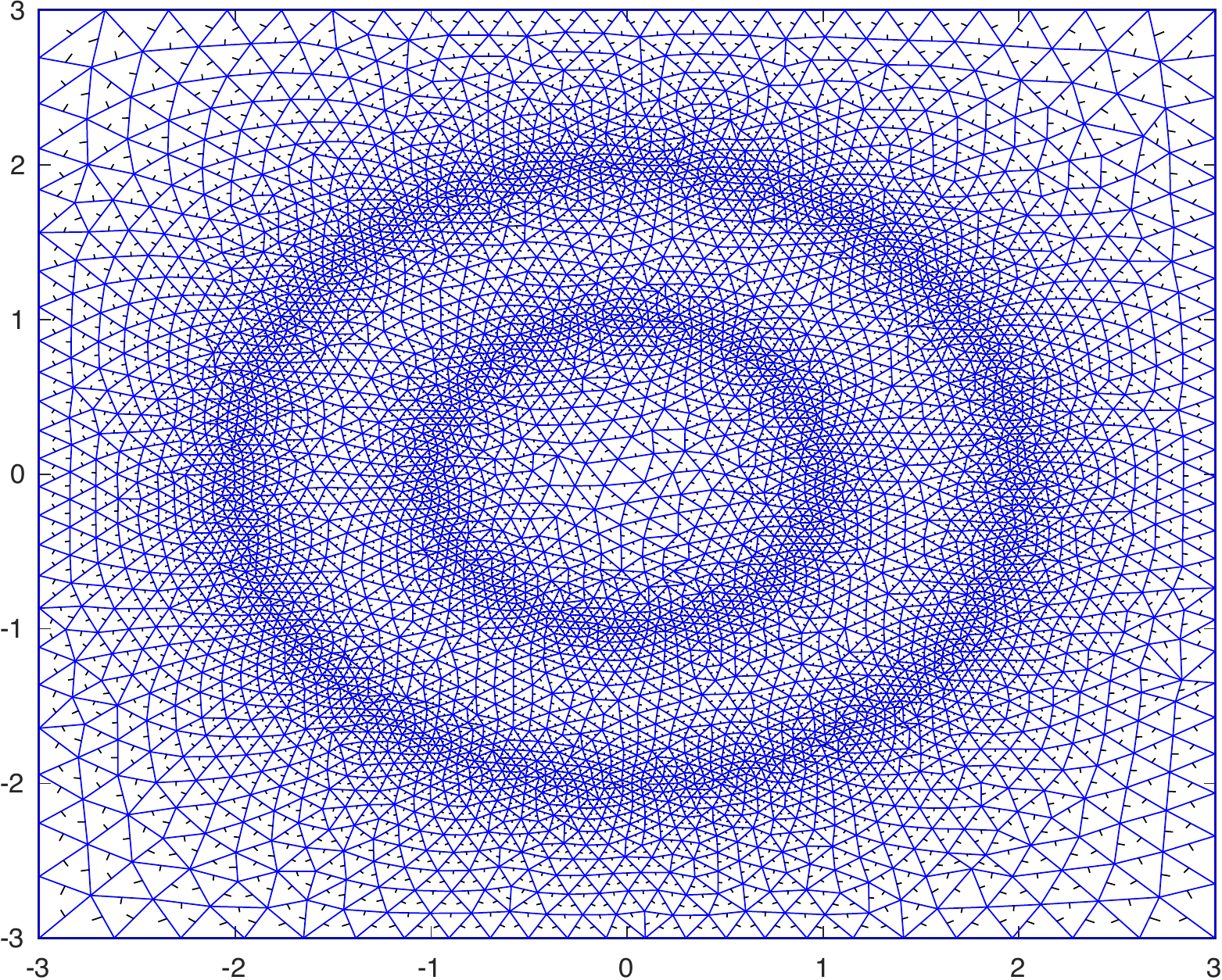}
  \hspace*{0.3cm}
  \includegraphics[width=0.48\textwidth]{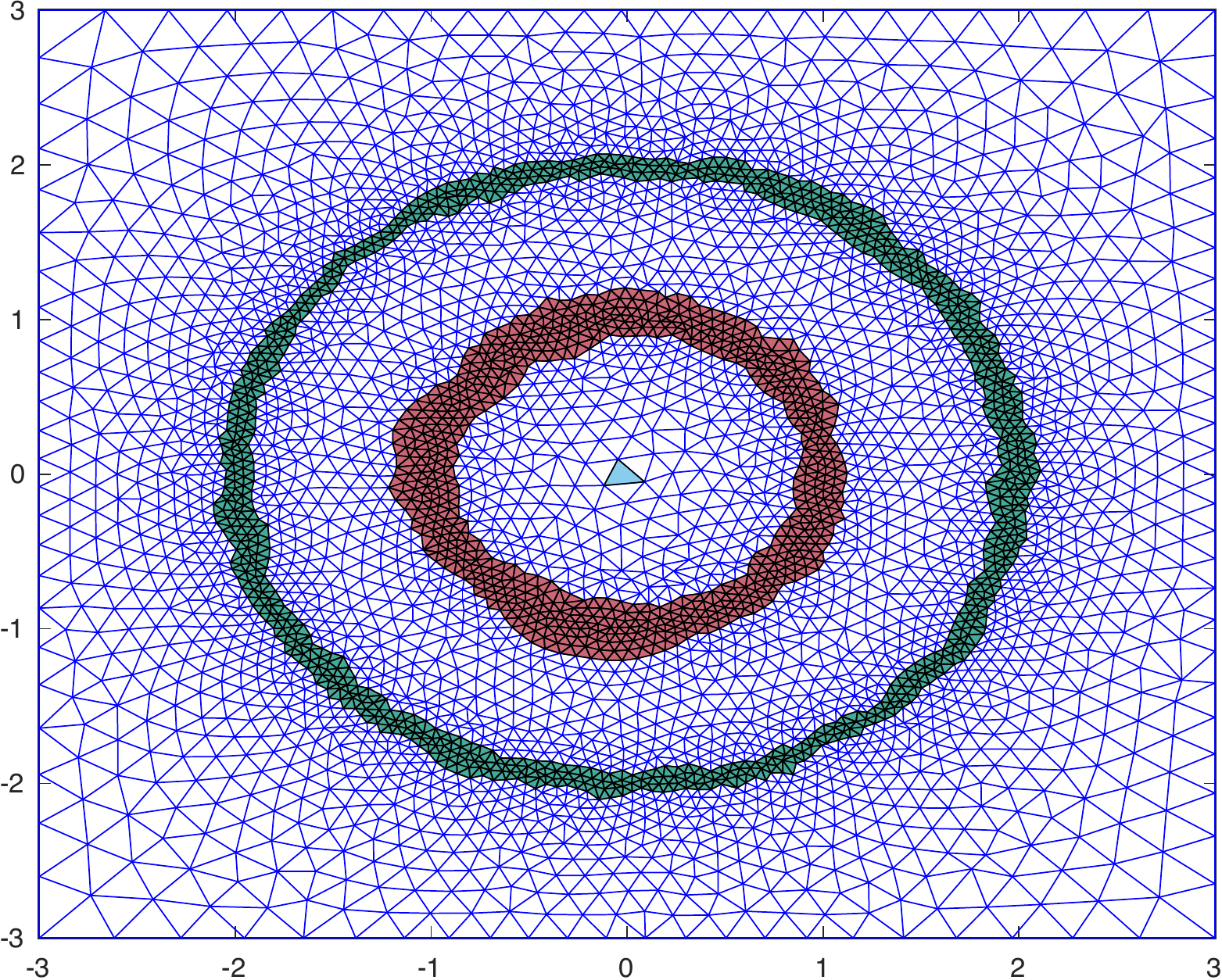}
  \caption{Flow transverse triangular grid for the flow described
           by~(\ref{eqn:excircles}). In the left image, small black
           line segments indicate in which direction the flow 
           crosses the respective edge, with the line segment pointing
           into the entered triangle. The right image shows three 
           isolated invariant sets for the combinatorial multivector 
           field generated by the phase space decomposition. While 
           the isolated triangle in the middle contains an equilibrium,
           the two ring-shaped regions do not. Those latter regions are 
           the ones that can be used for our approach.}
  \label{fig:excircles1}
  \end{center}
\end{figure}

For the first example~(\ref{eqn:excircles}), it is straightforward to create
a subdivision of the square~$[-3,3]$ into a triangular mesh such that
along each edge of the mesh the flow generated by~(\ref{eqn:excircles})
is transverse. While the details of this construction will be described
elsewhere, it is based on an adaptive initial meshing procedure, combined
with random perturbations of the vertices of the mesh to achieve flow
transversality. The resulting mesh is shown in the left diagram of
Figure~\ref{fig:excircles1}, where the flow directions across
the mesh edges are indicated by small black line segments which
point into the entered mesh triangle. Moreover, we used the interval
arithmetic package {\sc Intlab}~\cite{rump:99a} to rigorously verify
the flow transversality.

\begin{figure}[tb]
  \begin{center}
  \includegraphics[width=0.48\textwidth]{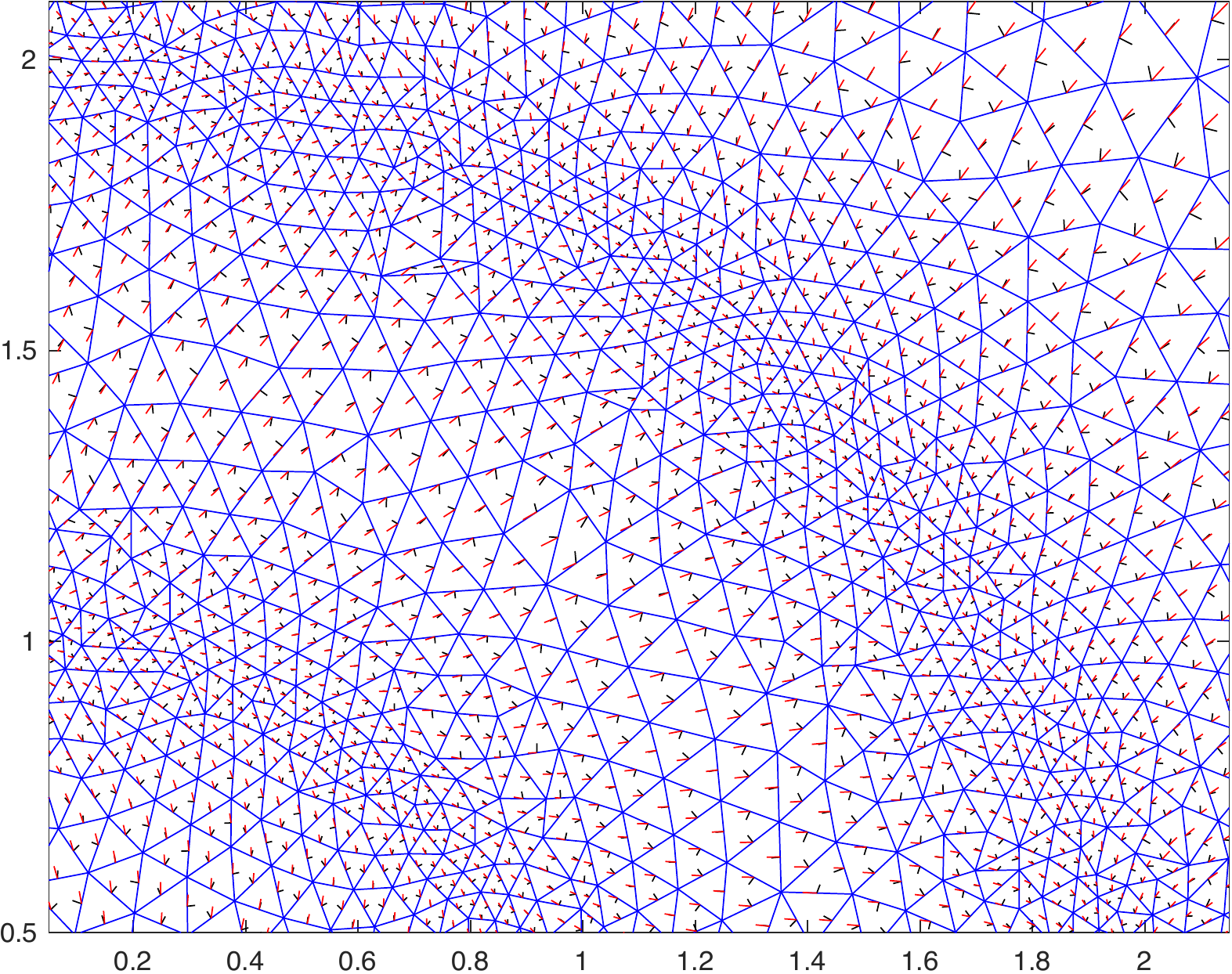}
  \hspace*{0.3cm}
  \includegraphics[width=0.48\textwidth]{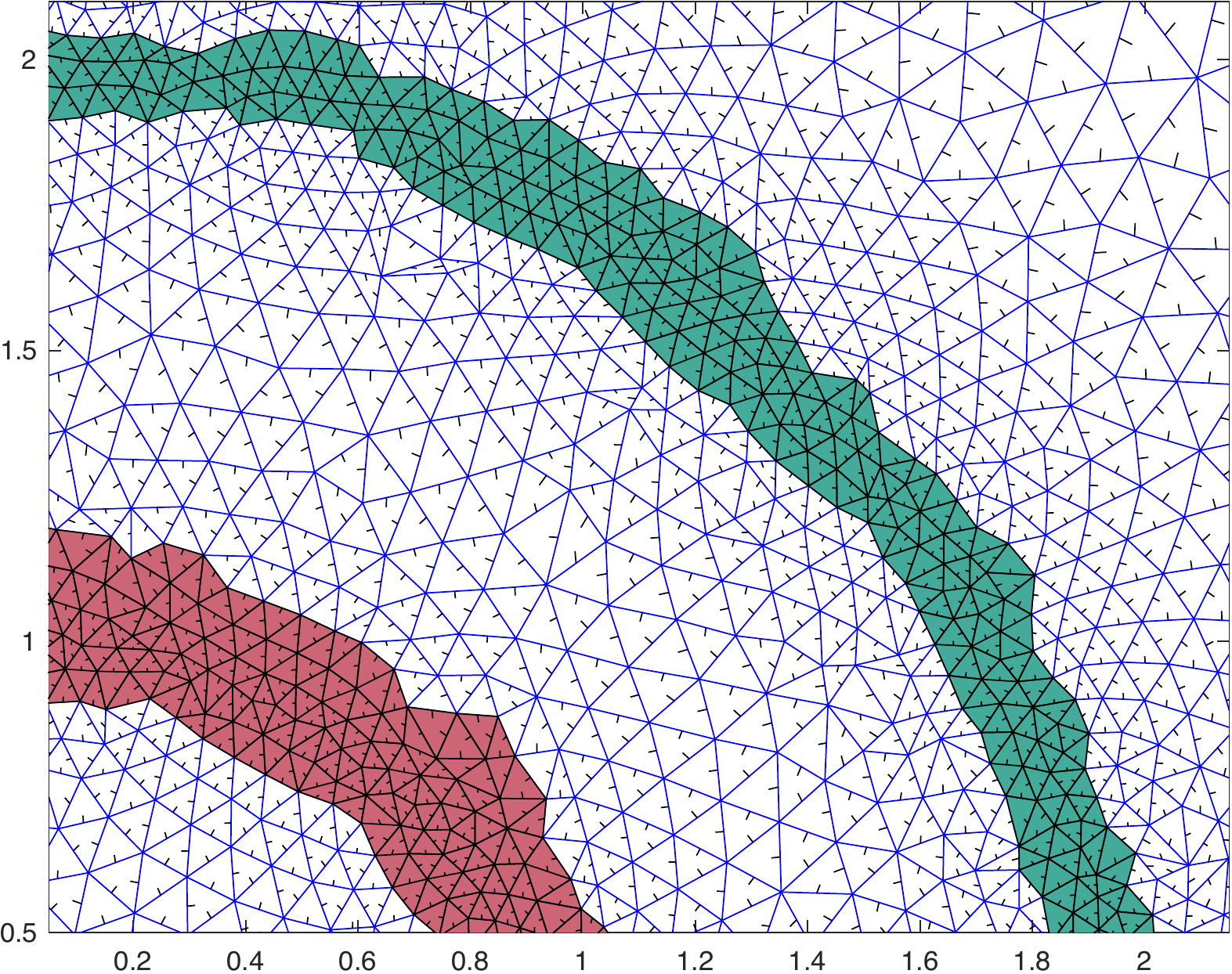}
  \caption{Closeups of the flow transverse grid in
           Figure~\ref{fig:excircles1}, for the flow described
           by~(\ref{eqn:excircles}). While the right image shows
           parts of the combinatorial isolated invariant sets 
           containing the periodic orbits, it also shows the transversality
           directions via small black line segments as in
           Figure~\ref{fig:excircles1}. In the left image, we 
           indicate the actual flow direction at the center points of
           the mesh edges via small red line segments.}
  \label{fig:excircles2}
  \end{center}
\end{figure}

Based on the mesh and the flow directions, one can then create a 
multivector field on the simplicial complex given by the mesh which
encodes the flow information in a combinatorial way. By employing the
results of~\cite{lipinski:etal:p20a}, one can then find isolated
invariant sets for the combinatorial multivector field. These are
shown in the right diagram of Figure~\ref{fig:excircles1}. In addition
to the triangle containing the equilibrium~$(0,0)$, one can also see two
regions enclosing the two periodic orbits. Closeups of these two diagrams
are contained in Figure~\ref{fig:excircles2}. Finally, one can divide
each of the two annular regions into suitable segments in such a way
that an application of Theorem~\ref{thm:bricks} establishes a periodic
orbit in each region.

\begin{figure}[tb]
  \begin{center}
  \includegraphics[width=0.48\textwidth]{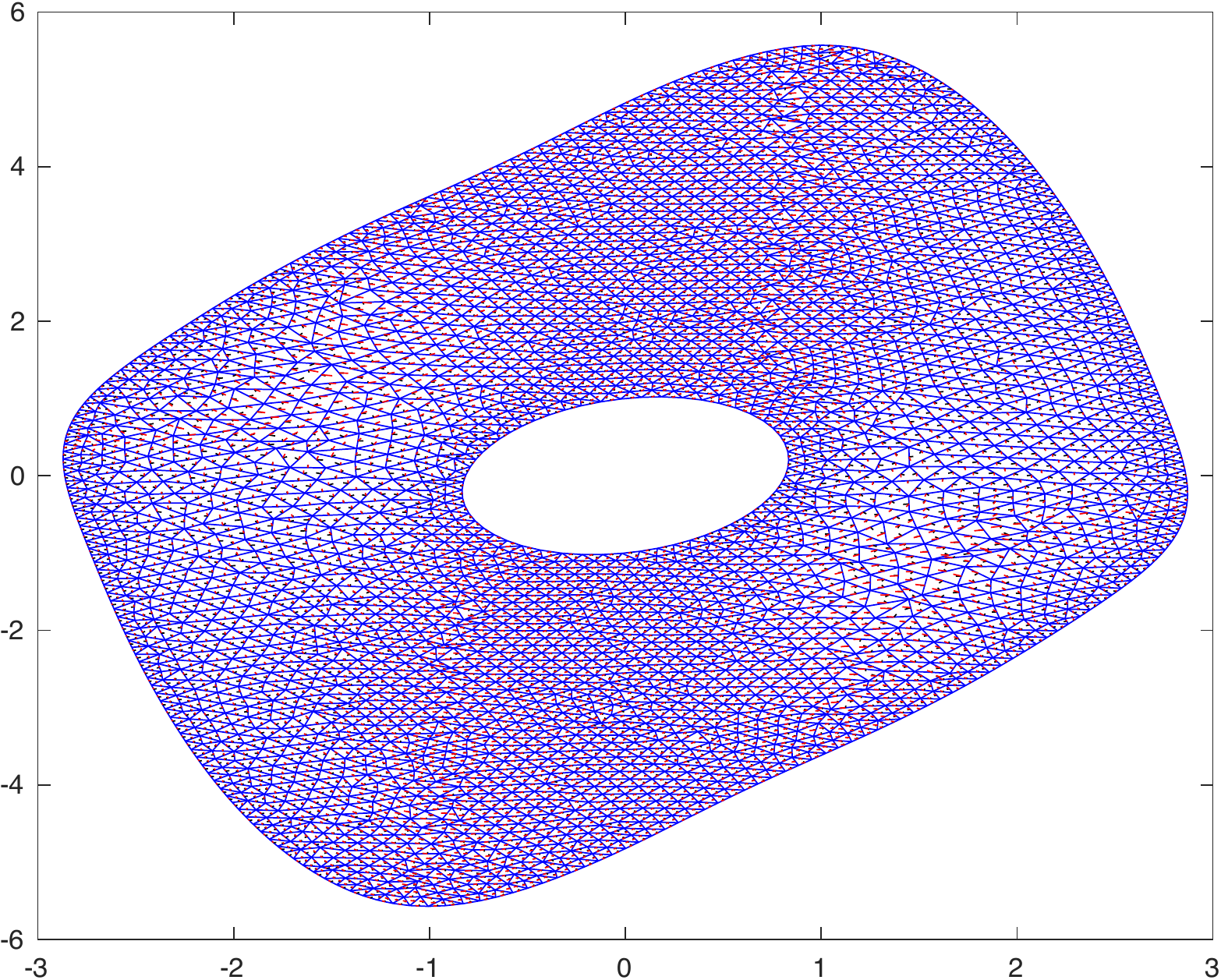}
  \hspace*{0.3cm}
  \includegraphics[width=0.48\textwidth]{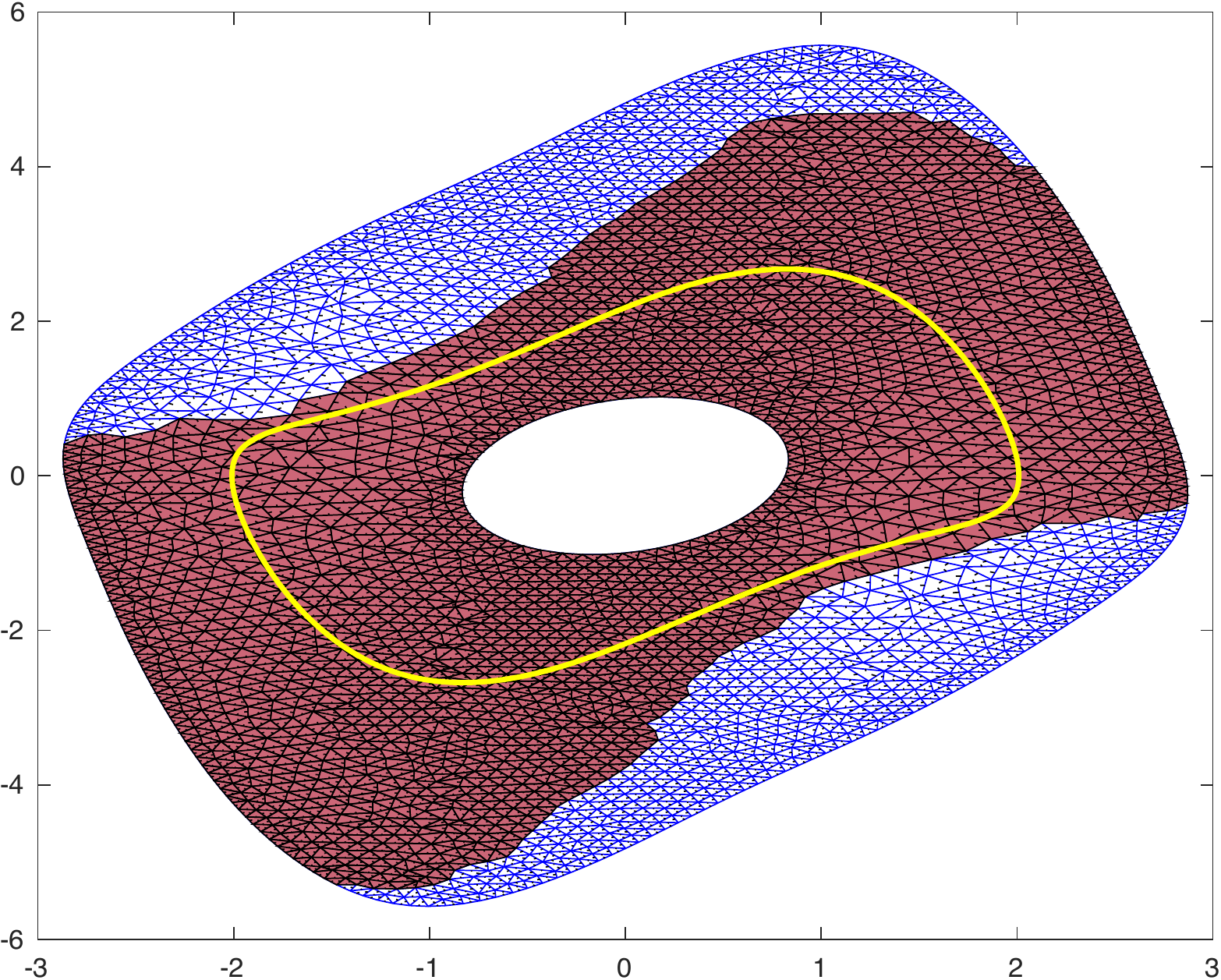}
  \caption{Flow transverse triangular grid for the flow described
           by the van der Pol equation~(\ref{eqn:exvdpol}) for $\mu = 1$.
           In the left image, small black line segments indicate in which
           direction the flow crosses the respective edge, and the actual
           flow direction at the center points of the mesh edges is given
           via small red line segments. The right image shows an isolated
           invariant set which encloses the limit cycle shown in
           yellow.}
  \label{fig:exvdpol1}
  \end{center}
\end{figure}

For the van der Pol example~\eqref{eqn:exvdpol} the first step in the
above procedure is somewhat more delicate, due to the strong rotational
dynamics exhibited by the system. Nevertheless, in Figure~\ref{fig:exvdpol1}
we show a triangulation of an attracting region, which has been verified
via interval arithmetic as being flow transverse across triangulation edges.
Within this region, the shaded red region is an isolating neighborhood for
the periodic orbit, which is indicated in yellow. This isolating
neighborhood was found through the associated combinatorial multivector
field as before, and again we can apply Theorem~\ref{thm:bricks} after
a segmentation of this region.

Our above discussion centered merely on the existence of periodic 
orbits in planar systems. It is natural to wonder whether one can
also establish the existence of connections between such orbits, 
or between periodic orbits and stationary states, in order
to obtain a more global picture of the underlying dynamics. This can
actually be accomplished using the notion of connection matrices,
which in the case of combinatorial dynamics on Lefschetz spaces were
recently introduced in~\cite{mrozek:wanner:p21a}. The detail of this
approach will be presented elsewhere. Moreover, while our approach
in its present form only provides the existence of a periodic orbit,
our method can be used in principle to also obtain lower bounds on
the period of the constructed periodic solutions. For this, one can
use vector field information to obtain lower bounds on the traversal
times of solutions across every cell in the phase space decomposition,
which can then be combined with combinatorial periodic solutions
in the digraph provided by our multivector field to establish the
bound. This will not be pursued further in this paper.

%%%%%%%%%%%%%%%%%%%%%%%%%%%%%%%%%%%%%%%%%%%%%%%%%%%%%%%%%%%%%%%%%
%%%%%%%%%%%%%%%%%%%%%%%%%%%%%%%%%%%%%%%%%%%%%%%%%%%%%%%%%%%%%%%%%
\subsection{Infinitely many periodic orbits on a branched manifold}

\begin{figure}[tb]
  \begin{center}
  \includegraphics[width=0.47\textwidth]{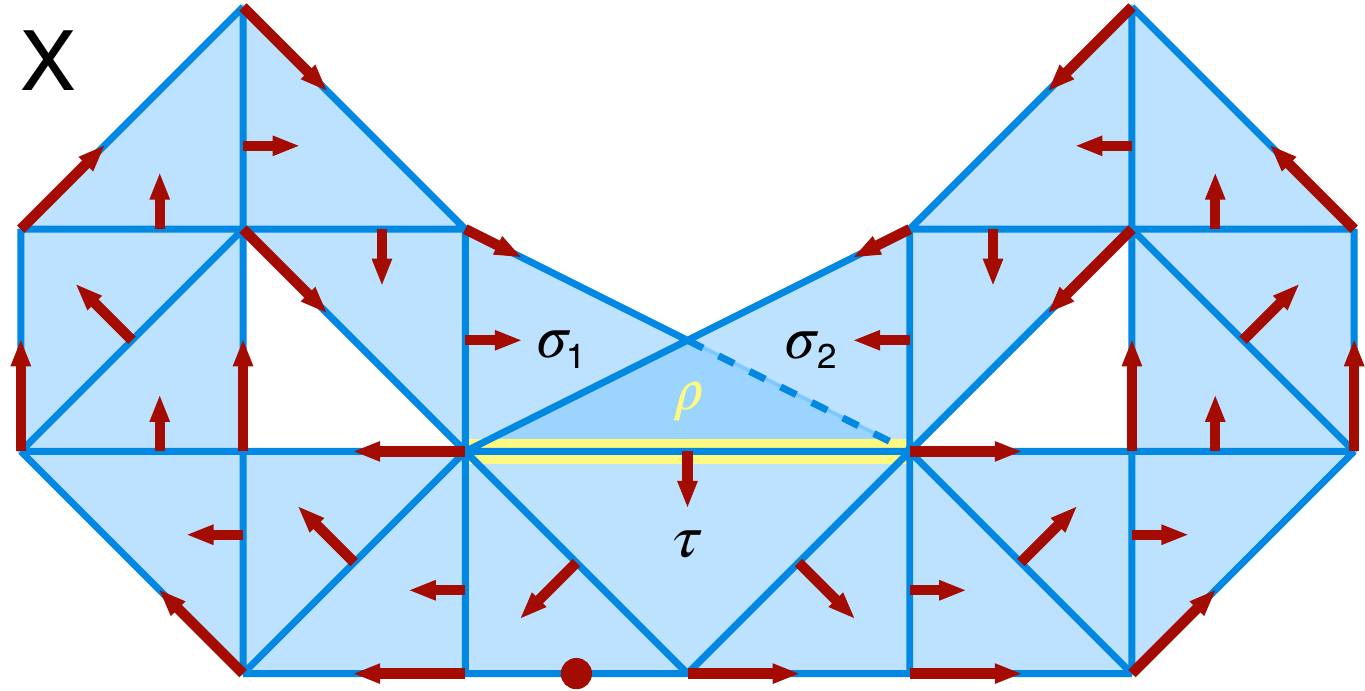}
  \hspace*{0.5cm}
  \includegraphics[width=0.47\textwidth]{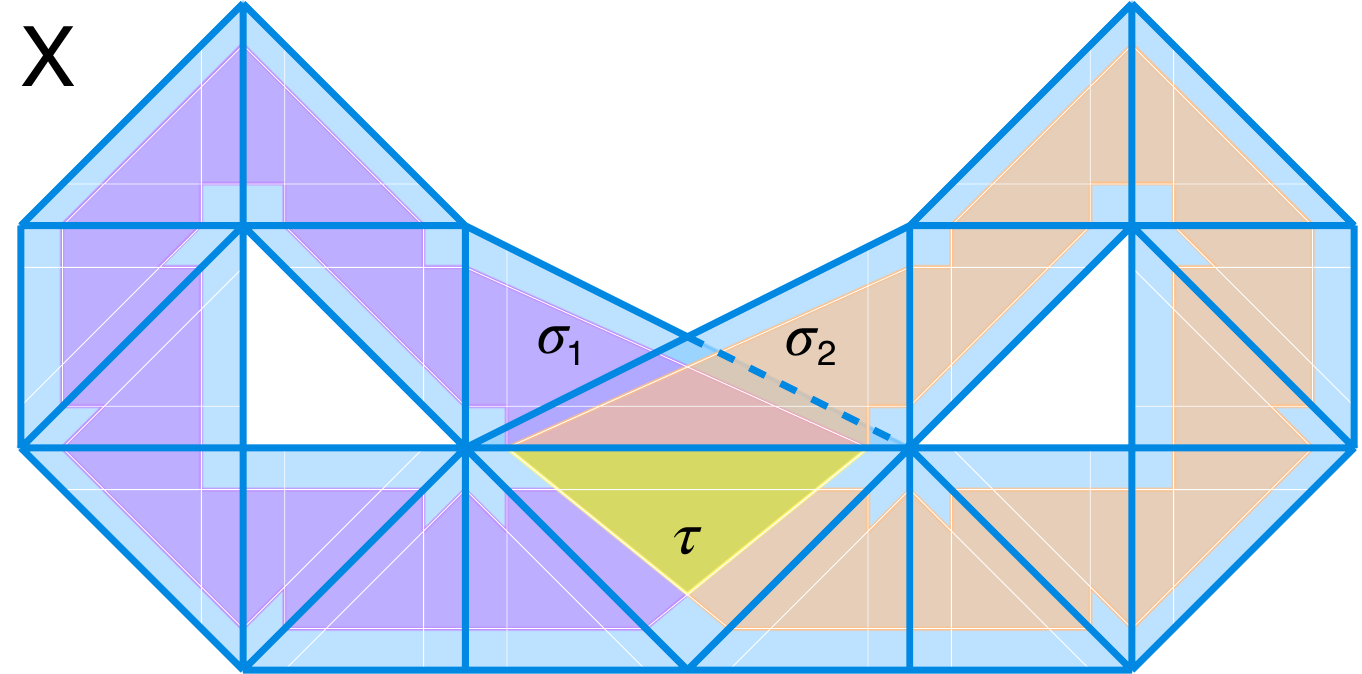}
  \caption{Combinatorial Lorenz model. The left panel shows a
           a combinatorial vector field on a branched manifold~$X$
           which can easily be seen as having chaotic dynamics.
           The union of the regions shown in violet, orange, and
           yellow in the right panel gives an isolating neighborhood
           for the semiflow constructed in~\cite{mrozek:wanner:21a},
           which has the same Conley-Morse graph as the combinatorial
           model. This region corresponds to the Morse set in the
           combinatorial setting which contains all triangles,
           and the associated semiflow exits everywhere along
           its boundary. Notice that while the image on the left
           indicates an equilibrium, we apply our approach for finding
           periodic solutions only to the union of the regions shown
           in violet, orange, and yellow in the right panel, and they
           will not contain any equilibrium solution.}
  \label{fig:clorenz}
  \end{center}
\end{figure}

For our second application we consider recent work on the relationship
between combinatorial vector fields in the sense of Forman and classical
dynamics, as developed in~\cite{batko:etal:20a, kaczynski:etal:16a,
mrozek:wanner:21a}. In these papers, we developed a Conley theory for
combinatorial vector fields, and showed that every combinatorial vector
field on a simplicial complex gives rise to both a discrete multivalued
dynamical system and a semiflow on the underlying polytope which exhibits
the same dynamics in the sense of Conley-Morse graphs. It was also 
shown that in the combinatorial setting one can easily construct 
systems with complex combinatorial dynamics. For example, it was
mentioned in~\cite{kaczynski:etal:16a} that the combinatorial vector
field shown in the left panel of Figure~\ref{fig:clorenz} exhibits
chaotic behavior in the sense of having infinitely many different
periodic orbits, which can easily be found by following the triangle
loops around each of the two eyes in a variety of sequence patterns.
In this model, the two triangles~$\sigma_1$ and~$\sigma_2$ only
intersect along the edge~$\rho$, and along this edge they are
attached to the triangle~$\tau$. Thus, the underlying space is a 
simplicial complex in the form of a branched two-dimensional 
manifold~$X$.

While the combinatorial vector field in the left panel of
Figure~\ref{fig:clorenz} exhibits infinitely many periodic orbits,
this is not clear when one only considers the associated Conley-Morse
graph. It can easily be shown that the combinatorial dynamics has 
four distinct Morse sets: The union of all triangles is contained in
one large Morse set, below which there is a Morse set of index one
given by the equilibrium on the bottom edge. Below this equilibrium,
there are two stable periodic orbits, which correspond to the
interior boundaries of the two eyes.

In~\cite{mrozek:wanner:21a} we showed that for any combinatorial
vector field, there always exists a semiflow on the underlying polytope
which has the same Conley-Morse graph. Thus, for the particular
vector field in Figure~\ref{fig:clorenz} there exists a semiflow 
with the four Morse sets discussed above. While we do not want to
go into the details of the construction in~\cite{mrozek:wanner:21a},
it provides an isolating neighborhood for the top-level Morse
set which is indicated in the right panel of Figure~\ref{fig:clorenz}
through additional shading in violet, orange, and yellow. In fact, the
semiflow provided by this construction has the whole boundary of the
shaded region as exit set.

But what about the dynamics of this semiflow inside this isolating
neighborhood? It is far from clear whether there still are infinitely
many periodic orbits in the colored region. Nevertheless, one can prove
exactly that. The following theorem relies heavily on the setup
of~\cite{mrozek:wanner:21a}, and we will only sketch the proof of
this result below.
\begin{thm}[Periodic solutions with prescribed itineraries]
\label{thm:clorenz}
Consider the combinatorial vector field shown in the left panel
of Figure~\ref{fig:clorenz}, and let~$\phi$ denote any strongly
admissible semiflow on the underlying branched manifold~$X$, in
the sense described in~\cite{mrozek:wanner:21a}. Furthermore,
let~$L$ describe the subset of~$X$ colored in violet and yellow
in the right panel of Figure~\ref{fig:clorenz}, and~$R$
the one colored in orange and yellow. Furthermore, consider an
arbitrary finite sequence~$\cS$ using the binary alphabet~$\{ L, R \}$.
Then there exists a periodic orbit of~$\phi$ which infinitely
many times traverses the sets~$L$ and~$R$ in clockwise and
counter-clockwise manner, respectively, as prescribed by the
sequence~$\cS$.
\end{thm}
For example, for the sequence~$\cS = (R,L)$ which will be discussed
in the proof below, the periodic orbit alternately rotates around the
right and left loop of~$X$. We would like to point out, however, that
we cannot guarantee that during the minimal period of the periodic
solution the periodic orbit only completes one~$(R,L)$ sequence. It 
might complete several of these. Nevertheless, one easily obtains
the following corollary.
\begin{cor}[Infinitely many periodic orbits on a branched manifold]
\label{cor:clorenz}
In the setting of Theorem~\ref{thm:clorenz}, any strongly admissible
semiflow~$\phi$ on~$X$ contains infinitely many periodic orbits in the
union~$L \cup R$.
\end{cor}
As mentioned earlier, we only briefly sketch the proof of
Theorem~\ref{thm:clorenz} in a special case, and leave the details
of the general case to the reader.
\medskip
\begin{figure}[tb]
  \begin{center}
  \includegraphics[width=0.47\textwidth]{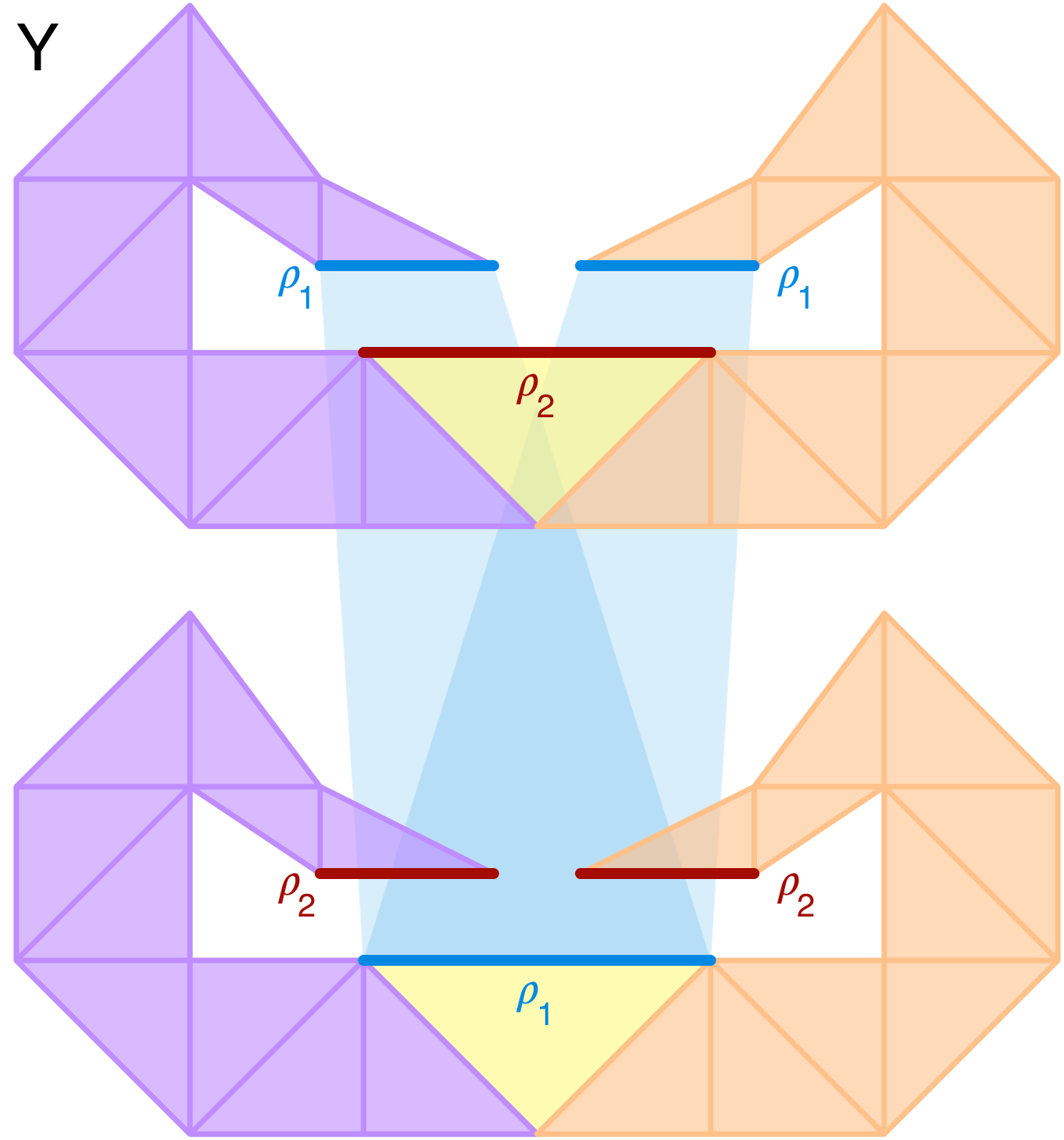}
  \hspace*{0.5cm}
  \includegraphics[width=0.47\textwidth]{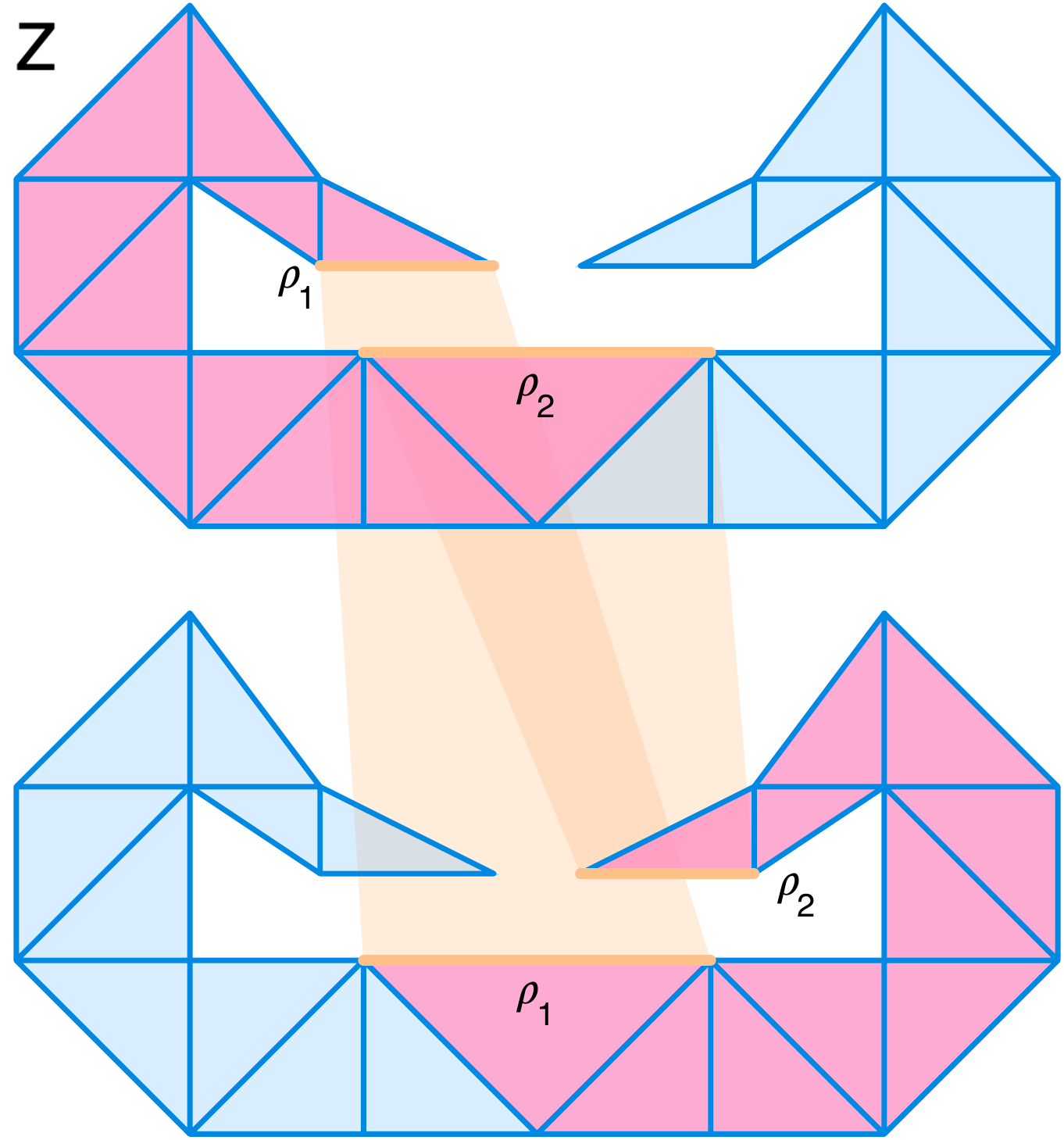}
  \caption{The spaces~$Y$ and~$Z$ used in the proof of
           Theorem~\ref{thm:clorenz}. The covering space~$Y$ is shown in
           the left panel, and it is a double cover of the space~$X$ from
           Figure~\ref{fig:clorenz}. It consists of two disjoint copies of~$X$,
           each of which has been cut along the edge~$\rho$. The three blue
           edges marked~$\rho_1$ are then identified, as indicated by the
           shaded blue regions. Similarly, the three red edges marked~$\rho_2$
           are identified, but this is not separately highlighted by additional
           shading. The space~$Z$ is shown in the right panel. It is obtained 
           from~$Y$ by cutting again along certain copies of the edge~$\rho$,
           but this time based on the underlying symbol sequence~$\cS = (R,L)$.
           At the lower copy of~$X$, we cut off the connection between the
           left loop and~$\rho_1$, while in the upper copy the right loop
           is severed from~$\rho_2$. The remaining edge identifications
           are shown in shaded orange. This leads to the pink set which
           encodes the dynamics prescribed by~$\cS$.}
  \label{fig:clorenzcover}
  \end{center}
\end{figure}
\proof
For the proof, we only consider the specific sequence~$\cS = (R,L)$.
The case of a general finite sequence~$\cS$ in~$\{ L,R \}$ can be 
discussed analogously.

In order to apply our main Theorem~\ref{thm:bricks} in this situation,
we cannot consider the space~$X$ directly. Rather, we consider a space~$Z$
which is based on a suitable covering space~$Y$. For our symbol sequence of
length two, the covering space~$Y$ will be a double cover, and it is
illustrated in the left panel of Figure~\ref{fig:clorenzcover}. The space
consists of two disjoint copies of~$X$, each of which has been cut along
the edge~$\rho$ shown in the left panel of Figure~\ref{fig:clorenz}. This
results in a total of six copies of this edge, three for each copy of~$X$,
and they are identified as shown in the figure. The three blue edges
marked~$\rho_1$ are identified, as indicated by the shaded blue regions.
Similarly, the three red edges marked~$\rho_2$ are identified, but in order
to keep the sketch simple, this is not separately highlighted by additional
shading.

Once the covering space~$Y$ has been created, we can construct the space~$Z$
which is actually used in the proof. The space~$Z$ is obtained  from~$Y$ by
cutting again along both copies of the edge~$\rho$, but this time based on the
underlying symbol sequence~$\cS = (R,L)$. At the lower copy of~$X$, we cut off
the connection between the left loop and~$\rho_1$, since the first symbol in~$\cS$
is given by~$R$. Similarly, in the upper copy the right loop is severed
from~$\rho_2$, since the second symbol in~$\cS$ is~$L$. The remaining edge
identifications are shown in shaded orange. This leads to the pink set which
now encodes the dynamics prescribed by~$\cS$.

After these preparations, we finally consider the restriction~$B \subset Z$ of
the isolating neighborhood indicated in the right panel of Figure~\ref{fig:clorenz},
and lifted to~$Z$ in a straightforward way. Notice that this set~$B$ is a 
proper subset of the pink region in the right panel of Figure~\ref{fig:clorenzcover}.
Then it follows from the results of~\cite{mrozek:wanner:21a} that~$B$ is indeed
an isolating block for the semiflow~$\phi$, whose exit set~$B^-$ equals the
boundary of~$B$ in~$Z$. Moreover, the set~$B$ has an underlying brick decomposition
as required by Theorem~\ref{thm:bricks}. We would like to point out that the
cuts that were performed while passing from~$Y$ to~$Z$ are necessary to deal with
the issues identified for semiflows in Proposition~\ref{prop:topl-intersection2}
below. Finally, the Conley index of the isolating block~$B$ is that of an unstable
periodic orbit of index one. An application of Theorem~\ref{thm:bricks} now easily
implies the existence of the desired periodic orbit. For this, one only has to project
the periodic orbit in~$Z$ guaranteed by the theorem down to the space~$X$. Since the
space~$Z$ was constructed via a covering space~$Y$ which contains two exact copies
of both~$X$ and the semiflow on~$X$, this projection leads to a periodic orbit
for the original semiflow~$\phi$ on~$X$.
\qed\medskip

We would like to point out that the above two results remain valid for
any semiflow on~$X$ which is strongly admissible in the sense
of~\cite{mrozek:wanner:21a}. Thus, the trivial combinatorial chaos
exhibited by the Forman vector field from Figure~\ref{fig:clorenz} gives
rise to a robust notion of chaos for a large class of semiflows on the
underlying polytope. In other words, the main results of this paper, in
combination with purely combinatorial constructions, can be used to 
create classical semiflows with extremely complicated recurrent dynamics.

One drawback of these results is that we cannot guarantee that the
constructed periodic orbits realize the symbol sequence~$\cS$ within 
their fundamental period. This can, however, be accomplished via
an approach based on the Lefschetz fixed point theorem and combinatorial
Poincar\'e maps. This will be presented elsewhere~\cite{mrozek:etal:p21b}.

%%%%%%%%%%%%%%%%%%%%%%%%%%%%%%%%%%%%%%%%%%%%%%%%%%%%%%%%%%%%%%%%%
%%%%%%%%%%%%%%%%%%%%%%%%%%%%%%%%%%%%%%%%%%%%%%%%%%%%%%%%%%%%%%%%%
\subsection{Chaos in a Lorenz-type isolating neighborhood}

It is well-known that while for planar systems the existence of periodic
orbits can also be established using the Poincar\'e-Bendixson theorem,
considering the case of three- or higher-dimensional flows is 
significantly more involved~\cite{fuller:52a}. Nevertheless, we demonstrate
in the present subsection that the main ideas from our branched manifold
example can easily be carried over to the three-dimensional situation. Our 
example is partially motivated by the discussion in~\cite{ghys:13a}, and
it provides a different approach to establishing chaos in Lorenz-type equations
from the approach developed in~\cite{mischaikow:mrozek:95a, mischaikow:mrozek:98a,
mischaikow:etal:01a}.

\begin{figure}[tb]
  \begin{center}
  \includegraphics[width=0.70\textwidth]{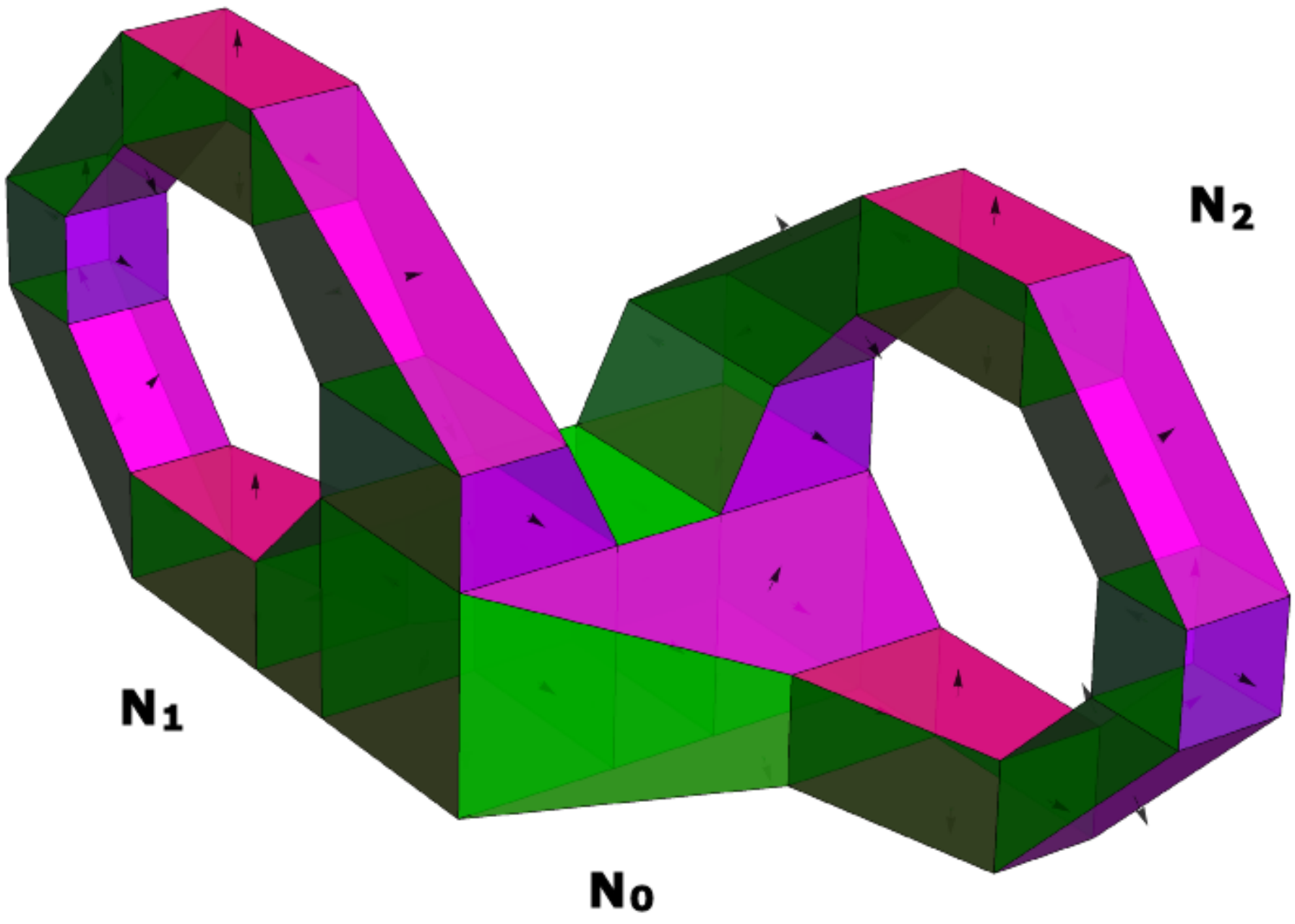}
  \caption{An isolated invariant set of a combinatorial multivector field on
  a cellular structure in three dimensions. Two-dimensional faces in the mouth
  of three-dimensional cells are marked with a small arrow. The union of the
  cells is an isolating block~$N$ for every flow which crosses cell faces
  transversally in the direction indicated by the arrows. The two-dimensional
  cells in the exit set~$N^-$ are marked in magenta. As in the branched 
  manifold example, every such flow contains infinitely many periodic
  orbits inside~$N$.}
  \label{fig:lorenz3d}
  \end{center}
\end{figure}

Figure~\ref{fig:lorenz3d} presents an example of an isolated invariant set~$\cS$
of a combinatorial multivector field~$\cV$ on a cellular structure~$\cX$ in three
dimensions. To improve visibility, only cells contained in~$\cS$ are presented.
Each three-dimensional cell is a cube-like polytope with convex, quadrilateral
faces. Each multivector $V\in\cV$ consists of a three-dimensional cell and its
faces except the faces indicated by a small arrow. Let~$N$ denote the union of
the closures of all three-dimensional cells in~$\cS$. Furthermore, let~$\phi$
denote a flow in~$\RR^3$ which is transversal to all two-dimensional faces.
Then one can easily verify that~$N$ is an isolating block for~$\phi$.

Note that~$N$ is the union of a rectangular cuboid~$N_0$ in the middle and
two handles attached to it --- the sets~$N_1$ on the left and~$N_2$ on the right.
It is easy to check that both~$N_0 \cup N_1$ and~$N_0 \cup N_2$ are isolating
neighborhoods for~$\phi$ with the Conley index of a hyperbolic periodic orbit.
Therefore, an argument analogous to the one discussed in the previous section
may be used to show that~$\phi$ has infinitely many periodic orbits in~$N$.
A more detailed discussion of this example based on the Lefschetz fixed point
theorem and combinatorial Poincar\'e maps will be presented
in~\cite{mrozek:etal:p21b}.

%%%%%%%%%%%%%%%%%%%%%%%%%%%%%%%%%%%%%%%%%%%%%%%%%%%%%%%%%%%%%%%%%
%%%%%%%%%%%%%%%%%%%%%%%%%%%%%%%%%%%%%%%%%%%%%%%%%%%%%%%%%%%%%%%%%
%%%%%%%%%%%%%%%%%%%%%%%%%%%%%%%%%%%%%%%%%%%%%%%%%%%%%%%%%%%%%%%%%
\section{Preliminaries}
\label{sec:prelims}

Throughout this paper, we denote the sets of positive integers, integers,
rational numbers, reals, and nonnegative reals by $\NN$, $\ZZ$, $\QQ$,
$\RR$, and $\RR_0^+$, respectively. In addition, given a set $X$ we denote
the family of subsets of $X$ by $\cP(X)$. Recall that a family $\cA\subset\cP(X)$
is a {\em partition} of $X$ if the elements of $\cA$ are mutually disjoint and
their union equals $X$.

For the remainder of this section, we collect basic notions from topology
and dynamical systems which will be needed later on.

%%%%%%%%%%%%%%%%%%%%%%%%%%%%%%%%%%%%%%%%%%%%%%%%%%%%%%%%%%%%%%%%%
%%%%%%%%%%%%%%%%%%%%%%%%%%%%%%%%%%%%%%%%%%%%%%%%%%%%%%%%%%%%%%%%%
\subsection{Basic notions in topology}

Given a topological space~$X$ and a subset $A\subset X$, we denote by $\inte A$,
$\cl A$, and $\bd A$ the interior, the closure, and the boundary of $A$,
respectively. Furthermore, we recall that the set~$A$ is called {\em locally
closed}, if the difference $\mo A:=\cl A \setminus A$ is closed. The
set~$\mo A$ is called the {\em mouth} of~$A$.

%%%%%%%%%%%%%%%%%%%%%%%%%%%%%%%%%%%%%%%%%%%%%%%%%%%%%%%%%%%%%%%%%
\subsubsection{Pairs of spaces and singular homology}

Throughout this paper, we consider \emph{topological pairs}, i.e., pairs
of topological spaces~$(X,A)$ for which $A\subset X$ and the topology
of~$A$ is induced by~$X$. In this setting, a single space~$X$ is treated
as the pair~$(X,\emptyset)$. In addition, we use the singular
homology~$H(X,A)$ of pairs with coefficients in the field of rational
numbers~$\QQ$, as defined for example in~\cite{dold:95a, spanier:81a}.

%%%%%%%%%%%%%%%%%%%%%%%%%%%%%%%%%%%%%%%%%%%%%%%%%%%%%%%%%%%%%%%%%
\subsubsection{Adjunctions, retracts, and ANRs}

For a topological pair~$(X,A)$, a topological space~$Y$,
and a continuous map $f\colon A\to Y$, the \emph{adjunction} $X\cup_f Y$ is defined
as the quotient space of the disjoint union of~$X$ and~$Y$ by the identification of
$x\in A$ with $f(x)\in Y$. We treat~$X\cup_f Y$ as the disjoint union of~$X\setminus A$
and~$Y$, with a dual treatment of the points in~$A$ and~$f(A)$.

Let~$(X,A)$ be a topological pair. A map $r\colon X\to A$ is called a
\emph{retraction} if $r(a)=a$ for all $a\in A$. In that case~$A$ is called
a \emph{retract} of~$X$. A metrizable space~$X$ is called an \emph{absolute
neighborhood retract} (which is abbreviated as ANR) if, whenever embedded as
a closed subset of a metric space, it is a retract of an ambient neighborhood.
If both~$X$ and~$Y$ are compact ANRs, then the same holds true for
the adjunction~$X\cup_f Y$, see for example \cite[V(2.9)]{borsuk:67a}
or~\cite[Theorem~5.6.1]{vanmill:89a}.

%%%%%%%%%%%%%%%%%%%%%%%%%%%%%%%%%%%%%%%%%%%%%%%%%%%%%%%%%%%%%%%%%
\subsubsection{Finite topological spaces}
\label{sec:top}

We close this discussion of topological prerequisites with a brief review
of finite topological spaces. For this, recall that any finite set~$X$ which
is equipped with a topology is called a \emph{finite topological space}. While
such spaces usually have very weak separation properties, those that satisfy
the~$T_0$ separation axiom are of particular interest. Recall that~$X$ is a
$T_0$-space, if for every pair of points~$x \neq y$ in~$X$ there exists an
open set~$U \subset X$ which contains exactly one of these two points. It is
a consequence of the celebrated Alexandrov Theorem~\cite{alexandrov:37a} that
a~$T_0$ topology on a finite set~$X$ may be functorially identified with a
partial order on~$X$. There are two dual ways to make the identification by
taking as open sets either the upper or the lower sets with respect to the
partial order. In this paper we identify open sets with upper sets, which
corresponds to
\begin{displaymath}
  x \le y
  \quad\Leftrightarrow\quad
  x \in \cl \{ y \}.
\end{displaymath}
Of particular interest for our applications is the fact that finite collections
of subsets of a given set automatically generate finite $T_0$-spaces in the
following way. Consider a finite family of subsets $\cX\subset\cP(X)$ of the
set~$X$. Since the inclusion relation on~$\cX$ defines a partial order, the
Alexandrov Theorem induces a unique $T_0$
topology on~$\cX$. In the case of such an Alexandrov topology, and to avoid
potential confusion if~$X$ is also a topological space, we denote the interior,
closure, and mouth of a subset $V \in \cX$ by $\Int V$, $\Cl V$ and $\Mo V$,
respectively. Since in the case of the Alexandrov topology the intersection of
an arbitrary family of open sets is open, the smallest open set containing a
given set $V \in \cX$ is open. We denote it by~$\Opn V$. Furthermore, if
$V=\{\sigma\}$ is a singleton, then we drop the braces and write $\Cl\sigma$
or $\Opn\sigma$ instead of $\Cl\{\sigma\}$ or $\Opn\{\sigma\}$. Finally,
observe that
\begin{eqnarray}
  \Opn\sigma&=&\setof{\tau\in\cX\mid \tau\supset\sigma},\label{eq:Opn}\\
  \Cl\sigma&=&\setof{\tau\in\cX\mid \tau\subset\sigma},\label{eq:Cl}\\
  \tau\in\Cl\sigma&\iff& \sigma\in\Opn\tau.\label{eq:Opn-Cl}
\end{eqnarray}
These characterizations will prove to be useful later on.

Questions of connectedness can easily be addressed for finite topological
spaces. It was shown in~\cite[Proposition~1.2.4]{barmak:11a} that a finite
topological space~$X$ is connected if and only if it is path-connected.
Moreover, $X$ is connected if and only if between any two points
$x,y \in X$ there exists a \emph{fence}, i.e., there exists a sequence
of points~$x = x_0, x_1, \ldots, x_n = y$ such that for every $k \in \{ 0,
\ldots, n-1 \}$ we have either $x_k \le x_{k+1}$ or $x_k \ge x_{k+1}$.
We would like to point out that one can view the length of a shortest
fence between~$x$ and~$y$ as a rudimentary measure of distance between
points in a connected finite topological space.

%%%%%%%%%%%%%%%%%%%%%%%%%%%%%%%%%%%%%%%%%%%%%%%%%%%%%%%%%%%%%%%%%
%%%%%%%%%%%%%%%%%%%%%%%%%%%%%%%%%%%%%%%%%%%%%%%%%%%%%%%%%%%%%%%%%
\subsection{Basic notions in dynamical systems}

We now collect some basic definitions and results from the theory of
dynamical systems.

%%%%%%%%%%%%%%%%%%%%%%%%%%%%%%%%%%%%%%%%%%%%%%%%%%%%%%%%%%%%%%%%%
\subsubsection{Semiflows and flows}

Suppose that $X$ is a locally compact metric space. By a \emph{semiflow}
on $X$ we mean a continuous map $\phi \colon X \times \RR_0^+ \to X$ such
that both
\[
  \phi(x,0)=x
  \quad\mbox{ and }\quad
  \phi(\phi(x,t),s)=\phi(x,t+s)
\]
are satisfied for all $x \in X$ and $s,t \in \RR_0^+$. If in the last identity 
the group $\RR$ replaces the semi-group $\RR_0^+$, then $\phi$ is called a \emph{flow}.
In both cases we call $X$ the \emph{phase space} of $\phi$. Furthermore, in order to
increase readability we sometimes abbreviate $\phi(x,t)$ by the short-hand $x\cdot t$,
as long as~$\phi$ is clear from the context.

Suppose now that $\TT = \RR_0^+$ or $\TT = \RR$, and let~$\phi : X \times \TT \to X$ denote
a semiflow or flow, respectively. Then a {\em solution} of~$\phi$ through $x\in X$ is a
map $\gamma: I\to X$ such that $I\subset\RR$ is a nontrivial interval containing zero, the
identity $\gamma(0)=x$ holds, and $\phi(\gamma(s),t) = \gamma(s+t)$ for all $t\in\TT$ and
$s\in I$ such that also $s+t\in I$. Note that in the case of a flow any two solutions through~$x$
coincide on the intersection of their domains. We denote the set of solutions through~$x$
by $\Sol(x,\phi)$. In the special case $I=\RR$, we further say that the solution~$\gamma$
is {\em full}. We say that a solution~$\gamma : I \to X$ through~$x$ is a {\em backward
solution}, if we have $\max I = 0$, and the collection of all backward solutions
through~$x$ is abbreviated by~$\Sol^-(x,\phi)$. Dually, the notion of {\em forward solution}
is for the case $\min I = 0$ and leads to the set~$\Sol^+(x,\phi)$. Note that every
point~$x \in X$ admits forward solutions, i.e., we always have $\Sol^+(x,\phi) \neq \emptyset$.
In general, however, backward solutions only exist in the case of a flow. If~$\phi$ is
a semiflow there may be points without backward solutions, so-called {\em start points}.
While such points can exist in many semiflows, the following result rules them out in a
number of interesting situations.
\begin{prop}[No start points on manifolds]
\label{prop:nostartpoint}
  Let~$\phi : X \times \RR_0^+ \to X$ be a semiflow on a metric space~$X$, and suppose
  that~$x \in X$ has an open neighborhood which is homeomorphic to~$\RR^n$ for some
  $n \in \NN$. Then~$x$ is not a start point, i.e., we have~$\Sol^-(x,\phi) \neq
  \emptyset$.
\end{prop}
This proposition is originally due to H.~Halkin, and its proof can be found for example
in~\cite[Theorem~11.8]{bhatia:hajek:69a} or~\cite[Theorem~4.1]{ciesielski:11a}.
For later use, we also have to introduce the following
convenient notation. Suppose that $\gamma \in \Sol^\pm(x,\phi)$ is a backward or forward
solution through~$x$ with domain~$I$. Then if $J \subset \RR$ is any interval we slightly
abuse notation and define
\begin{equation} \label{def:gammaval}
  \gamma(J) := \gamma( J \cap I )
  \quad\mbox{ where }\quad
  \gamma : I \to X.
\end{equation}
This convention will be useful for our definitions of exit and entry sets below.

While the notion of solution incorporates time dependence, the following terms focus on the
traversed image. The \emph{positive semi-trajectory} through~$x$ is the set $x \cdot \RR_0^+$,
and a \emph{negative semi-trajectory} through~$x$ is~$\gamma((-\infty,0])$, where
$\gamma : I \to X$ denotes a solution through~$x$ for which $(-\infty,0] \subset I$. Note
that if~$\phi$ is a flow, then the negative semi-trajectory through~$x$ is uniquely determined
and equal to~$x \cdot (-\infty,0]$. A \emph{trajectory} through~$x$ is the union of the
positive semi-trajectory and a negative semi-trajectory through~$x$. A point $x$ is called
\emph{stationary}, if its positive semi-trajectory is equal to the one-point set~$\{x\}$.
The point $x$ is called \emph{periodic}, if it is not stationary and if there exists a $T>0$ 
such that $\phi(x,T)=x$. In this case the trajectory through~$x$ is called a \emph{periodic
trajectory} or a \emph{periodic orbit}. The minimal positive time~$T$ for which
$\phi(x,T) = x$ is called the \emph{basic period} of the trajectory. By the
\emph{$\omega$-limit set of $x$} we mean  
\[
\omega(x):=\bigcap_{t>0}\clos{\phi(x,[t,\infty))}.
\]
If the space~$X$ is compact, then $\omega(x)\neq\emptyset$ for each $x\in X$. For an
arbitrary set $A\subset X$, the \emph{invariant part} $\Inv(A,\phi)$ of~$A$ is the set
of all $x\in A$ such that there exists a full solution through~$x$ with values in~$A$.
If the invariant part of~$A$ is equal to~$A$, then the set~$A$ is called \emph{invariant}.
In particular, each $\omega$-limit set is invariant.

%%%%%%%%%%%%%%%%%%%%%%%%%%%%%%%%%%%%%%%%%%%%%%%%%%%%%%%%%%%%%%%%%
\subsubsection{Sections}

For proving the existence of periodic orbits, we will make extensive use of the notion
of sections as considered in~\cite{mccord:etal:95a}. For this, let $\epsilon>0$, and for
$A\subset X$ define the \emph{$\epsilon$-left collar of $A$} as
\[
  L_\epsilon(A) :=
  \left\{ x\in X\mid \phi(x,(0,\epsilon))\cap A\neq \emptyset \right\}.
\]
Then a \emph{section for~$\phi$} is a subset $\Sigma \subset X$ such that there 
exists an $\epsilon > 0$ with the following properties:
\begin{itemize}
\item[(i)] The $\epsilon$-left collar~$L_\epsilon(\Sigma)$ of~$\Sigma$ is open, and
\item[(ii)] if $x\in L_\epsilon(\clos \Sigma)$, then the set 
$\phi(x,(0,\epsilon)) \cap \clos{\Sigma}$ consists of exactly one element.
\end{itemize}
Let~$\Sigma$ be a section. Then for every $x\in L_\epsilon(\Sigma)$ we can define
$\xi(x) \in \Sigma$ as the unique point given by property~(ii) above. According to
\cite[Proposition~3.1]{mccord:etal:95a}, the so-defined map~$\xi : L_\epsilon(\Sigma)
\to \Sigma$ is continuous, which immediately implies that an open subset of a section
is a section as well. 

Let $N\subset X$ and let $\Sigma$ be a section. We call $\Sigma$ a \emph{global section
in $N$} (or, as in  \cite{mccord:etal:95a}, a \emph{Poincar\'e section in $N$} if $N$
is an isolating neighborhood) if for each $x\in N$ its positive semi-trajectory
intersects $\Sigma$. If $X=N$, we simply call $\Sigma$ a \emph{global section}.

Being able to construct global sections for a semiflow is crucial for
proving the existence of periodic orbits. For our results, we make use of 
an approach due to Fuller~\cite{fuller:65a}, which can be described as follows.
Denote the unit circle in the complex plane by $S^1$, and let $f\colon X \to S^1$
be continuous. For every point $x \in X$ let $F_x \colon \RR_0^+ \to \RR$ denote
an arbitrary lift of the map
\[
  \RR_0^+\ni t \mapsto f(\phi(x,t))\in S^1
\]
with respect to the covering map $\exp  i (\cdot)$, i.e., suppose that
$\exp i F_x(t)=f(\phi(x,t))$ holds for all~$t \in \RR_0^+$. 
Following \cite{fuller:65a} we then define the map
\[
  \Delta f:=\Delta_\phi f\colon X\times \RR_0^+\to \mathbb R
  \quad\mbox{ as }\quad
  \Delta f(x,t):=F_x(t)-F_x(0).
\]
One can easily see that this map is independent of the choice of lift~$F_x$
and that it is continuous. Moreover, we have both $\Delta f(x,0)=0$ and
$f(\phi(x,t)) = f(x) \exp i \Delta f(x,t)$ for $t \in \RR_0^+$, as well as
\begin{equation}
\label{eq:delta}
  \Delta f(x,t+s)=\Delta f(x,t)+\Delta f(\phi (x,t),s)
  \quad\mbox{ for all }\quad
  t,s \in \RR_0^+.
\end{equation}
The function~$\Delta f$ provides the following mechanism for the
construction of a global section.
\begin{prop}[Global sections via angular function positivity]
\label{prop:sec}
Assume that~$X$ is compact, and consider a continuous angular function
$f\colon X \to S^1$. If the above-defined function~$\Delta f$ satisfies
the strict inequality $\Delta f(x,t)>0$ for every $x\in X$ and $t> 0$,
then for each point $z\in S^1$ the inverse image $f^{-1}(z)$ is a
global section.
\end{prop}
\proof
We need to show that the set~$\Sigma:=f^{-1}(z)$ is a global section.
For this, note that in view of~\eqref{eq:delta} and the assumption of the
proposition the map $\Delta f(x,\cdot)$ is strictly increasing for every $x\in X$.

We claim first that for every $x \in X$ the positive semi-trajectory of~$x$ has to
intersect $\Sigma$. Suppose this is not the case. Then the map $\Delta f(x,\cdot)$
has to be bounded, and therefore it converges as $t \to \infty$. But then the set
$f(\omega(x))$, which is non-empty due to the assumed compactness of~$X$, can consist
of only one element, which in fact is equal to $\lim_{t\to\infty} f(x)\exp  i
\Delta f(x,t)$. This in turn implies that the function $\Delta f(y,\cdot)$ is
constant for each $y\in\omega(x)$, contrary to our assumption. This establishes
the claim.

We now return to the discussion of~$\Sigma$. Due to the compactness of~$\Sigma$
there exists an $\epsilon>0$ such that 
\begin{equation}
\label{eq:edd}
\Delta f(\Sigma, [0,\epsilon])\subset [0,\pi).
\end{equation}
It follows immediately from~\eqref{eq:edd} that the $\epsilon$-left
collar~$L_\epsilon(\Sigma) = L_\epsilon(\clos \Sigma)$ satisfies (ii).

For the proof of (i), assume finally that $x\in L_\epsilon (\Sigma)$. Then the
inclusion~\eqref{eq:edd} implies that there exists a neighborhood~$U$ of~$x$
such that $f(\phi(y,[0,\epsilon])) \subset S^1\setminus \{-z\}$ for every $y\in U$.
Consider now the map $\alpha\colon S^1\setminus \{-z\}\to (\theta-\pi,\theta+\pi)$
defined as the inverse of $\exp i (\cdot)$, where we use the abbreviation
$\alpha(z)=\theta$. Then for each $y\in U$ one has
\[
\Delta f(y,t)=\alpha(f(\phi(y,t)))-\alpha(f(y)),
\] 
and therefore $\alpha(f(\phi(y,\cdot)))$ is strictly increasing on the interval
$[0,\epsilon]$. It follows, in particular, that both $\alpha(f(x))<\theta$ and
$\alpha(f(\phi(x,\epsilon)))>\theta$ are satisfied --- and the same inequalities
hold if~$x$ is replaced by~$y$ from some neighborhood~$V$ of~$x$. Therefore, for
each $y\in U\cap V$ there exists a $t_y\in (0,\epsilon)$ such that $f(\phi(y,t_y))=z$.
This establishes the inclusion $U\cap V\subset L_\epsilon(\Sigma)$, which proves~(i).
Thus, the set~$\Sigma$ satisfies all the required conditions and the result follows.
\qed\medskip

The above Proposition~\ref{prop:sec} shows that if the phase space~$X$ is compact,
then the notion of a \emph{surface of section} considered in~\cite{fuller:65a} is a
particular case of the notion of a global section as it was defined above. Therefore,
one can now apply~\cite[Theorem~1]{fuller:65a} to obtain the following important
result on the construction of global sections, which significantly weakens the necessary
assumptions.
\begin{prop}[Global sections via positivity at discrete times]
\label{prop:fuller}
  Assume that~$X$ is compact, and consider a continuous angular function
  $f\colon X \to S^1$. If for every $x\in X$ there exists a time $t_x>0$
  such that $\Delta_\phi f(x,t_x)>0$, then there exists a global section
  for the semiflow~$\phi$.
\end{prop}

%%%%%%%%%%%%%%%%%%%%%%%%%%%%%%%%%%%%%%%%%%%%%%%%%%%%%%%%%%%%%%%%%
\subsubsection{Isolated invariant sets and isolating blocks}

Recall that a subset~$A \subset X$ is called invariant
if it equals the invariant part of~$A$, i.e., if we have $A = \Inv(A,\phi)$.
In the study of dynamical systems, invariant sets take a prominent role. Unfortunately,
however, general invariant sets can be extremely sensitive to perturbations and
change dramatically. To remedy this, Conley~\cite{conley:78a} introduced the concept
of isolated invariant set. For this, let~$N \subset X$ be a compact nonempty subset.
Then~$N$ is called an \emph{isolating neighborhood} if $\Inv(N,\phi) \subset \inte N$, and
we call $S = \Inv(N,\phi)$ an \emph{isolated invariant set}. In other words, the 
invariant set~$S$ is maximal in its neighborhood~$\inte N$. As it turns out, this
maximality property makes it much easier to study isolated invariant sets, as one
can employ topological tools such as the Conley index~\cite{conley:78a}.

In general, an isolated invariant set~$S$ can have many associated isolating 
neighborhoods~$N$. For the purposes of this paper, we will focus our attention
to specific isolating neighborhoods called isolating blocks. For this,
let~$B \subset X$ be a given compact set and define its \emph{weak exit} and
\emph{weak entry set} respectively as
\begin{displaymath}
  \begin{array}{rcl}
    \DS B^{w-} & := & \DS \left\{ x \in \bd B \; \mid \;
      \forall \epsilon > 0 \mbox{ we have }
      \phi(x,(0,\epsilon))\cap (X \setminus B) \neq \emptyset \right\}, \\[1ex]
    \DS B^{w+} & := & \DS \left\{ x \in \bd B \; \mid \;
      \forall \gamma \in \Sol^-(x,\phi) \mbox{ and }
      \forall \epsilon > 0 \mbox{ we have }
      \gamma((-\epsilon,0)) \cap (X \setminus B) \neq \emptyset \right\},
  \end{array}
\end{displaymath}
where we use the convention introduced in~(\ref{def:gammaval}). We would like to
point out that if~$x \in \bd B$ is a start point, then one automatically has 
$x \in B^{w+}$ in view of $\Sol^-(x,\phi) = \emptyset$. Using these notions,
we call a compact set~$B \subset X$ an \emph{isolating block}, if
\begin{equation} \label{def:isoblock1}
  \bd B = B^{w-} \cup B^{w+}
  \qquad\mbox{ and }\qquad
  B^{w-} \mbox{ is closed}.
\end{equation}
In order to illustrate this definition, note that $x \in \bd B \setminus
(B^{w-} \cup B^{w+})$ if and only if there exists a solution $\gamma :
(-\epsilon_1, \epsilon_2) \to X$ through~$x$ with $\gamma((-\epsilon_1,
\epsilon_2)) \subset B$, i.e., this solution creates an internal semiflow
or flow tangency. Thus, the first condition in~(\ref{def:isoblock1}) rules
out such internal tangencies. As for the second condition, consider
the \emph{exit-time function} of a set~$B$ defined as
\begin{equation} \label {def:exittime}
  \tau_B \;\; \colon \;\;
  B \ni x \mapsto \sup\{t \geq 0 \;\mid\; \phi(x,[0,t]) \subset B \}
  \in [0,\infty].
\end{equation}
For general sets, the exit time function does not need to be continuous.
This changes if we consider isolating blocks, as the following celebrated
result shows.
\begin{prop}[Wa\.zewski's theorem]
\label{prop:block}
  Let $B \subset X$ be an isolating block for the semiflow~$\phi$. Then the
  following hold:
  \begin{itemize}
  \item[(a)] The exit-time function $\tau_B : B \to [0,\infty]$ is continuous.
  \item[(b)] We have $x \in B^{w-}$ if and only if $\tau_B(x) = 0$.
  \item[(c)] For every~$x \in B$ with $\tau_B(x) < \infty$ one has
             $\phi(x,\tau_B(x)) \in B^{w-}$.
  \item[(d)] The exit-time function~$\tau_B$ is bounded if and only if we have
             $\Inv(B,\phi) = \emptyset$.
  \end{itemize}
\end{prop}
The proof of parts~(a), (b), and~(c) can be found in Lemmas~2.1 and~2.2
in~\cite{srzednicki:04a}, and in~(d) it is clear that the boundedness of~$\tau_B$ 
implies $\Inv(B,\phi) = \emptyset$. Moreover, if $\tau_B(x) = \infty$,
then~$\emptyset \neq \omega(x) \subset \Inv(B,\phi)$ due to the compactness
of~$B$ and~\cite[Proposition~2.6(iii)]{diekmann:etal:95a}.

We note here that if~$B$ is an isolating block, then~$(B,B^{w-})$ is a very special
case of an index pair, a concept used to define the Conley index. In particular,
we have the following proposition.
\begin{prop}[Conley index of an isolating block]
\label{prop:Conley-index}
If~$B \subset X$ is an isolating block, then the homological Conley index
of~$\Inv(B,\phi)$ is given by the relative homology~$H_*(B,B^{w-})$.
\qed
\end{prop}

Despite the label ``weak,'' the above notions for exit and entry set are the ones
that are usually used to define isolating blocks, see for example~\cite{conley:78a,
mischaikow:mrozek:02a, srzednicki:04a}. However, for our applications in this paper
we need to consider a more restrictive notion of isolating blocks, which is 
characterized by more regular exit and entry behavior. For this, assume again
that~$B \subset X$ is a compact set and define the \emph{strong exit} and
\emph{strong entry sets} respectively via
\begin{equation} \label{def:exitentry}
  \begin{array}{rcl}
    \DS B^- & := & \DS \left\{ x \in \bd B \; \mid \;
      \exists \epsilon > 0 \mbox{ such that }
      \phi(x,(0,\epsilon)) \cap B = \emptyset \right\}, \\[1ex]
    \DS B^+ & := & \DS \left\{ x \in \bd B \; \mid \;
      \forall \gamma \in \Sol^-(x,\phi) \;
      \exists \epsilon > 0 \mbox{ such that }
      \gamma((-\epsilon,0)) \cap B = \emptyset \right\}.
  \end{array}
\end{equation}
While it is clear from the definitions that both $B^- \subset B^{w-}$ and
$B^+ \subset B^{w+}$ are satisfied, one can easily see that in general equality
does not hold. Moreover, it is not immediately obvious that assuming the identity
$\bd B = B^- \cup B^+$ and the closedness of~$B^-$ always imply that a compact
set~$B \subset X$ is an isolating block. Nevertheless, the following
proposition shows that this is in fact the case.
\begin{prop}[Properties of stronger exit and entry sets]
\label{prop:strongexit}
  Let $B \subset X$ be a compact set. Then the following hold:
  \begin{itemize}
  \item[(a)] If the strong exit set~$B^-$ is closed, then we have $B^{w-} = B^-$.
  \item[(b)] If~$B^-$ is closed and $\bd B = B^- \cup B^+$, then~$B$ is
             an isolating block.
  \end{itemize}
\end{prop}
\proof
It is clear that~(b) follows from~(a) and the inclusions $B^- \subset B^{w-}$
and $B^+ \subset B^{w+}$. In order to prove~(a), note first that we only have
to establish the inclusion $B^{w-} \subset B^-$. Suppose therefore that there
exists a point $x \in B^{w-} \setminus B^-$. Then $x \notin B^-$ implies that
\begin{displaymath}
  \phi(x,(0,\epsilon)) \cap B \neq \emptyset
  \quad\mbox{ for all }\quad
  \epsilon > 0,
\end{displaymath}
and in view of $x \in B^{w-}$ one has
\begin{displaymath}
  \phi(x,(0,\epsilon)) \cap (X \setminus B) \neq \emptyset
  \quad\mbox{ for all }\quad
  \epsilon > 0.
\end{displaymath}
Thus, for every $n \in \NN$ there exist times $0 < a_n < b_n < 1/n$ such that
\begin{displaymath}
  \phi(x,a_n) \in B
  \quad\mbox{ and }\quad
  \phi(x,b_n) \notin B.
\end{displaymath}
Let $c_n = \sup\{ t \in [a_n,b_n] \mid \phi(x,t) \in B \}$. Then we have
$\phi(x,c_n) \in B$ since~$B$ is closed, and this in turn implies the strict
inequality $c_n < b_n$. Moreover, $\phi(x,(c_n,b_n]) \subset X \setminus B$
due to the properties of the supremum, and therefore $\phi(x,c_n) \in B^-$
for all $n \in \NN$. Since we assumed that~$B^-$ is closed, the inequality
$0 < c_n < 1/n$ shows that $x = \lim_{n \to \infty} \phi(x,c_n) \in B^-$,
which contradicts our choice of~$x$ and proves the result.
\qed\medskip

As we will see in the next sections, the classical notion of isolating
block does not suffice for our purposes. We therefore make use of the
following stronger notion, which is motivated by the above result.
\begin{defn}[Strong isolating block]
\label{def:strongisoblock}
Let~$\phi$ denote a semiflow on a locally compact metric space~$X$,
and let~$B \subset X$ be a compact set. Then~$B$ is called a
\emph{strong isolating block} if 
\begin{displaymath}
  \bd B = B^- \cup B^+
  \qquad\mbox{ and }\qquad
  B^\pm \mbox{ are both closed},
\end{displaymath}
where~$B^\pm$ denote the strong entry and exit sets defined
in~\eqref{def:exitentry}.
\end{defn}

In view of Proposition~\ref{prop:strongexit} every strong isolating
block is also an isolating block in the classical sense. But, as we 
will see in the next section, the additional assumptions lead to
better-defined exit and entry behavior.

%%%%%%%%%%%%%%%%%%%%%%%%%%%%%%%%%%%%%%%%%%%%%%%%%%%%%%%%%%%%%%%%%
%%%%%%%%%%%%%%%%%%%%%%%%%%%%%%%%%%%%%%%%%%%%%%%%%%%%%%%%%%%%%%%%%
%%%%%%%%%%%%%%%%%%%%%%%%%%%%%%%%%%%%%%%%%%%%%%%%%%%%%%%%%%%%%%%%%
\section{Periodic orbits via bricks and the Conley index}
\label{sec:periodicex}

In this section we present a method for establishing the existence
of periodic orbits in a semiflow or flow. Our approach is based
on a decomposition of part of the phase space into smaller pieces,
which are specific types of strong isolating blocks whose invariant
part is empty. The main result is the subject of Theorem~\ref{thm:bricks},
where we use this decomposition in combination with the Conley index to
obtain a periodic orbit. We would like to point out, however, that a
different approach based on transfer maps in homology and the Lefschetz
number will be presented in~\cite{mrozek:etal:p21b}.

This section is organized as follows. The first subsection is devoted
to a study of the underlying concept of bricks, while the second subsection
introduces paths of bricks induced by a semiflow. Finally, the last
subsection contains the statement and proof of our main result.

%%%%%%%%%%%%%%%%%%%%%%%%%%%%%%%%%%%%%%%%%%%%%%%%%%%%%%%%%%%%%%%%%
%%%%%%%%%%%%%%%%%%%%%%%%%%%%%%%%%%%%%%%%%%%%%%%%%%%%%%%%%%%%%%%%%
\subsection{Bricks and brick paths}

The main goal of this subsection is the formulation of combinatorial properties
that imply the existence of a periodic orbit. For this, it is essential to
have a decomposition of part of the phase space which allows one to track the
flow behavior, which is based on the notion of strong isolating blocks defined
above.
\begin{defn}[Bricks, proper pairs of bricks, and brick paths]
\label{def:bricks}
Let~$\phi$ denote a semiflow on a locally compact metric space~$X$,
and let~$B \subset X$. Then~$B$ is called a \emph{brick} if the following
conditions are satisfied:
\begin{itemize}
\item[(i)]  The set~$B$ is a strong isolating block in the sense of
            Definition~\ref{def:strongisoblock}.
\item[(ii)] The invariant part of~$B$ is empty, i.e., we have
            $\Inv(B,\phi) = \emptyset$.
\end{itemize}
A pair~$(B,\tilde{B})$ of two bricks~$B$ and~$\tilde{B}$ is called a
\emph{proper pairs of bricks}, if both
\begin{displaymath}
  B \neq \tilde{B}
  \quad\mbox{ and }\quad
  B \cap \tilde{B} = B^- \cap \tilde{B}^+ \neq \emptyset
\end{displaymath}
are satisfied. Finally, a finite sequence $(B_0,\ldots,B_k)$ or an infinite
sequence $(B_0,B_1,\ldots)$ of bricks is called a \emph{brick path}, if either
one has $k=0$ or else the pair $(B_j,B_{j+1})$ is a proper pair of bricks for
all $j = 0,\ldots,k-1$ or all $j \in \NN_0$, respectively. If, moreover, we
have $B_i \neq B_j$ for all $i \neq j$, then the path is called \emph{injective}.
\end{defn}
The above notions are central for the results of this paper. Note that in view
of Proposition~\ref{prop:strongexit} every brick~$B$ is an isolating block
with~$B^{w-} = B^-$. Furthermore, Proposition~\ref{prop:block} shows that
the exit-time function~$\tau_B$ is continuous and bounded, and that we have
$x \in B^-$ if and only if $\tau_B(x) = 0$. One can also easily show that
if~$(B,\tilde{B})$ is a proper pair of bricks, then one has
\begin{displaymath}
  B \cap \tilde{B} = B^-\cap \tilde{B}^+ =
  B \cap \tilde{B}^+ = B^- \cap \tilde{B} \neq \emptyset,
\end{displaymath}
since due to the assumed compactness of~$B$ and~$\tilde{B}$ both the
inclusions~$B^- \subset B$ and $\tilde{B}^+ \subset \tilde{B}$ are satisfied.

The concept of bricks is a strengthening of isolating blocks that will
prove to be extremely useful. For us, they are the building blocks of a
phase space decomposition which allows for the use of combinatorial techniques
in a flow or semiflow setting. More precisely, we consider the following
decompositions.
\begin{defn}[Brick decomposition]
\label{def:brickdecomp}
Let~$\phi$ denote a semiflow on a locally compact metric space~$X$,
and let~$A \subset X$ be a compact flow region of interest which is not
necessarily invariant, i.e., its exit set~$A^{w-}$ could be nonempty. Then
a \emph{brick decomposition} of~$A$ is a finite collection~$\mathcal A$ of
subsets of~$A$ such that the following hold:
\begin{itemize}
\item[(i)]   The collection~$\mathcal A$ covers~$A$, i.e., we have
             $A = \bigcup_{\sigma\in{\mathcal A}}\sigma$.
\item[(ii)]  Every set $\sigma \in {\mathcal A}$ is a brick in the sense of
             Definition~\ref{def:bricks}.
\item[(iii)] If the bricks $\sigma_1,\sigma_2 \in {\mathcal A}$ are such that both
             $\sigma_1 \neq \sigma_2$ and $\sigma_1 \cap \sigma_2 \neq \emptyset$
             are satisfied, then either $(\sigma_1,\sigma_2)$ or $(\sigma_2,\sigma_1)$
             is a proper pair of bricks, i.e., we have either
             \begin{displaymath}
               \sigma_1\cap \sigma_2 = \sigma_1^- \cap \sigma_2^+
               \quad\mbox{ or }\quad
               \sigma_1\cap \sigma_2 = \sigma_1^+ \cap \sigma_2^-,
             \end{displaymath}
             respectively.
\end{itemize}
\end{defn}
For later use, we note the following consequence of the last of these assumptions.
\begin{lem}[Nontrivial solution pieces are contained in at most one brick]
\label{lem:uniquebrick}
  Let~$\mathcal A$ be a brick decomposition of the compact set~$A \subset X$
  in the sense of Definition~\ref{def:brickdecomp}. Then for every~$x \in A$
  and every~$t > 0$ with $\phi(x,[0,t]) \subset A$ there exists at most
  one brick~$\sigma \in {\mathcal A}$ which satisfies the inclusion
  $\phi(x,[0,t]) \subset \sigma$.
\end{lem}
\proof
If~$\phi(x,[0,t])$ is not contained in any brick we are clearly done. Suppose
now that there are two bricks $\sigma_1 \neq \sigma_2$ such that both $\phi(x,[0,t])
\subset \sigma_1$ and $\phi(x,[0,t]) \subset \sigma_2$ are satisfied. This implies
$x \in \sigma_1 \cap \sigma_2$, and Definition~\ref{def:brickdecomp}(iii) further
yields either $x \in \sigma_1^-$ or $x \in \sigma_2^-$. But then one either has
$\phi(x,(0,\epsilon)) \subset X \setminus \sigma_1$ or $\phi(x,(0,\epsilon))
\subset X \setminus \sigma_2$, respectively, for some $\epsilon > 0$, which
contradicts $\phi(x,[0,t]) \subset \sigma_1 \cap \sigma_2$.
\qed\medskip

We would like to point out that in Definition~\ref{def:brickdecomp}, we do not pose
any further dynamical requirements on the set~$A$. Nevertheless, having a brick
decomposition does impose constraints on~$A$, as the following result shows.
\begin{prop}[Exit set characterization via a brick decomposition]
  \label{prop:xxx}
  Let~$\mathcal A$ be a brick decomposition of the compact set~$A \subset X$
  in the sense of Definition~\ref{def:brickdecomp}. Then we have:
  \begin{itemize}
  \item[(a)] The weak exit set and the exit set of~$A$ coincide, i.e.,
             one has $A^{w-} = A^-$.
  \item[(b)] For every $x \in A$ the inclusion $x \in A^-$ is satisfied if and
             only if we have $\tau_\sigma(x) = 0$ for every $\sigma \in \mathcal A$
             with $x\in\sigma$.
  \end{itemize}
\end{prop}
\proof
To prove~(a), suppose that $x \in A^{w-}$. If $\sigma \in {\mathcal A}$ is a
brick which does not contain~$x$, then due to the compactness of~$\sigma$ there
exists an $\epsilon_\sigma > 0$ with $\phi(x,(0,\epsilon_\sigma)) \cap \sigma =
\emptyset$. On the other hand, if $x \in \sigma \in {\mathcal A}$, then
$x \in A^{w-}$ and Definition~\ref{def:brickdecomp}(i) imply $x \in \sigma^{w-}$,
and therefore~$x \in \sigma^-$ due to Proposition~\ref{prop:strongexit}(a).
Thus, also in this case there exists an $\epsilon_\sigma > 0$ such that
$\phi(x,(0,\epsilon_\sigma)) \cap \sigma = \emptyset$ is satisfied. If we now
define $\epsilon = \min\{ \epsilon_\sigma \mid \sigma \in {\mathcal A} \} > 0$,
then we immediately obtain $\phi(x,(0,\epsilon)) \cap A = \cup_{\sigma \in {\mathcal A}}
(\phi(x,(0,\epsilon)) \cap \sigma) = \emptyset$, i.e., one has $x \in A^-$.
Together with $A^- \subset A^{w-}$ this completes the proof of~(a).

As for~(b), suppose first that $x \in A^-$, and let $\sigma \in {\mathcal A}$
be such that $x \in \sigma$. Then there exists an $\epsilon > 0$ with
$\phi(x,(0,\epsilon)) \subset X \setminus A \subset X \setminus \sigma$, and
this shows $x \in \sigma^-$, as well as $\tau_\sigma(x) = 0$ due to
Proposition~\ref{prop:block}(b) and Proposition~\ref{prop:strongexit}(b).

Finally, suppose that $x \in A$ and that $\tau_\sigma(x) = 0$ holds for every
$\sigma \in \mathcal A$ with $x\in\sigma$.  In view of
Propositions~\ref{prop:block}(b) and~\ref{prop:strongexit}(b), combined with the
finiteness of the collection~$\mathcal A$, there exists an $\epsilon > 0$
such that $\phi(x,(0,\epsilon)) \subset X \setminus \sigma$ for every~$\sigma$
with $x \in \sigma \in \mathcal A$. Moreover, the union $\bigcup_{x \notin \sigma
\in \mathcal A} \sigma$ is compact, and it does not contain~$x$. Hence, one can
find an $0 < \epsilon_0 < \epsilon$ such that the set $\phi(x,(0,\epsilon_0))$
is disjoint from each $\sigma \in \mathcal A$, i.e., it is disjoint from~$A$.
This implies $x \in A^-$, and completes the proof of the proposition.
\qed\medskip

Despite the implications of Proposition~\ref{prop:xxx}, in general it is not
true that a set~$A$ which has a brick decomposition is an isolating block
itself. For example, one can easily see that if the bricks create a
non-convex corner in~$A$, then the exit set~$A^-$ might no longer be
closed and the semiflow~$\phi$ could have an internal tangency.

%%%%%%%%%%%%%%%%%%%%%%%%%%%%%%%%%%%%%%%%%%%%%%%%%%%%%%%%%%%%%%%%%
%%%%%%%%%%%%%%%%%%%%%%%%%%%%%%%%%%%%%%%%%%%%%%%%%%%%%%%%%%%%%%%%%
\subsection{Brick paths induced by a semiflow}

Suppose now that we are given a forward solution of a semiflow~$\phi$,
which originates at some point in a region with a brick decomposition.
As the following result shows, this leads to a well-defined sequence of
bricks that are traversed by this solution.
\begin{prop}[Brick paths induced by a brick decomposition]
\label{prop:yyy}
  Let~$\mathcal A$ be a brick decomposition of the compact set~$A \subset X$
  in the sense of Definition~\ref{def:brickdecomp}. Then for every 
  $x \in A \setminus A^-$ exactly one of the following two statements holds:
  \begin{itemize}
  \item[(a)] There exists a unique path $(\sigma_0,\sigma_1,\ldots,\sigma_n)$
             of bricks $\sigma_k \in {\mathcal A}$ with $\sigma_k \neq \sigma_{k+1}$
             and a finite collection of times $0 = t_0 < t_1 < \ldots < t_{n+1}$ with
             \begin{equation} \label{prop:yyy:1}
               \phi(x,[t_j,t_{j+1}])\subset \sigma_j
               \qquad\mbox{ and }\qquad
               \phi(x,t_{j+1}) \in \sigma_j^-
             \end{equation}
             for all $j=0,1,\ldots,n$, as well as $\phi(x,t_{n+1}) \in A^-$ and
             $\tau_A(x) = t_{n+1}$.
  \item[(b)] There exists a unique path $(\sigma_0,\sigma_1,\ldots)$ of infinitely
             many bricks $\sigma_k \in {\mathcal A}$ with $\sigma_k \neq \sigma_{k+1}$
             and a sequence of times $0 = t_0 < t_1 < \ldots$ such that~(\ref{prop:yyy:1})
             is satisfied for all $j \in \NN_0$. Furthermore, we have $t_j \to \infty$
             as $j \to \infty$, as well as $\tau_A(x) = \infty$.
  \end{itemize}
\end{prop}
\proof
Due to Proposition~\ref{prop:xxx}(b) and our choice $x \in A \setminus A^-$
there exists a brick $\sigma_0 \in {\mathcal A}$ with $\tau_{\sigma_0}(x) > 0$.
Furthermore, in view of Lemma~\ref{lem:uniquebrick} this brick is
uniquely defined. We now define $\sigma_0 = 0$ and set $t_1 = \tau_{\sigma_0}(x)$.
Then Proposition~\ref{prop:block}(c) and Proposition~\ref{prop:strongexit} imply
the inclusion $\phi(x,t_1) \in \sigma_0^-$, and we clearly have $\phi(x,[t_0,t_1])
\subset \sigma_0$.

Suppose now that we have constructed times~$0 = t_0 < \ldots < t_{n+1}$ and
bricks~$\sigma_0,\ldots,\sigma_n$ such that~(\ref{prop:yyy:1}) holds
for all $j = 0,\ldots,n$, where $n \in \NN_0$. Let~$x_{n+1} = \phi(x,t_{n+1})$.
According to~(\ref{prop:yyy:1}) one has $x_{n+1} \in \sigma_n^- \subset A$.
If in addition we have $x_{n+1} \in A^-$, then we have established alternative~(a)
of the proposition.

On the other hand, if $x_{n+1} \notin A^-$, then $x_{n+1} \in A \setminus A^-$, and
as above one can construct a unique brick~$\sigma_{n+1}$ with $x_{n+1} \in \sigma_{n+1}$
and $\tau_{\sigma_{n+1}}(x_{n+1}) > 0$, as well as $\sigma_{n+1} \neq \sigma_n$. If we
now define the time $t_{n+2} = t_{n+1} + \tau_{\sigma_{n+1}}(x_{n+1})$, then we have both
\begin{displaymath}
  \phi(x,t_{n+2}) =
  \phi(\phi(x,t_{n+1}),\tau_{\sigma_{n+1}}(x_{n+1})) =
  \phi(x_{n+1},\tau_{\sigma_{n+1}}(x_{n+1})) \in \sigma_{n+1}^-
\end{displaymath}
and
\begin{displaymath}
  \phi(x,[t_{n+1},t_{n+2}]) =
  \phi(x_{n+1},[0,\tau_{\sigma_{n+1}}(x_{n+1})]) \subset
  \sigma_{n+1}.
\end{displaymath}
Thus, (\ref{prop:yyy:1}) holds for $j = 0,\ldots,n+1$. If one proceeds
recursively, this leads to an infinite sequence of bricks and times as
in~(b), unless of course one has $x_j \in A^-$ for some~$j \in \NN$ as noted
above.

In order to complete the proof of the proposition, it only remains to be
shown that if the above procedure leads to an infinite sequence of times,
then we have $t_j \to \infty$ as $j \to \infty$. Suppose this is not
the case, i.e., the increasing sequence $(t_j)_{j \in \NN}$ is bounded.
Then the limit
\begin{displaymath}
  t_* = \lim_{j \to \infty} t_j \;\mbox{ exists }
  \qquad\mbox{ and }\qquad
  x_* = \phi(x,t_*) = \lim_{j \to \infty} x_j = 
    \lim_{j \to \infty} \phi(x,t_j) .
\end{displaymath}
Since the brick collection~${\mathcal A}$ is finite, there exists a
brick~$\sigma_* \in {\mathcal A}$ and a subsequence $(t_{k_j})_{j \in \NN}$
of times such that $\sigma_{k_j} = \sigma_*$ for all $j \in \NN$. Furthermore,
according to our above construction we have $x_{k_j} \in \sigma_{k_j} \cap
\sigma_{k_j - 1}^-$ and $\phi(x,[t_{k_j-1},t_{k_j}]) \subset \sigma_{k_j - 1}$,
which immediately yields $x_{k_j} \in \sigma_{k_j}^+ = \sigma_*^+$ in view of
Lemma~\ref{lem:uniquebrick} and Definition~\ref{def:brickdecomp}(iii). Now
recall that due to Definition~\ref{def:bricks}(i) the entry set~$\sigma_*^+$
is closed, and therefore $x_* = \lim_{j \to \infty} x_{k_j} \in \sigma_*^+$.
Thus, there exists an $\epsilon > 0$ such that on the one hand we have
\begin{displaymath}
  \phi(x,(t_* - \epsilon, t_*)) \subset X \setminus \sigma_*,
\end{displaymath}
while on the other hand the inclusions
\begin{displaymath}
  \phi(x,[t_{k_j},t_{k_j+1}]) \subset \sigma_{k_j} = \sigma_*
\end{displaymath}
are satisfied for all $j \in \NN$. Together with $t_* = \lim_{j \to \infty}
t_{k_j}$ this leads to a contradiction and completes the proof of the
proposition.
\qed\medskip

Notice that the above proof heavily relies on the assumed closedness of the
entry sets of bricks, even though this closedness is not required in order to
ensure that every brick is an isolating block.

%%%%%%%%%%%%%%%%%%%%%%%%%%%%%%%%%%%%%%%%%%%%%%%%%%%%%%%%%%%%%%%%%
%%%%%%%%%%%%%%%%%%%%%%%%%%%%%%%%%%%%%%%%%%%%%%%%%%%%%%%%%%%%%%%%%
\subsection{Periodic orbits via brick decompositions}

The last two subsections have laid the foundation for establishing the existence
of periodic orbits using combinatorial techniques. By introducing the notions of
brick paths and brick decompositions in Definitions~\ref{def:bricks}
and~\ref{def:brickdecomp} we have a purely combinatorial way of encoding 
potential semiflow or flow solution behavior inside a region of interest.
In particular, the notion of proper pair of bricks encodes the direction in
which solutions can move between bricks. Needless to say, however, not every 
combinatorial brick path necessarily leads to a semiflow or flow solution which
traverses the bricks in the given order. As the last subsection showed, one can
in general only create a unique brick path from a solution, not vice versa.

In the present subsection, we will show that under suitable conditions it is
possible to deduce the existence of a periodic orbit for the given semiflow 
or flow from the existence of an appropriate brick decomposition. Our result is
motivated by and based on the work~\cite{mccord:etal:95a}. While the main result
in~\cite{mccord:etal:95a} contains a fairly strong assumption concerning the
existence of a global section, in our approach we can deduce the existence of
this section from a fairly mild combinatorial requirement on brick paths. This
condition in turn relies on an appropriate coarsening of the given brick
decomposition. More precisely, we have the following main result of this
section. For this, recall that $\ZZ_p = \{ 0, 1, \ldots, p-1 \}$ denotes
the finite cyclic additive group of order~$p$ with respect to addition 
modulo~$p$.
\begin{thm}[Periodic orbits via brick decompositions]
  \label{thm:bricks}
  Let~$\mathcal A$ be a brick decomposition of the compact set~$A \subset X$
  in the sense of Definition~\ref{def:brickdecomp}. Assume that the following
  hold:
  \begin{itemize}
  \item[(a)] The set~$A$ is an isolating block for the semiflow or flow~$\phi$,
             and both the isolating block~$A$ and its exit set~$A^- = A^{w-}$
             is an ANR.
  \item[(b)] For either $r = 0$ or $r = 1$ we have
             \begin{displaymath}
               \dim H_{2n + r}(A,A^-) = \dim H_{2n+1 + r}(A,A^-)
               \quad\mbox{ for all }\quad
               n \in \ZZ,
             \end{displaymath}
             and not all of these homology groups are trivial.
  \end{itemize}
  In addition, we assume that the brick decomposition~${\mathcal A}$ has
  a coarsening into~$p \in \NN$ subcollections in the following sense.
  Suppose that for $i \in \ZZ_p$ there are subcollections $\emptyset \neq
  {\mathcal A}_i \subset {\mathcal A}$ which are pairwise disjoint,
  and such that $\bigcup_{i=0}^p {\mathcal A}_i = {\mathcal A}$. Let
  $A_i := \bigcup_{\sigma \in {\mathcal A}_i} \sigma \subset A$, then we
  clearly have $A = \bigcup_{i=0}^p A_i$. In the following, suppose further that
  whenever~$i$ or~$j$ are indices of these sets or subcollections,
  additions involving them are always performed in the cyclic
  group~$\ZZ_p$. For example, if $i = p-1$ then $A_{i+1} = A_0$.
  Then we assume that this coarsening has the following properties:
  \begin{itemize}
  \item[(c)] The number~$p$ of subcollections~${\mathcal A_i}$ satisfies $p \ge 3$.
  \item[(d)] We have $A_i \cap A_{i+1} \neq \emptyset$ for all $i \in \ZZ_p$.
  \item[(e)] We have $A_i \cap A_j = \emptyset$ for all $i,j \in \ZZ_p$
             with $j \notin \{ i+1, \; i, \; i-1 \}$.
  \item[(f)] For every index $i \in \ZZ_p$, every brick $\sigma_0\in\mathcal A_i$,
             and every maximal brick path~$(\sigma_0,\sigma_1,\ldots)$ in~${\mathcal A}$
             there exists a $k \geq 0$ such that either $\sigma_k^- \subset A^-$,
             or one has $\sigma_k \in {\mathcal A}_{i+1}$, while for all
             $\ell=1,\ldots,k-1$ we have $\sigma_\ell \notin {\mathcal A}_{i+2}$.
  \end{itemize}
  Then the isolating block~$A$ contains a non-trivial periodic orbit.
\end{thm} 
\proof
As mentioned before the formulation of the theorem, we prove our result by applying
Theorem~1.3 in~\cite{mccord:etal:95a}. The latter result guarantees the existence of a
periodic orbit as long as~(b) is satisfied, an isolating neighborhood of the invariant part
of~$A$ has a Poincar\'e section, the phase space is an ANR, and the semiflow has compact
attraction. Since~(b) is already assumed, in order to verify the other conditions we embed
the semiflow~$\phi$ restricted to~$A$ into a semiflow on a new compact ANR which has a
global section. 

For this, we use the following convention. If~$w$ and~$z$ are points of the unit circle~$S^1$,
then~$(w,z)$ and~$[w,z]$ denote the open and closed arc joining~$w$ and~$z$ in the positive
direction, respectively. Furthermore, for $m \in \ZZ$ we set $z_m := e^{2\pi i m / p}$.
For every $m = 0,\ldots,p-1$ we denote by~$g_m$ a Urysohn lemma map which transfers the
compact set~$A_m$ to the arc~$[z_m,z_{m+1}] \subset S^1$ such that $g_m|_{A_{m-1}\cap A_m}$
is constant equal to~$z_m$, and $g_m|_{A_m \cap A_{m+1}}$ is constant equal to~$z_{m+1}$.
Note that the existence of~$g_m$ follows from~(d) and~(e). Gluing all maps~$g_m$, one then
obtains a continuous map $g \colon A \to S^1$.

Let $Y := A \cup_{g|_{A^-}} S^1$ denote the adjunction. Furthermore, consider the
semiflow~$\Phi$ on the space~$Y$ defined as
\begin{displaymath}
  \Phi(y,t) :=
  \left\{ \begin{array}{ccccc}
    \DS \phi(y,t) & \mbox{ for } & y \in A & \mbox{ and } &
      0\leq t \leq \tau_A(y), \\[1ex]
    \DS e^{i (t-\tau_A(y))} g(\phi(y,\tau_A(y)) & \mbox{ for } & y \in A &
      \mbox{ and } & t \geq \tau_A(y), \\[1ex]
    \DS e^{i t} y & \mbox{ for } & y \in S^1 & \mbox{ and } & t \geq 0.
  \end{array} \right.
\end{displaymath}
One can readily verify that the continuity of~$\Phi$ is a consequence of
Proposition~\ref{prop:block}(a). Clearly, the adjunction~$Y$ is the required
new phase space, which is in fact a compact ANR. Moreover, the motion of~$\phi$
inside~$A$ is preserved by the semiflow~$\Phi$. Thus, it remains to show
that~$\Phi$ has a global section, since the one restricted to the interior
of~$A$ then satisfies our requirements. For this, we extend the map~$g$
to a map $G \colon Y \to S^1$ given by
\begin{displaymath}
  G(y) :=
  \left\{ \begin{array}{ccc}
    \DS g(y) & \mbox{ for } & y \in A, \\[1ex]
    \DS y    & \mbox{ for } & y \in S^1 .
  \end{array} \right.
\end{displaymath}
According to Proposition~\ref{prop:fuller}, the required section exists
if for every $y \in Y$ one can find a time $t_y > 0$ which satisfies
$\Delta_\Phi G(y,t_y) > 0$. This is immediately clear for $y \in S^1$,
since in this case $t_y > 0$ can be chosen arbitrarily.

Consider now a point $y \in A \setminus A^-$ with $\tau_A(y) < \infty$.
Then the identity $\Delta_\Phi G(y,t) = \Delta_\phi g(y,t)$ holds for
all $0 \leq t \leq \tau_A(y)$. If in addition one has the strict inequality
$\Delta_\phi g(y,\tau_A(y)) > 0$, then we can choose $t_y := \tau_A(y)$.
On the other hand, if this inequality does not hold, then in view
of~(\ref{eq:delta}) we can choose any $t_y > \tau_A(y) - \Delta_\phi
g(y,\tau_A(y))$. Either way, also in this case one obtains
$\Delta_\Phi G(y,t_y) > 0$, as required.

Finally, consider the remaining case of a point $y \in A \setminus A^-$
with $\tau_A(y) = \infty$. According to our assumptions there exists
an integer $m \in \ZZ_p$ such that $y \in A_m \setminus A_{m+1}$, which
in turn implies $g(y) = e^{2\pi i T}$ for some $m/p \le T < (m+1)/p$.
In order to prove the existence of a time~$t_y$ for which $\Delta_\phi
g(y,t_y) > 0$ we consider the infinite path $(\sigma_0,\sigma_1,\ldots)$ 
associated to~$y$ which is guaranteed by Proposition~\ref{prop:yyy}(b),
with associated transition times $0 = t_0 < t_1 < \ldots$. Moreover,
Proposition~\ref{prop:yyy}(a) guarantees that none of the exit
sets~$\sigma_j^-$ is contained in~$A^-$. Since $y \in \sigma_0$, we
either have $\sigma_0 \in {\mathcal A}_m$ or $\sigma_0 \in {\mathcal A}_{m-1}$.
Note that in the latter case one necessarily has $g(y) = z_m$.

We now distinguish two cases, namely $\sigma_0 \in \mathcal A_m$ and
$\sigma_0 \in \mathcal A_{m-1} \setminus \mathcal A_m$. To begin with,
suppose that $\sigma_0\in\mathcal A_m$, which immediately yields
\begin{equation} \label{eq:phi}
  \phi(y,[0,t_1]) \subset [z_m,z_{m+1}].
\end{equation}
According to our assumption~(f) there exists an index~$k$ such that 
$\sigma_k \in \mathcal A_{m+1}$, which gives
\begin{equation} \label{eq:ykk}
  g(\phi(y,t_k)) \in [z_{m+1},z_{m+2}],
\end{equation}
and none of the bricks $\sigma_1,\ldots,\sigma_{k-1}$ belongs to the
collection~$\mathcal A_{m+2}$. Together with \eqref{eq:phi} this implies
\begin{equation} \label{eq:emp}
  g(\phi(y,[0,t_k]) \cap (z_{m+2},z_{m+3}) = \emptyset.
\end{equation}
If additionally we select~$k$ as the minimal one satisfying~(f), then
we have $\sigma_{k-1} \in \mathcal A_m$. Together with \eqref{eq:ykk}
this implies $g(\phi(y,t_k))=z_{m+1}$. But then \eqref{eq:emp} shows
that 
\begin{displaymath}
  \Delta_\phi g(y,t_k)=2\pi\left(\tfrac{m+1}{p}-T\right)>0,
\end{displaymath}
and we therefore can choose $t_y := t_k$ in the first case.

Finally, we consider the second case $\sigma_0 \in \mathcal A_{m-1}
\setminus \mathcal A_m$. Under this assumption one obtains the
inclusion $g(y,[0,t_1]) \subset [z_{m-1},z_m]$. As before, it follows
from~(f) that there exists an index~$k_\ast$ such that $\sigma_{k_\ast}
\in {\mathcal A}_m$, and one also has
\begin{displaymath}
  g(\phi(y,[0,t_{k_\ast}])) \cap (z_{m+1},z_{m+2}) = \emptyset.
\end{displaymath}
If we again select~$k_\ast$ as minimal, then one further has $g(y_\ast) = z_m$
for $y_\ast := \phi(y,t_{k_\ast})$. Recall that $g(y) = z_m$, which immediately
shows that $\Delta_\phi g(y,t_{k_\ast}) = 0$. Thus, if we apply the previously
considered case to the point~$y_\ast$, then we have $T = m/p$, and we get the
existence of a time $t > 0$ such that $\Delta_\phi g(y_\ast,t) = 2\pi ((m+1)/p
- T) = 2\pi / p$. Finally, we can now choose the time $t_y := t_{k_\ast}+t$,
for which \eqref{eq:delta} yields the strict inequality
\begin{displaymath}
  \Delta_\phi g(y,t_y) = \Delta_\phi g(y_\ast,t) > 0,
\end{displaymath}
and this completes the proof of the theorem.
\qed\medskip
\begin{figure}[tb]
  \begin{center}
  \includegraphics[width=0.9\textwidth]{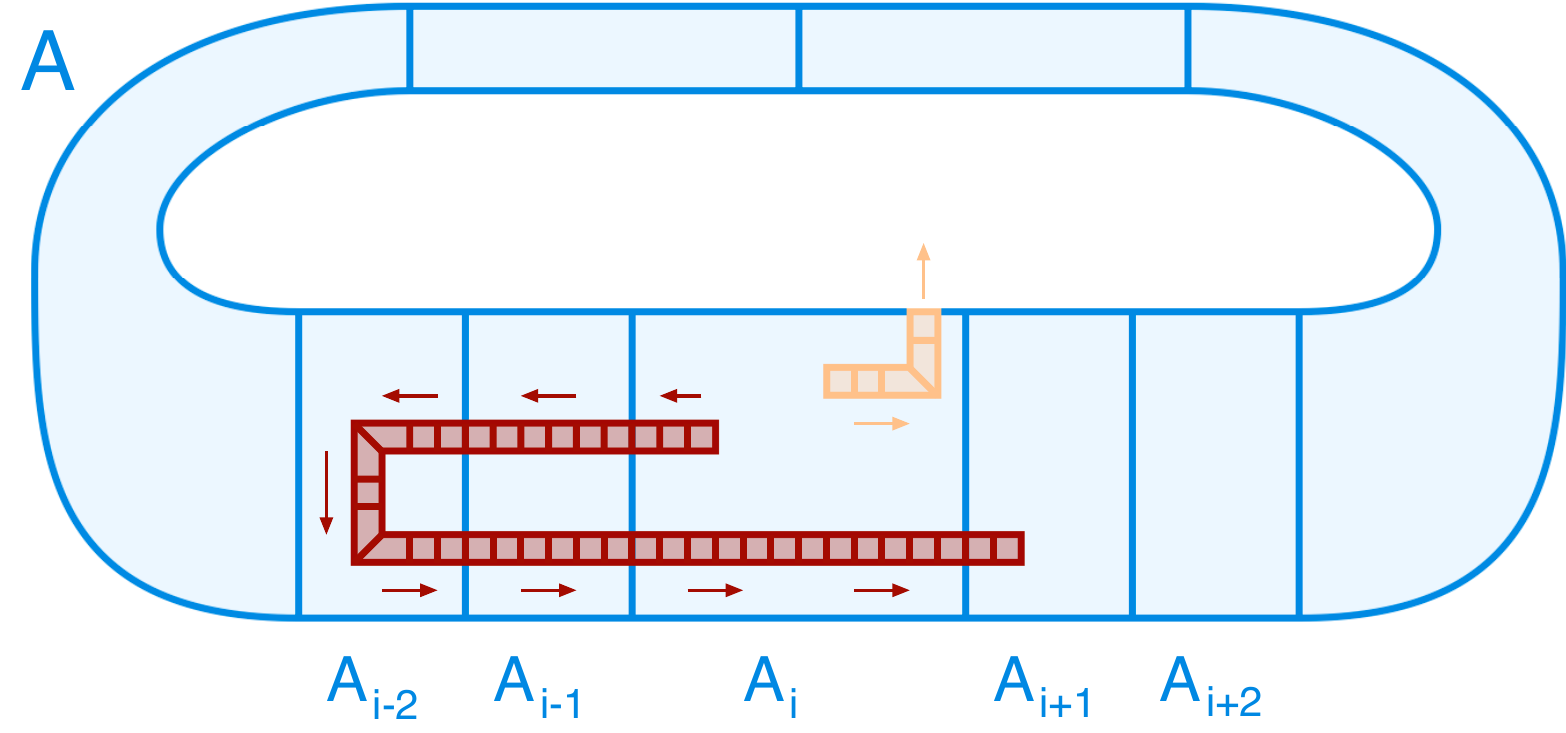}
  \caption{Illustration of assumption~(f) in Theorem~\ref{thm:bricks}.
           It is not necessary that the intersections~$A_i \cap
           A_{i+1}$ of successive segments are flow transverse. There 
           may be brick paths, such as the one shown in red, that 
           first backtrack to a segment~$A_j$ with $j < i$ before
           reaching~$A_{i+1}$, as long as it does not loop around to 
           reach~$A_{i+2}$ first. Also, brick paths are allowed to 
           reach the exit set~$A^-$, such as the one shown in orange.}
  \label{fig:brickpathcond}
  \end{center}
\end{figure}

For an illustration of the assumptions of the theorem we refer the reader
to Figure~\ref{fig:brickpathcond}. Notice that for the theorem to apply,
one needs to be able to subdivide the isolating block~$A$ into at least
three segments~$A_0, \ldots, A_{p-1}$ in such a way that there is a
well-defined overall semiflow motion in the direction of increasing
indices. We would like to point out, however, that the specific
assumptions are fairly mild:
\begin{itemize}
\item While the individual bricks in the collection~$\mathcal{A}$ have to
be isolating blocks that satisfy additional assumptions, the coarsened 
segments~$A_i$ do not need to be isolating blocks. In particular, the semiflow
behavior of~$\phi$ across the intersections~$A_i \cap A_{i+1}$ does not have
to be unidirectional, and there certainly can be internal flow tangencies 
to these intersections.
\item In order to ensure an overall semiflow direction, one only needs to
verify that bricks paths starting in the collection~$\mathcal{A}_i$ cannot
reach~$\mathcal{A}_{i+2}$ by looping around in the direction of decreasing
indices and thereby avoiding~$\mathcal{A}_{i+1}$. But they most certainly 
are allowed to first backtrack before reaching~$\mathcal{A}_{i+1}$.
\end{itemize}
In other words, as long as one can generate a sufficiently fine brick
decomposition, the coarsening which is required for the application of
Theorem~\ref{thm:bricks} will be achievable fairly easily.

%%%%%%%%%%%%%%%%%%%%%%%%%%%%%%%%%%%%%%%%%%%%%%%%%%%%%%%%%%%%%%%%%
%%%%%%%%%%%%%%%%%%%%%%%%%%%%%%%%%%%%%%%%%%%%%%%%%%%%%%%%%%%%%%%%%
%%%%%%%%%%%%%%%%%%%%%%%%%%%%%%%%%%%%%%%%%%%%%%%%%%%%%%%%%%%%%%%%%
\section{Recurrence in combinatorial topological dynamics}
\label{sec:ctdrecurrence}

In the last section we presented an approach for proving the
existence of periodic orbits in semiflows, which was based
on a suitable brick decomposition of part of the phase space. In the
present section, we will show that such decompositions can be 
constructed through combinatorial multivector fields on cellular
spaces.

Multivector fields as introduced in \cite{lipinski:etal:p20a, mrozek:17a}
constitute a combinatorial counterpart of classical vector fields. They
provide an abstract language and methods to design, analyze, and use
algorithms facilitating the automated investigation of classical dynamical
systems. Traditionally, graph theory is used in this role. Multivector fields
may be viewed as directed graphs whose vertices constitute a finite topological
space. In this way algorithms based on the methods of topological dynamics may
be invented, designed, and tested without the need of relating them to
conventional flows. Classical dynamics is needed only in applications to
conventional vector fields: First on input when the combinatorial multivector
field is constructed, and then on output when the results are interpreted.
But we would like to emphasize that in order to derive the brick decomposition
used in the last section, the language of multivector fields proves to be
an ideal tool.

The remainder of this section is organized as follows.
We begin by introducing the combinatorial framework which provides
the link between semiflows and the results of the previous section.
For this, Section~\ref{subsec51} recalls basic definitions and results
for cell complexes. This is followed in Section~\ref{subsec52} by the
notion of strongly admissible semiflows, whose behavior can be encoded by
studying the flow across cell boundaries. Section~\ref{subsec53} then
describes the concept of combinatorial multivector fields and shows how
they can be generated in a canonical way from strongly admissible
semiflows. After that, we turn our attention to how the combinatorial
multivector field~$\cV$ associated with an admissible semiflow~$\phi$
can be used to establish the existence of periodic orbits for~$\phi$.
This is accomplished by first showing in Section~\ref{subsec54} how
isolated invariant sets in the combinatorial setting give rise to
associated isolated invariant sets for admissible semiflows, and then
introducing combinatorial versions of Poincar\'e sections and Lyapunov
functions in Section~\ref{subsec55}. This last section concludes by
showing how the latter concepts can be used to verify the assumptions
of our main result.

%%%%%%%%%%%%%%%%%%%%%%%%%%%%%%%%%%%%%%%%%%%%%%%%%%%%%%%%%%%%%%%%%
%%%%%%%%%%%%%%%%%%%%%%%%%%%%%%%%%%%%%%%%%%%%%%%%%%%%%%%%%%%%%%%%%
\subsection{Cellular complexes}
\label{subsec51}

Throughout this subsection, let~$X$ denote a topological space.
Then a {\em $d$-dimensional cell}, or briefly a {\em cell}, is a subset
$\sigma\subset X$ which is homeomorphic to the closed unit ball in the
Euclidean space~$\RR^d$. In this case, we write $\dim\sigma := d$ for the
dimension of the cell. The associated {\em open cell} and {\em cellular
boundary} of the cell~$\sigma$ are defined as the image of the open unit
ball and of the unit sphere in~$\RR^d$, respectively, with respect to the
same homeomorphism used in the definition of~$\sigma$. Since a homeomorphism
of a closed ball in~$\RR^d$ has to preserve the associated open ball and
sphere, the open cell and the cellular boundary do not depend on the choice
of a particular homeomorphism. We denote the open cell by~$\opc{\sigma}$,
and the cellular boundary by~$\cebd\sigma$. We note that these concepts
have the following straightforward properties, whose proof will be left
to the reader.
\begin{prop}[Cell properties]
\label{prop:opc-cebd-properties}
Let $\sigma\subset X$ be a $d$-dimensional cell in a topological space.
Then the following identities are satisfied:
\begin{eqnarray}
  \opc{\sigma}\cup\cebd\sigma&=&\sigma,       \label{eq:sigma1}\\
  \opc{\sigma}\cap\cebd\sigma&=&\emptyset,    \label{eq:sigma2}\\
  \cebd\sigma&=&\sigma\setminus\opc{\sigma},  \label{eq:sigma2a}\\
  \cl\opc{\sigma}&=&\sigma.                   \label{eq:sigma3}
\end{eqnarray}
Note that the cell boundary~$\cebd\sigma$ is in general not the topological
boundary~$\bd\sigma$, and that the open cell $\opc{\sigma}$ is not
necessarily the topological interior of~$\sigma$.
\qed
\end{prop}
The notion of cell forms the basic building block for our combinatorial
study of dynamical systems. By combining cells in a suitable way, one can
arrive at a combinatorial description of certain topological spaces. More
precisely, we say that a finite family~$\cX$ of cells in a topological
space~$X$ is a {\em cellular decomposition} of~$X$, if the following
two conditions are met:
\begin{itemize}
   \item[(CD1)] If $\opc{\sigma}\cap\tau\neq\emptyset$, then we have the
                inclusion $\sigma\subset\tau$, for all cells $\sigma,\tau\in\cX$.
   \item[(CD2)] We have the representation $X=\bigcup_{\sigma\in\cX} \opc{\sigma}$.
\end{itemize}
In the following, we write~$\cX_k$ for the subfamily of all at most $k$-dimensional
cells in~$\cX$. We refer to $X_k:=\bigcup_{\sigma \in \cX_k}\sigma$ as the $k$-skeleton
of $X$. Given two cells $\tau,\sigma\in\cX$ such that $\tau\subset\sigma$
we say that~$\tau$ is a {\em face} of~$\sigma$, and that~$\sigma$ is a {\em coface}
of~$\tau$. The face or coface is called {\em proper}, if $\sigma\neq\tau$.
A face~$\tau$ of~$\sigma$ is a {\em facet} of~$\sigma$, if $\dim\sigma=\dim\tau+1$.
A cell is called {\em top-dimensional}, if it is not a proper face of a
higher-dimensional cell. We will refer to a top-dimensional cell also as
a {\em toplex}, and denote the family of all toplexes by~$\cXtop$. One can
easily see that a space~$X$ which has a cellular decomposition is equal to the
union of its toplexes. Finally, the {\em star} of a cell $\sigma\in\cX$ is the
family of all cofaces of~$\sigma$, and it is denoted by~$\Star\sigma$. We note
that the star of a toplex contains only this toplex.

Of course, not every topological space admits a cellular decomposition, and even
if it does, a cellular decomposition need not be unique. By a {\em cellular space}
we mean a metrizable topological space~$X$ which admits a cellular decomposition.
Since~$X$ is the union of its toplexes, the finiteness of~$\cX$ in turn implies
that~$X$ is automatically compact. 

A simple induction argument with respect to the cardinality of $\cX$ shows that
a cellular space is a finite, regular CW-complex. Vice versa, every finite regular
CW-complex is a cellular space~\cite[Section~I.3]{cohen:73a}.

In the sequel we assume that~$X$ is a fixed cellular space and that~$\cX$
denotes its fixed cellular decomposition. The following proposition is an
immediate consequence of assumption (CD1), and its simple proof is omitted.
\begin{prop}[Open cells form a partition]
\label{prop:open-cells-intersection}
Let~$X$ denote a cellular space. Then the collection of open cells of~$X$ forms
a partition of~$X$. In particular, for any two cells $\sigma,\tau\in\cX$ the
inequality $\opc{\sigma}\cap\opc{\tau}\neq\emptyset$ implies the equality
$\sigma=\tau$.
\qed
\end{prop}
It follows immediately from Proposition~\ref{prop:open-cells-intersection}
that for every point $x\in X$ there exists exactly one cell $\sigma\in\cX$
such that $x\in\opc{\sigma}$. We denote this unique cell by~$\cell{x}$.
Furthermore, as an easy consequence of Proposition~\ref{prop:open-cells-intersection}
one obtains the following result, whose straightforward proof is also omitted.
\begin{prop}[Cell representations]
\label{prop:cel-decomp-prop}
Suppose that~$\cX$ is a cellular decomposition of the metric space~$X$. Then
we have the representations
\begin{displaymath}
  \mbox{\rm (i)} \quad \sigma =
    \bigcup_{\cX \ni \rho \subset\sigma} \opc{\rho}, \qquad\quad
  \mbox{\rm (ii)} \quad \cebd\sigma = 
      \bigcup_{\cX \ni \rho \subsetneq\sigma} \rho, \qquad\quad
    \mbox{\rm (iii)} \quad \sigma\cap\tau = 
      \bigcup_{\cX \ni \rho \subset\sigma\cap\tau} \rho,
\end{displaymath}
where~$\sigma, \tau \in \cX$ are arbitrary cells.
\qed
\end{prop}
As we explained in Section~\ref{sec:top}, we will consider~$\cX$ as an
Alexandrov topological space with the Alexandrov topology induced by the
inclusion relation on~$\cX$. We would like to point out that the star of
a cell $\sigma\in\cX$ coincides with $\Opn\sigma$ in the Alexandrov topology.
In particular, a singleton which consists of a toplex is open in this topology.
It can easily be seen that the family of all stars of cells in~$\cX$ forms a
basis of the Alexandrov topology. For a family $\cA\subset\cX$ of cells we
define its {\em support} by
\begin{displaymath}
  |\cA|:=\bigcup_{\sigma\in\cA} \opc{\sigma}.
\end{displaymath}
In this case we say that~$\cA$ is a {\em combinatorial representation}
of~$|\cA|$. Moreover, a set $A\subset X$ is called {\em representable},
if there exists a family $\cA\subset\cX$ such that $A=|\cA|$. Note that
if such a family exists, it is uniquely determined by the set~$A$. As an
immediate consequence of Proposition~\ref{prop:open-cells-intersection}
we obtain the following result, which summarizes the properties of
support and representable sets.
\begin{prop}[Characterization of representable sets]
\label{prop:representable}
Suppose that~$\cX$ is a cellular decomposition of the metric space~$X$.
Then a subset $A\subset X$ is representable if and only if
\begin{equation} \label{eq:representable}
  \forall \sigma\in\cX \; : \;\;
  \opc{\sigma}\cap A\neq\emptyset \; \implies \; \opc{\sigma}\subset A.
\end{equation}
\end{prop}
\proof
Assume first that the set~$A$ is representable. Then we have 
$A=\bigcup_{\tau\in\cA} \opc{\tau}$ for some family $\cA\subset\cX$. 
Let $\sigma\in\cX$ be such that $\opc{\sigma}\cap A\neq\emptyset$. We
claim that $\sigma\in\cA$. Indeed, if not, then by
Proposition~\ref{prop:open-cells-intersection} we get
$\opc{\sigma} \cap \opc{\tau} = \emptyset$ for every $\tau\in\cA$,
and therefore
\begin{displaymath}
  \opc{\sigma}\cap A =
  \bigcup_{\tau\in\cA} \left( \opc{\sigma}\cap\opc{\tau} \right) =
  \emptyset,
\end{displaymath}
which is a contradiction. Hence, $\sigma\in\cA$ and, in consequence,
$\opc{\sigma}\subset \bigcup_{\tau\in\cA} \opc{\tau}=A$.

Suppose now that \eqref{eq:representable} holds. Let $\cA :=
\setof{\tau\in\cX\mid \opc{\tau}\cap A\neq\emptyset}$. From
\eqref{eq:representable} one immediately obtains $|\cA|\subset A$.
To verify the opposite inclusion, consider an arbitrary point $x\in A$.
Then in view of $x \in \opc{\cell{x}}$ we have $\opc{\cell{x}} \cap
A \neq \emptyset$. Hence, $\cell{x} \in \cA$ and $x \in |\cA|$.
\qed\medskip

The following simple result shows how set representations are
affected by set operations. Its simple proof is omitted.
\begin{prop}[Properties of set representations]
\label{prop:support}
Suppose that~$\cX$ is a cellular decomposition of the metric space~$X$,
and that $\cA,\cB\subset\cX$ are families of cells. Then
\begin{eqnarray}
  |\cA\cup\cB|&=&|\cA|\cup|\cB|,\label{eq:support1}\\
  |\cA\cap\cB|&=&|\cA|\cap|\cB|,\label{eq:support2}\\
  |\cA|\setminus|\cB|&=&|\cA\setminus\cB|,\label{eq:support3}\\
  X&=&|\cX|.\label{eq:support4}
\end{eqnarray}
In particular, the space~$X$ is representable, and the union,
intersection and difference of representable sets is again
representable.
\qed
\end{prop}
\begin{cor}[Representability of cell components]
\label{cor:sigma-representable}
Suppose that~$\cX$ is a cellular decomposition of the metric space~$X$.
Then the sets $\sigma$, $\opc{\sigma}$, and $\cebd\sigma$ are representable
for each $\sigma\in\cX$.
\end{cor}
\proof
The representability of $\sigma$ and $\opc{\sigma}$ follows immediately
from~(CD1) and Proposition~\ref{prop:representable}. Hence, the
representability of $\cebd\sigma$ follows from Proposition~\ref{prop:support}.
\qed\medskip

In our combinatorial approach to the study of dynamics the toplexes in
a cellular decomposition play a special role, since --- in some sense ---
they pass on the orbits amongst each other. This ``passing along'' is realized
along the intersections of toplexes, which can be described using the following
notion. We set
\begin{displaymath}
  \Fr(X):=|\cX\setminus\cXtop|
\end{displaymath}
and refer to~$\Fr(X)$ as the {\em frame} of~$X$. Then one can readily see that
the frame contains the union of all toplex intersections.

We close this subsection on cell complexes and cellular decompositions
with a sequence of results, which will be useful in subsequent sections.
\begin{prop}[Openness of the star]
\label{prop:star-open}
Suppose that~$\cX$ is a cellular decomposition of the metric space~$X$.
Then for every $\sigma\in\cX$ the support of its star is open in $X$.
\end{prop}
\proof
Let $\sigma\in\cX$. It suffices to prove that $A:=X\setminus |\Star\sigma|$
is closed. Hence, assume that $x\in\cl A$. Let $(x_n)$ be a sequence of points
in~$A$ converging to~$x$. Set $\sigma_n:=\cell{x_n}$. Then we clearly have
$\sigma_n\not\in\Star\sigma$. Since $\cX$ is finite, without loss of generality
we may assume that the identity $\sigma_n=\tau$ holds for some $\tau\in\cX$
and all $n\in\NN$. It follows that $x_n\in \opc{\tau}$ for all $n\in\NN$,
and by \eqref{eq:sigma3} we get $x\in\cl\opc{\tau}=\tau$. It follows that
$x \in \tau \cap \opc{\cell{x}}$. Therefore, property~(CD1) of a cellular
decomposition implies that $\cell{x}\subset\tau$. Since $\tau\not\in\Star\sigma$,
we get $\cell{x}\not\in\Star\sigma$. It follows that $x\not\in |\Star\sigma|$,
which means that $x\in A$. Hence, $A$ is closed and $|\Star\sigma|$ is open.
\qed
\begin{prop}[Support of the combinatorial closure of a cell]
\label{prop:supp-comb-cl-cell}
Suppose that~$\cX$ is a cellular decomposition of the metric space~$X$.
Then for every $\sigma\in\cX$ 
\[
|\Cl\sigma|=\sigma.
\]
\end{prop}
\proof
   Consider an arbitrary point $x\in |\Cl\sigma|$. Then $x\in\opc{\tau}$ for some $\tau\in\Cl\sigma$. It follows that
   $x\in \opc{\tau}\subset\tau\subset\sigma$. To see the opposite inclusion, consider now a point $x\in\sigma$.
   Let $\tau\in\cX$ be such that $x \in\opc{\tau}$. Then one has $\opc{\tau}\cap\sigma\neq\emptyset$,
   which in turn implies $\tau\subset\sigma$. Consequently, we have both $\tau\in\Cl\sigma$ and $x\in|\Cl\sigma|$.
\qed

\begin{prop}[Topological properties of toplex components]
\label{cor:sigma-comb-vs-top}
Suppose that~$\cX$ is a cellular decomposition of the metric space~$X$.
Let $\sigma,\tau \in \cXtop$ be such that $\sigma\neq\tau$. Then we have:
\begin{align}
   \opc{\sigma}\subset\inte\sigma,\label{eq:sigma-comb-vs-top-1}\\
   \bd\sigma\subset\cebd\sigma,   \label{eq:sigma-comb-vs-top-2}\\
   \opc{\sigma}\cap\cebd\tau=\emptyset, \label{eq:sigma-comb-vs-top-3}\\
   \bd\sigma\cap\bd\tau=\cebd\sigma\cap\cebd\tau.\label{eq:sigma-comb-vs-top-4}
\end{align}
\end{prop}
\proof
Since $\sigma\in\cXtop$, we have $\Star\sigma=\{\sigma\}$. Hence
$\opc{\sigma} = |\Star\sigma|$ is open in view of
Proposition~\ref{prop:star-open}, and it follows that $\opc{\sigma}
\subset\inte\sigma$, thereby establishing~\eqref{eq:sigma-comb-vs-top-1}.
Property~\eqref{eq:sigma-comb-vs-top-2} follows immediately
from~\eqref{eq:sigma-comb-vs-top-1} and~\eqref{eq:sigma2a}. To
verify~\eqref{eq:sigma-comb-vs-top-3}, suppose that $\opc{\sigma} \cap
\cebd\tau \neq \emptyset$. According to Corollary~\ref{cor:sigma-representable}
the set~$\tau$ is representable, and Proposition~\ref{prop:representable}
implies the inclusion $\sigma\subset\tau$. Since $\sigma\in\cXtop$, this in
turn yields $\sigma=\tau$, which contradicts the choice of~$\sigma$ and~$\tau$,
and therefore proves~\eqref{eq:sigma-comb-vs-top-3}. Finally, in order to
prove~\eqref{eq:sigma-comb-vs-top-4}, note that the choice of~$\tau$ and~$\sigma$,
in combination with~(CD1), implies the identity $\opc{\tau}\cap\sigma=\emptyset$,
which in turn gives $\cl\opc{\tau}\subset \cl(X\setminus\sigma)$. Hence,
property~\eqref{eq:sigma3} implies that $\cebd\sigma \cap \cebd\tau \subset
\sigma \cap \tau = \sigma\cap\cl\opc{\tau} \subset \sigma\cap
\cl(X\setminus\sigma) = \bd\sigma$. A symmetric argument proves that
$\cebd\sigma \cap \cebd\tau \subset \bd\tau$. Therefore, one has
$\cebd\sigma\cap\cebd\tau\subset\bd\sigma\cap\bd\tau$. The opposite
inclusion follows immediately from~\eqref{eq:sigma-comb-vs-top-2}.
This completes the proof.
\qed
\begin{prop}[The point-to-cell mapping]
\label{prop:sigma-map}
Suppose that~$\cX$ is a cellular decomposition of the metric space~$X$,
and consider the point-to-cell map
\begin{displaymath}
  \Cell \;\colon\; X \ni x \mapsto \cell{x} \in \cX.
\end{displaymath}
Then the following hold:
\begin{itemize}
   \item[(i)]   For every $\cA\subset\cX$ we have $\Cell^{-1}(\cA)=|\cA|$.
   \item[(ii)]  For every $\cA\subset\cX$ we have $\Cell(|\cA|)=\cA$.
   \item[(iii)] The point-to-cell map~$\Cell$ is a continuous, open surjection.
\end{itemize}
\end{prop}
\proof
Properties (i), (ii), and the surjectivity of~$\Cell$ follow immediately
from the definitions. From Proposition~\ref{prop:star-open} and property~(i)
we see that $\Cell^{-1}(\Star\sigma)$ is open in~$X$ for every $\sigma\in\cX$.
Since the stars of cells in~$\cX$ form a basis of the topology in~$\cX$, this
proves the continuity of the point-to-cell map~$\Cell$.

To prove that~$\Cell$ is an open map, consider an open set $U\subset X$.
We have to prove that~$\Cell(U)$ is open in~$\cX$. Thus, consider a cell
$\tau\in\Cell(U)$. We will be done if we prove that
\begin{equation} \label{eq:sigma-map}
  \Star\tau\subset \Cell(U).
\end{equation}
Choose a point $x\in U$ such that $\cell{x}=\tau$ is satisfied. Then we
clearly have $x\in\opc{\tau}$. In order to verify~\eqref{eq:sigma-map},
consider also a cell $\rho\in\Star\tau$. Since $\tau \subset \rho =
\cl\opc{\rho}$, and since the set~$U$ is open and contains~$x$, there exists
a point $x'\in U\cap\opc{\rho}$. But this inclusion immediately implies
$\rho = \Cell(x')$, as well as $\Cell(x') \in \Cell(U)$. Combined, they
yield $\rho \in \Cell(U)$, which in turn proves \eqref{eq:sigma-map}.
\qed\medskip

As an immediate consequence of Proposition~\ref{prop:sigma-map} and
property \eqref{eq:support3} of Proposition~\ref{prop:support} we
obtain the following corollary.
\begin{cor}[Topological properties of the support]
\label{cor:rep-opn-cls}
Suppose that~$\cX$ is a cellular decomposition of the metric space~$X$,
and let $\cA\subset\cX$ denote a family of cells. Then the support~$|\cA|$
is open (respectively closed) in~$X$ if and only if~$\cA$ is open (respectively
closed) in~$\cX$. In other words, a representable set is open (respectively
closed) in~$X$ if and only if its representation is open (respectively
closed) in~$\cX$.
\qed
\end{cor}
Needless to say, while on the level of representable sets openness and
closedness can be verified either in the metric space~$X$ or in the
Alexandrov space~$\cX$, in general the topology of the former space
is significantly larger.
\begin{prop}[Generating cellular decompositions]
\label{prop:sub-decompositions}
Suppose that~$\cX$ is a cellular decomposition of the metric space~$X$. 
Then the following hold:
\begin{itemize}
   \item[(i)]   $\cA$ is a cellular decomposition of~$|\cA|$ for every
                closed $\cA\subset\cX$.
   \item[(ii)]  $X\setminus\opc{\sigma}$ is a cellular space with
                $\cX\setminus\{\sigma\}$ as a cellular decomposition
                for every $\sigma\in\cXtop$.
   \item[(iii)] $\Fr(X)$ is a cellular space with $\cX\setminus\cXtop$
                as a cellular decomposition.
   \item[(iv)]  $X_k$, the $k$-skeleton of~$X$, is a cellular space
                with~$\cX_k$ as a cellular decomposition.
\end{itemize}
\end{prop}
\proof
  Observe that for every $\sigma\in\cA$ we have $\opc{\sigma}\subset|\cA|$, which
  together with~\eqref{eq:sigma3} implies the inclusion $\sigma\subset\cl |\cA|=|\cA|$,
  since~$|\cA|$ is closed in $X$ by Corollary~\ref{cor:rep-opn-cls}.
  Hence, the collection~$\cA$ is a family of cells not only in $X$ but also in~$|\cA|$.
  To see that~$\cA$ is a cellular decomposition of~$|\cA|$, observe that~(CD1) is inherited
  from~$\cX$ and~(CD2) follows immediately from the definition of support. This proves~(i).
  The remaining three assertions are straightforward consequences of~(i).
\qed\medskip

As an immediate consequence of Proposition~\ref{prop:sub-decompositions} and the
characterization of ANRs contained in~\cite[Theorem~IV.6.1]{borsuk:67a}, we obtain
the following proposition.
\begin{prop}[Cellular decompositions as ANRs]
\label{prop:anrs}
Suppose that~$\cX$ is a cellular decomposition of the metric space~$X$. 
If $\cA\subset\cX$ is closed, then~$|\cA|$ is a compact ANR.
\qed
\end{prop}
\begin{prop}[Acyclicity of the star]
\label{prop:star-acyclicity}
Suppose that~$\cX$ is a cellular decomposition of the metric space~$X$ and $\tau\in\cX$.
Then the set~$|\Star\tau| \subset X$ is acyclic.
\end{prop}
\proof
  Let~$n$ denote the cardinality of  $S:=\Star\tau$. We establish the result by
  induction on~$n$. If $n=1$, then we have $|\Star\tau|=\opc{\tau}$. Hence, it is
  acyclic as an open ball. Now, assume that $n>1$ and the claim holds for stars of
  cardinality less than $n$. Select a $\sigma\in \Star\tau$ of maximal dimension.
  Since $n>1$, we have $\sigma\neq\tau$. By Proposition~\ref{prop:sub-decompositions},
  the family $\cY:=\cX\setminus\{\sigma\}$ is a cellular decomposition of cellular
  space $Y:=|\cY|$, and the family $\cZ:=\Bd\sigma\setminus\{\sigma\}$ is a cellular
  decomposition of cellular space $Z:=|\cZ|=\cebd\sigma$. Clearly, $\tau\in\cY\cap\cZ$.
  Let~$S_Y$ denote the support of the star of~$\tau$ in~$\cY$, and let~$S_Z$ denote the
  support of the star of~$\tau$ in~$\cZ$. Without loss of generality we may assume
  that~$\sigma$ is a unit ball centered at $o\in\sigma$. Let $p:\sigma\setminus\{o\}\to\cebd\sigma$
  denote the radial projection and set $V_0:=p^{-1}(S_Z)$. Then one can easily verify that
  $U:=\opc{\sigma}$ and $V:=S_Y\cup V_0$ are both open in~$S$.
  Moreover, $p_{|V_0}:V_0\to S_Z$ is a deformation retraction of~$V_0$ onto~$S_Z$
  and $\id_{S_Y}\cup p_{|V_0}:V\to S_Y$ is a deformation retraction of~$V$ onto~$S_Y$.
  Since~$S_Y$ is the support of the star of~$\tau$ in~$\cY$ and~$S_Z$ is the support
  of the star of~$\tau$ in~$\cZ$, both of these supports are acyclic by our induction
  assumption. Hence, also~$V$ and $U\cap V=V_0$ are acyclic. Clearly, $U=\opc{\sigma}$
  is acyclic as an open ball. Therefore, the conclusion follows by applying the
  Mayer-Vietoris theorem to~$U$ and~$V$, since one easily verifies that $S=U\cup V$. 
\qed
\begin{prop}[Point-to-cell mapping induces an isomorphism in homology]
\label{prop:Sigma-iso}
Suppose that~$\cX$ is a cellular decomposition of the metric space~$X$ and
$\cB\subset\cA$ are closed sets in~$\cX$. Then the restriction of the map~$\Sigma$
defined in Proposition~\ref{prop:sigma-map} and given by
\begin{displaymath}
  \Sigma_{|\cA|}:  (|\cA|,|\cB|) \ni x \mapsto \cell{x}\in(\cA,\cB)
\end{displaymath}
is a well-defined and continous map of pairs. Moreover, the map induced in
singular homology and given by
\begin{displaymath}
  H_*(\Sigma_{|\cA|}): H_* (|\cA|,|\cB|)  \to  H_*(\cA,\cB)
\end{displaymath}
is an isomorphism.
\end{prop}
\proof
   The well-definedness of the definition and the continuity of~$\Sigma_{|\cA|}$
   follow immediately from Proposition~\ref{prop:sigma-map}.
   To prove that~$\Sigma_{|\cA|}$ induces an isomorphism in singular homology
   assume first that $\cB=\emptyset$. In this case, we need to verify that
   $H_*(\Sigma_{|\cA|}): H_*(|\cA|) \to  H_*(\cA)$ is an isomorphism. Note that
   for each $\sigma\in\cA$ we have $\Sigma^{-1}(\Star\sigma)=|\Star\sigma|$ by
   Proposition~\ref{prop:sigma-map}(i). Hence, $\Sigma^{-1}(\Star\sigma) =
   \Sigma^{-1}(\Opn\sigma)$ is acyclic by Proposition~\ref{prop:star-acyclicity} .
   Therefore, the conclusion for  $\cB=\emptyset$ follows
   from \cite[Theorem~2.2]{barmak:etal:20a}. The case of general~$\cB$
   then follows easily by applying the five lemma.
\qed

%%%%%%%%%%%%%%%%%%%%%%%%%%%%%%%%%%%%%%%%%%%%%%%%%%%%%%%%%%%%%%%%%
%%%%%%%%%%%%%%%%%%%%%%%%%%%%%%%%%%%%%%%%%%%%%%%%%%%%%%%%%%%%%%%%%
\subsection{Strongly admissible semiflows}
\label{subsec52}

Equipped with our toolset for studying cellular spaces, we can now turn
our attention to the study of semiflows on such spaces. For this, it
will be crucial to study the sequence of cells that is traversed by
solutions of the semiflow, and for this, the following notions will
be extremely useful.
\begin{defn}[Immediate future and past]
\label{def:ift:ipt}
Let~$X$ denote a cellular space with cellular decomposition~$\cX$, and
let~$\phi : X \times \RR_0^+ \to X$ denote a semiflow on~$X$.
\begin{itemize}
\item We say that a cell $\sigma\in\cXtop$ is the \emph{immediate future}
of a point~$x\in X$, if we have
\begin{displaymath}
  \phi(x,(0,\epsilon)) \subset \opc{\sigma}
  \quad\mbox{ for some }\quad
  \epsilon>0.
\end{displaymath}
In this case, we write $\ift(x) = \sigma$. 
\item We say that $\sigma\in\cXtop$ is the \emph{immediate past} of a
backward solution $\gamma \in \Sol^-(x,\phi)$ through the point $x\in X$, if
\begin{displaymath}
  \gamma((-\epsilon,0)) \subset \opc{\sigma}
  \quad\mbox{ for some }\quad
  \epsilon>0,
\end{displaymath}
and we write~$\ipt(\gamma) = \sigma$.
\end{itemize}
\end{defn}
Note that in the case of a flow~$\phi : X \times \RR \to X$ the immediate
past of a solution through~$x$ depends only on the point~$x$. Hence, since
in this case every point admits a backward solution, we can define the
immediate past of a point~$x$ as the immediate past of any of its backward
solutions, and we denote it by~$\ipt(x)$. Observe also that we have the
implication
\begin{equation} \label{eq:ift-x-sigma}
  \sigma=\ift{x} \quad \implies \quad x\in\sigma.
\end{equation}
Finally, we would like to point out that not every point~$x$ has to have an
immediate future, since the solution starting at~$x$ might oscillate between
two neighboring cells if~$x$ lies in their intersection. Thus, the existence
of an immediate future does impose certain regularity assumptions on a
semiflow, and similarly the existence of the immediate past. This is 
formalized in the following definition.
\begin{defn}[Admissible and strongly admissible semiflows]
\label{def:admissible}
Let~$X$ denote a cellular space with cellular decomposition~$\cX$, and
let~$\phi : X \times \RR_0^+ \to X$ denote a semiflow on~$X$. We say
that~$\phi$ is \emph{admissible with respect to~$\cX$}, if the following
three conditions are satisfied.
\begin{enumerate}
\item[(A1)] The immediate future of every point in~$X$ and the immediate
past of every backward solution are well-defined.
\item[(A2)] For every point $x \in \Fr(X)$ the immediate future of~$x$ is
different from the immediate past of any backward solution through~$x$.
\item[(A3)] For every $\sigma\in\cXtop$ the strong exit and entry
sets~$\sigma^-$ and~$\sigma^+$ are closed, where these sets are defined
in~\eqref{def:exitentry}.
\end{enumerate}
In addition, if for every cell~$\sigma$ its invariant part is empty unless
$H_*(\sigma,\sigma^-)$ is non-trivial, then we say that~$\phi$ is
\emph{strongly admissible}.
\end{defn}
The above notions are motivated by our earlier work~\cite{mrozek:wanner:21a},
where we used combinatorial vector fields to construct strongly admissible
semiflows on simplicial complexes. Note, however, that the underlying phase
space decomposition used in~\cite{mrozek:wanner:21a} is not always a
cellular decomposition.

The conditions required by Definition~\ref{def:admissible} impose a number
of constraints on the semiflow behavior. These allow us to draw several
useful conclusions. To begin with, we have the following characterizations
of strong exit and entry sets~$\sigma^-$ and~$\sigma^+$, respectively,
for every toplex $\sigma\in\cXtop$.
\begin{prop}[Exit and entry set characterizations]
\label{prop:sigma-pm}
Let~$X$ denote a cellular space with cellular decomposition~$\cX$,
and assume that~$\phi$ is an admissible semiflow with respect to~$\cX$.
Furthermore, let $\sigma\in\cXtop$ be arbitrary. Then
\begin{equation} \label{eq:sigma-minus}
  \sigma^- = \setof{x \in \bd\sigma\mid \ift(x)\neq\sigma}
\end{equation}
and
\begin{equation} \label{eq:sigma-plus}
  \sigma^+ = \setof{x \in \bd\sigma\mid
    \forall\gamma\in\Sol^-(x,\phi):\, \ipt(\gamma)\neq\sigma}.
\end{equation}
Moreover, if~$\phi$ is an admissible flow, then we have
\begin{equation} \label{eq:sigma-i}
  \sigma^+ = \setof{x \in \bd\sigma\mid \ipt(x)\neq\sigma}.
\end{equation}
\end{prop}
\proof
Assume first that $x\in\sigma^-$. Then the inclusion $x \in \bd\sigma$ is
satisfied, and there exists an $\epsilon>0$ such that $\phi(x, (0,\epsilon))
\cap \sigma=\emptyset$. This in turn implies $\ift(x)\neq\sigma$. Suppose
now conversely that $x \in \bd\sigma$ and that $\tau := \ift(x) \neq \sigma$.
Then one has $\phi(x, (0,\epsilon)) \subset \opc{\tau}$ for some constant
$\epsilon>0$. It follows from properties~\eqref{eq:sigma-comb-vs-top-2}
and~\eqref{eq:sigma-comb-vs-top-3} that then $\phi(x,(0,\epsilon)) \cap \sigma
= \emptyset$, which means that $x\in\sigma^-$. This proves \eqref{eq:sigma-minus}.

In order to verify \eqref{eq:sigma-plus}, we first assume that $x\in\sigma^+$
and let $\gamma\in\Sol^-(x,\phi)$. Then there exists an $\epsilon>0$ such that
$\gamma((-\epsilon,0)) \cap \sigma = \emptyset$. This yields $\ipt(\gamma)
\neq \sigma$ and proves that~$x$ belongs to the right hand side of
\eqref{eq:sigma-plus}. Now suppose conversely that~$x$ belongs to the
right-hand side of \eqref{eq:sigma-plus}. Let $\gamma \in \Sol^-(x,\phi)$
and assume that $\tau := \ipt(\gamma) \neq \sigma$. Then $\gamma((-\epsilon,0))
\subset \opc{\tau}$ for some constant $\epsilon>0$. As above, it follows from
properties~\eqref{eq:sigma-comb-vs-top-2} and~\eqref{eq:sigma-comb-vs-top-3}
that $\gamma((-\epsilon,0)) \cap \sigma = \emptyset$, which means $x\in\sigma^+$.
This proves \eqref{eq:sigma-plus}, and~\eqref{eq:sigma-i} follows immediately.
\qed\medskip

In addition to the straightforward exit and entry set characterizations 
derived above, the notion of admissible semiflow also shows that toplexes
have to be strong isolating blocks. This is the subject of the next result.
\begin{prop}[Toplexes are strong isolating blocks]
\label{prop:small-iso-blocks}
Let~$X$ denote a cellular space with cellular decomposition~$\cX$,
and assume that~$\phi$ is an admissible semiflow with respect to~$\cX$.
Then every $\sigma\in\cXtop$ is a strong isolating block for~$\phi$
in the sense of Definition~\ref{def:strongisoblock}.
\end{prop}
\proof
Condition (A3) of admissibility assures that both~$\sigma^-$ and~$\sigma^+$
are closed, and according to their definition in~\eqref{def:exitentry}
we have $\sigma^- \cup \sigma^+ \subset \bd\sigma$. In order to see the
reverse inclusion, suppose that $x \in \bd\sigma$ and that $x \not\in \sigma^-$.
Then $\ift(x)=\sigma$ in view of~\eqref{eq:sigma-minus}. But then it follows
from~(A2) that for every backward solution~$\gamma$ through~$x$ we have to have
the inequality $\ipt(\gamma) \neq \ift(x) = \sigma$. Hence, by~\eqref{eq:sigma-i},
we finally obtain $x\in\sigma^+$. This implies the desired identity
$\sigma^- \cup \sigma^+ = \bd\sigma$, and the result follows.
\qed\medskip

Our next result shows that as long as we stay in the interior of a cell,
the immediate future of a point cannot change. Note that this result holds
for any cell, not just top-dimensional ones.
\begin{prop}[Immediate future is constant on cell interiors]
\label{prop:imm-fut}
Let~$X$ denote a cellular space with cellular decomposition~$\cX$,
and assume again that~$\phi$ is an admissible semiflow with respect to~$\cX$.
Then for every cell $\sigma\in\cX$ and any two points $x,y\in\opc{\sigma}$
we have $\ift(x)=\ift(y)$.
\end{prop}
\proof
If $\sigma\in\cXtop$, then property~\eqref{eq:sigma-comb-vs-top-1} immediately
implies $\ift(x)=\sigma=\ift(y)$. Consider now the case $\sigma\not\in\cXtop$,
and assume that $\ift(x)\neq\ift(y)$. Let $\sigma_1 := \ift(x)$ and $\sigma_2 :=
\ift(y)$. Then we have both $x\in\bd\sigma_1$ and $y\in\bd\sigma_2$. Since~$\cX$
is a cellular decomposition, we further have $\sigma\subset \sigma_1\cap\sigma_2$
due to~(CD1). Clearly, the inclusions $x\in\sigma_2^-$ and $y\in\sigma_1^-$ are
satisfied, since both $\ift(x)=\sigma_1\neq\sigma_2$ and $\ift(y)=\sigma_2\neq
\sigma_1$.

Now recall that $\opc{\sigma}$ is path connected, as the homeomorphic image of
an open ball in Euclidean space. We therefore can consider a path $\mu:[0,1]
\to \opc{\sigma}$ which satisfies both $\mu(0)=x$ and $\mu(1)=y$. By shortening the
path and potentially renaming the point~$y$, if necessary, we may assume without
loss of generality that $\setof{\ift(\mu(s))\mid s\in [0,1]} \subset
\{\sigma_1,\sigma_2\}$. Now consider the argument $s_*:=\sup\setof{s\in [0,1]\mid
\ift(\mu(s))=\sigma_1}$ and let $z:=\mu(s_*)$. Then $z\in\cl\sigma_2^- =
\sigma_2^-$, as well as $z\in\cl\sigma_1^- = \sigma_1^-$. It follows that
we then have $\ift(z)\neq \sigma_2$ and $\ift(z)\neq \sigma_1$ --- which in
turn contradicts the inclusion $\ift(z)\in \{\sigma_1,\sigma_2\}$ and proves 
the result.
\qed\medskip

With Proposition~\ref{prop:small-iso-blocks} it becomes clear that the
top-dimensional cells in the cellular decomposition~$\cX$ are natural
candidates for bricks, as defined in Definition~\ref{def:bricks}, especially
if the semiflow~$\phi$ is strongly admissible. In order to apply our 
results from Section~\ref{sec:periodicex}, we also need a mechanism for
verifying that a pair of bricks is proper. For the flow case, the
following result significantly simplifies this task.
\begin{prop}[Intersections of toplexes for admissible flows]
\label{prop:topl-intersection}
Let~$X$ denote a cellular space with cellular decomposition~$\cX$,
and assume that~$\phi$ is an admissible flow with respect to~$\cX$.
Then for every pair of cells $\sigma,\tau\in\cXtop$ such that
$\sigma \neq \tau$ we have
\begin{displaymath}
  \sigma \cap \tau \; = \;
  \left( \sigma^- \cap \tau^+ \right) \cup
    \left( \sigma^+ \cap \tau^- \right).
\end{displaymath}
\end{prop}
\proof
In view of $\tau\neq\sigma$, we deduce from
Proposition~\ref{prop:open-cells-intersection} that $\opc{\sigma}
\cap \opc{\tau} = \emptyset$. Moreover, property~\eqref{eq:sigma-comb-vs-top-3}
implies both $\opc{\tau} \cap \cebd{\sigma} = \emptyset$ and
$\opc{\sigma} \cap \cebd{\tau} = \emptyset$. Therefore, we have
\begin{displaymath}
  \sigma \cap \tau =
  (\opc{\sigma} \cup \cebd{\sigma}) \cap (\opc{\tau} \cup \cebd{\tau}) =
  \cebd{\sigma} \cap \cebd{\tau},
\end{displaymath}
and together with property~\eqref{eq:sigma-comb-vs-top-4} this further
implies
\begin{displaymath}
  \sigma\cap\tau=\bd\sigma\cap\bd\tau.
\end{displaymath}
Thus, Proposition~\ref{prop:small-iso-blocks} establishes the identities
\begin{eqnarray*}
  \sigma \cap \tau & = & \bd\sigma \cap \bd\tau \;\; = \;\;
    (\sigma^-\cup\sigma^+) \cap (\tau^-\cup\tau^+) \\
  & = & (\sigma^-\cap\tau^-)\cup(\sigma^-\cap\tau^+) \cup
    (\sigma^+\cap\tau^-)\cup(\sigma^+\cap\tau^+).
\end{eqnarray*}
By the very definition of entry and exit sets we get 
\begin{displaymath}
  \left( \sigma^- \cap \tau^+ \right) \cup
    \left( \sigma^+ \cap \tau^- \right) \; \subset \;  \sigma \cap \tau.
\end{displaymath}
Hence, to complete the proof it suffices to verify the inclusion
\begin{displaymath}
  (\sigma^-\cap\tau^-)\cup(\sigma^+\cap\tau^+)
  \;\; \subset \;\;
  (\sigma^-\cap\tau^+)\cup(\sigma^+\cap\tau^-).
\end{displaymath}
For this, suppose that the inclusion does not hold. Then there exists a
point $x\in X$ which belongs to the left-hand side, but not to the right-hand
side. Consider first the case that we have $x\in\sigma^-\cap\tau^-$. Since
$x\not\in(\sigma^-\cap\tau^+)\cup(\sigma^+\cap\tau^-)$, we then get
$x\not\in\tau^+$ and $x\not\in\sigma^+$. It follows from~\eqref{eq:sigma-i}
that $\ipt(x)=\sigma$ and $\ipt(x)=\tau$, and therefore $\sigma=\tau$, which
is a contradiction. If on the other hand $x\in\sigma^+\cap\tau^+$, an analogous
argument based on~\eqref{eq:sigma-minus} shows that $\sigma=\ift(x)=\tau$, which
gives again a contradiction.
\qed\medskip

As the above result shows, in the flow case a nontrivial intersection 
of two toplexes has to be  the union of~$\sigma^- \cap \tau^+$
and~$\sigma^+ \cap \tau^-$. We say that $\sigma$ and $\tau$ have a
\emph{circular intersection}, if their intersection $\sigma \cap \tau$ 
is non-empty and does not equal just~$\sigma^- \cap \tau^+$ or~$\sigma^+ \cap \tau^-$.
One can easily see that in general circular intersections are possible. 
However, in practical situations of interest they are very unlikely to appear
and their absence may easily be verified algorithmically.
\begin{figure}[tb]
  \begin{center}
  \includegraphics[width=0.9\textwidth]{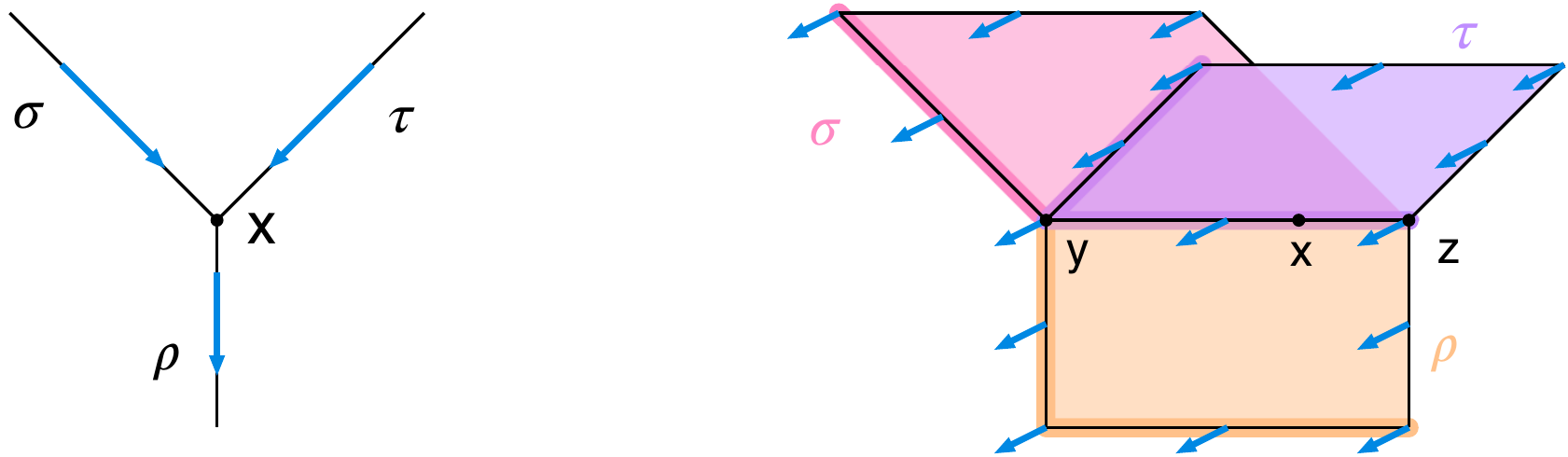}
  \caption{Counterexamples to Proposition~\ref{prop:topl-intersection} in the
           semiflow case. In both cases, the point~$x$ is not contained
           in either~$\sigma^- \cap \tau^+$ or~$\sigma^+ \cap \tau^-$.
           Note, however, that one can use Proposition~\ref{prop:topl-intersection2}
           to show that both~$(\sigma,\rho)$ and~$(\tau,\rho)$ are proper brick
           pairs. While this is immediate in the left figure, in the right image
           the point~$y$ is the unique point in~$\sigma^- \cap \rho^-$ or
           $\tau^- \cap \rho^-$, respectively, but there is no backward
           solution inside~$\rho$ in either case.}
  \label{fig:semiflowbrick}
  \end{center}
\end{figure}

Unfortunately, the situation is different in the semiflow case. As demonstrated
in the two examples shown in Figure~\ref{fig:semiflowbrick}, it is not true in
general that $\sigma \cap \tau \subset (\sigma^- \cap \tau^+) \cup (\sigma^+
\cap \tau^-)$. Nevertheless, even in the semiflow case one can fairly easily
determine whether two specifically given toplexes satisfy this inclusion or not,
since the next result identifies the exceptional situation.
\begin{prop}[Intersection of toplexes for admissible semiflows]
\label{prop:topl-intersection2}
Let~$X$ denote a cellular space with cellular decomposition~$\cX$,
and assume that~$\phi$ is an admissible semiflow with respect to~$\cX$.
Then for every pair of cells $\sigma,\tau\in\cXtop$ such that
$\sigma \neq \tau$ exactly one of the following two alternatives holds:
\begin{itemize}
\item[(i)] The equality $\sigma \cap \tau = (\sigma^- \cap \tau^+)
\cup (\sigma^+ \cap \tau^-)$ is satisfied.
\item[(ii)] There exists a point~$x \in \sigma^- \cap \tau^-$ and a toplex
$\rho \in \cXtop$ which is different from~$\sigma$ and~$\tau$ such that
$\ift(x) = \rho$, and there are backward solutions $\gamma, \eta \in
\Sol^-(x,\phi)$ such that both $\gamma((-\epsilon,0)) \subset \sigma$
and $\eta((-\epsilon,0)) \subset \tau$ are satisfied for some $\epsilon > 0$.
\end{itemize}
\end{prop}
\proof
As in the proof of Proposition~\ref{prop:topl-intersection} one can show that
\begin{displaymath}
  \sigma \cap \tau \;\; = \;\;
  (\sigma^-\cap\tau^-)\cup(\sigma^-\cap\tau^+) \cup
    (\sigma^+\cap\tau^-)\cup(\sigma^+\cap\tau^+),
\end{displaymath}
and the same argument used in that proof can still be used to establish
\begin{displaymath}
  \sigma^+\cap\tau^+
  \;\; \subset \;\;
  (\sigma^-\cap\tau^+)\cup(\sigma^+\cap\tau^-).
\end{displaymath}
Thus, the only way that the inclusion in~(i) cannot be satisfied is if there
exists a point $x\in X$ which belongs to $(\sigma^- \cap \tau^-) \setminus
((\sigma^-\cap\tau^+)\cup(\sigma^+\cap\tau^-))$. But such a point has to 
have the properties stated in~(ii).
\qed\medskip
%

%%%%%%%%%%%%%%%%%%%%%%%%%%%%%%%%%%%%%%%%%%%%%%%%%%%%%%%%%%%%%%%%%
%%%%%%%%%%%%%%%%%%%%%%%%%%%%%%%%%%%%%%%%%%%%%%%%%%%%%%%%%%%%%%%%%
\subsection{Combinatorial multivector fields induced by admissible semiflows}
\label{subsec53}

In this subsection we demonstrate that every admissible semiflow
on a cellular space gives rise to a well-defined combinatorial multivector
field on~$\cX$. This construction will be used in the next subsection to show
that, conversely, periodic orbits in the combinatorial setting give rise
to corresponding periodic orbits in the underlying semiflow under suitable
conditions.

To set the stage, we recall some basic definitions from~\cite{lipinski:etal:p20a}.
Suppose we are given a cellular space~$X$ with cellular decomposition~$\cX$, which
will be considered as a finite Alexandrov space. Then a \emph{combinatorial
multivector field} on the finite topological space~$\cX$ is a partition~$\cV$
of~$\cX$ into locally closed, non-empty subsets which we call \emph{multivectors}.
A multivector $V\in\cV$ is $\emph{regular}$ if $H_*(\Cl V,\Mo V)=0$, otherwise it is \emph{critical}.
This distinction is motivated by the Wa\.zewski Theorem (Proposition~\ref{prop:block}).
A subset $\cA\subset\cX$ is \emph{$\cV$-compatible} if $V\cap\cS\neq\emptyset$
for some multivector $V\in\cV$ implies $V\subset\cS$. Note that given a $\cV$-compatible
subset  $\cA\subset\cX$, there is an induced multivector field 
$\cV_\cA:=\setof{V\in\cV\mid \cV\subset\cA}$ on $\cA$.

With a multivector field $\cV$ we associate the \emph{combinatorial dynamical system}
generated by the multivalued map $F_\cV : \cX \mto \cX$ whose value for the cell
$\sigma\in\cX$ is defined as
\begin{equation} \label{eq:F-cV}
  F_\cV(\sigma) := \Cl\sigma\cap[\sigma]_\cV.
\end{equation}
In this formula, $\Cl\sigma$ denotes the closure of~$\sigma$ in the finite
topological space~$\cX$, and~$[\sigma]_\cV$ denotes the unique multivector
in~$\cV$ which contains~$\sigma$.

The combinatorial dynamics associated with the combinatorial multivector
field~$\cV$ is given by the dynamics of the multivalued map~$F_\cV$. Thus,
a partial map $\gamma:\ZZ\pto\cX$ is called a \emph{solution} of $F_\cV :
\cX \mto \cX$, if the inclusion $\gamma(i+1) \in F_\cV(\gamma(i))$ is
satisfied for all $i,i+1 \in \dom\gamma$. A solution~$\gamma$ \emph{through
$\sigma\in\cX$} is a solution such that $0 \in \dom\gamma$ and $\gamma(0) =
\sigma$. A solution~$\gamma$ \emph{in $\cA\subset\cX$} is a solution with
values in the set~$\cA$. As in the classical case, a solution~$\gamma$ is called
\emph{full}, if we have $\dom\gamma=\ZZ$. Finally, a solution~$\gamma$ is
called a \emph{path}, if its domain $\dom\gamma$ is a bounded interval
in~$\ZZ$. In this case, the \emph{endpoints} of the path are the values
of the path at the endpoints of the domain.  Formula  \eqref{eq:F-cV}
for the multivalued map ~$F_\cV$ induced by the multivector field $\cV$
has one drawback: Every cell $\sigma\in\cX$  admits a stationary (constant)
solution, because one has $\sigma\in\Cl\sigma$. To overcome this, we say that 
a full solution $\gamma:\ZZ\to\cX$ is  \emph{essential}, if 
for every $t\in\ZZ$ for which the multivector $V:=[\gamma(t)]_\cV$ is regular,
there exist positive offsets $k,k'>0$ such that $\gamma(t+k)$ and $\gamma(t-k')$
are not contained in~$V$.

We say that a full solution   $\rho:\ZZ\to X$  is \emph{recurrent}, if
$\rho^{-1}(x)$ is unbounded from both the left and the right for every
$x\in \im\rho$. We say that a multivector field $\cV$ is \emph{gradient-like},
if for every recurrent solution $\rho:\ZZ\to X$ there exists a $V\in\cV$
such that $\im\rho\subset V$.

The combinatorial setting also has notions of invariance. We first define
the \emph{invariant part} of a subset $\cN\subset\cX$, denoted by $\Inv\cN$, 
as the set of cells $\sigma\in\cN$ such that there exists an essential
solution through~$\sigma$ in~$\cN$.
We say that a subset $\cA\subset\cX$ is  \emph{invariant}, if $\cA=\Inv \cA$.
 A closed
set~$\cN\subset\cX$ is called an \emph{isolating set} for an invariant
set~$\cS$, if every path in~$\cN$ with endpoints in~$\cS$ is a solution
in~$\cS$. If such an isolating set exists, then the set~$\cS$ is called
an \emph{isolated invariant set}.

For later reference, we recall the following two properties of isolated
invariant sets in the combinatorial setting. Their proof can be found
in~\cite[Propositions~4.10 and~4.12]{lipinski:etal:p20a}.
\begin{prop}[Properties of isolated invariant sets]
\label{prop:comb-iso-inv-sets}
Assume that~$\cS$ is an isolated invariant set for a combinatorial multivector
field~$\cV$. Then
\begin{itemize}
   \item[(i)]  $\cS$ is locally closed, and
   \item[(ii)]   $\cS$ is $\cV$-compatible.  \qed
\end{itemize}
\end{prop}
The above notions form the foundation for a dynamical theory of
combinatorial multivector fields on finite topological spaces, which
includes concepts such as index pair, the Conley index, and Morse decompositions.
As in the classical case, the homological Conley index is defined as the homology 
of an index pair. However, unlike the classical case, there is a very special index
pair, namely  $(\Cl\cS,\Mo\cS)$. This leads to the following proposition.
\begin{prop}[Conley index of an isolated invariant set]
\label{prop:Conley-index-combimnatorial}
The Conley index of an isolated invariant set $\cS$ in a combinatorial
vector field $\cV$ is $H_*(\Cl\cS,\Mo\cS)$.
\qed
\end{prop}
For more details, we refer the reader again to~\cite{lipinski:etal:p20a}.
For us, however, combinatorial multivector fields will be a useful
language for determining periodic orbits in semiflows, more precisely,
admissible semiflows on the cellular space~$X$. This in turn makes it
necessary to have a canonical way of creating a combinatorial multivector
field from a given admissible semiflow --- and this will be described 
in the remainder of this subsection.

Suppose now that~$\phi$ is an admissible semiflow on the cellular
space~$X$ with cellular decomposition~$\cX$. Then we say that~$\sigma
\in \cX$ is the \emph{immediate future} of a cell $\tau\in \cX$, if we
have $\ift(x) = \sigma$ for every $x \in \opc{\tau}$. In this case, we
write $\ift(\tau) = \sigma$. Notice that in view of
Proposition~\ref{prop:imm-fut}, the immediate future~$\ift(\tau)$
is well-defined for every $\tau \in \cX$. With this concept, we can 
now associate a combinatorial multivector field on~$\cX$ to the
semiflow~$\phi$.
\begin{prop}[Combinatorial multivector field generated by an admissible semiflow]
\label{prop:cmvf}
Let~$X$ denote a cellular space with cellular decomposition~$\cX$,
and assume that~$\phi$ is an admissible semiflow with respect to~$\cX$.
For every $\sigma \in \cXtop$ we define
\begin{displaymath}
  V_\sigma := \setof{ \tau \in \cX \mid \ift(\tau) = \sigma}.
\end{displaymath}
Then the collection
\begin{displaymath}
  \cV := \setof{V_\sigma \mid \sigma \in \cXtop}
\end{displaymath}
is a combinatorial multivector field on~$\cX$. We call it the \emph{combinatorial
multivector field~$\cV$ associated with the semiflow~$\phi$}.
\end{prop}
\proof
We begin by showing that~$\cV$ is a partition of the cellular
decomposition~$\cX$. To see that the elements of~$\cV$ are mutually
disjoint, assume that $\sigma,\sigma'\in\cXtop$ are toplexes such
that there exists a cell $\tau \in V_\sigma \cap V_{\sigma'}$. In
view of the definition of the immediate future of a cell this implies
$\sigma = \ift(\tau) = \sigma'$, which in turn proves that
$V_\sigma = V_{\sigma'}$. Furthermore, it is clear from the definitions
that every $\tau \in \cX$ belongs to~$V_{\ift{\tau}}$. Hence, one has
$\cX = \bigcup_{\sigma\in\cXtop} V_\sigma$, and this shows
that~$\cV$ is a partition of~$\cX$.

In order to establish~$\cV$ as a combinatorial multivector field, we
still need to verify that the sets~$V_\sigma$ are locally closed for
all $\sigma\in\cXtop$. For this, we first verify that
\begin{equation} \label{eq:cmvf-1}
  V_\sigma\subset\Cl\sigma.
\end{equation}
Let $\tau\in V_\sigma$. Then, by~\eqref{eq:ift-x-sigma}, for any
$x\in\opc{\tau}$ we have $x\in\sigma$. Hence, assumption~(CD1)
implies that $\tau\subset\sigma$. According to the definition of
the Alexandrov topology this implies $\tau\in\Cl\sigma$, and
establishes the validity of~\eqref{eq:cmvf-1}.

Using~\eqref{eq:cmvf-1} one immediately obtains $\Cl V_\sigma
\subset \Cl\sigma$. Moreover, due to the definition of the immediate
future we always have $\sigma\in V_\sigma$, which in turn implies the
inclusion $\Cl\sigma \subset \Cl V_\sigma$. It follows that
$\Cl V_\sigma = \Cl\sigma$, as well as $\Mo V_\sigma = \Cl\sigma
\setminus V_\sigma$. We will now verify that
\begin{equation} \label{eq:cmvf-2}
  |\Mo V_\sigma| =
  \setof{x\in\sigma \mid \ift(x)\neq\sigma} =
  \setof{x\in\bd\sigma \mid \ift(x)\neq\sigma}.
\end{equation}
To this end, consider a point $x\in X$ and let $\tau=\cell x$.
Assume first that $x \in |\Mo V_\sigma|$. Then we have both
$\tau \in \Cl\sigma$ and $\tau\not \in V_\sigma$.
It follows that $\tau\subset\sigma$ and $\ift(\tau)\neq\sigma$,
which implies that~$x$ belongs to the middle set of~\eqref{eq:cmvf-2}.
Assume now, conversely, that the point~$x$ belongs to the middle 
set of~\eqref{eq:cmvf-2}. Then $x \in \opc{\tau} \cap \sigma$, and
we obtain from property~(CD1) of a cellular decomposition that
$\tau \subset \sigma$. Hence, we have $\tau\in\Cl\sigma$, and
since $\ift(\tau) = \ift(x) \neq \sigma$, we also get
$\tau \not \in V_\sigma$. In consequence, $\tau\in\Mo V_\sigma$
and $x\in |\Mo V_\sigma|$. This proves \eqref{eq:cmvf-2}, since
the second equality in this equation is immediate.

If we finally combine~\eqref{eq:cmvf-2} with assumption~(A3) in
Definition~\ref{def:admissible} and with Proposition~\ref{prop:sigma-pm},
then one immediately sees that the set~$|\Mo V_\sigma|$ is closed in~$X$.
Therefore, by Corollary~\ref{cor:rep-opn-cls}, the set~$\Mo V_\sigma$
is closed in the finite topological space~$\cX$, which means
that~$V_\sigma$ is locally closed. This completes the proof of
the proposition.
\qed\medskip
\begin{figure}[tb]
  \begin{center}
  \includegraphics[width=0.9\textwidth]{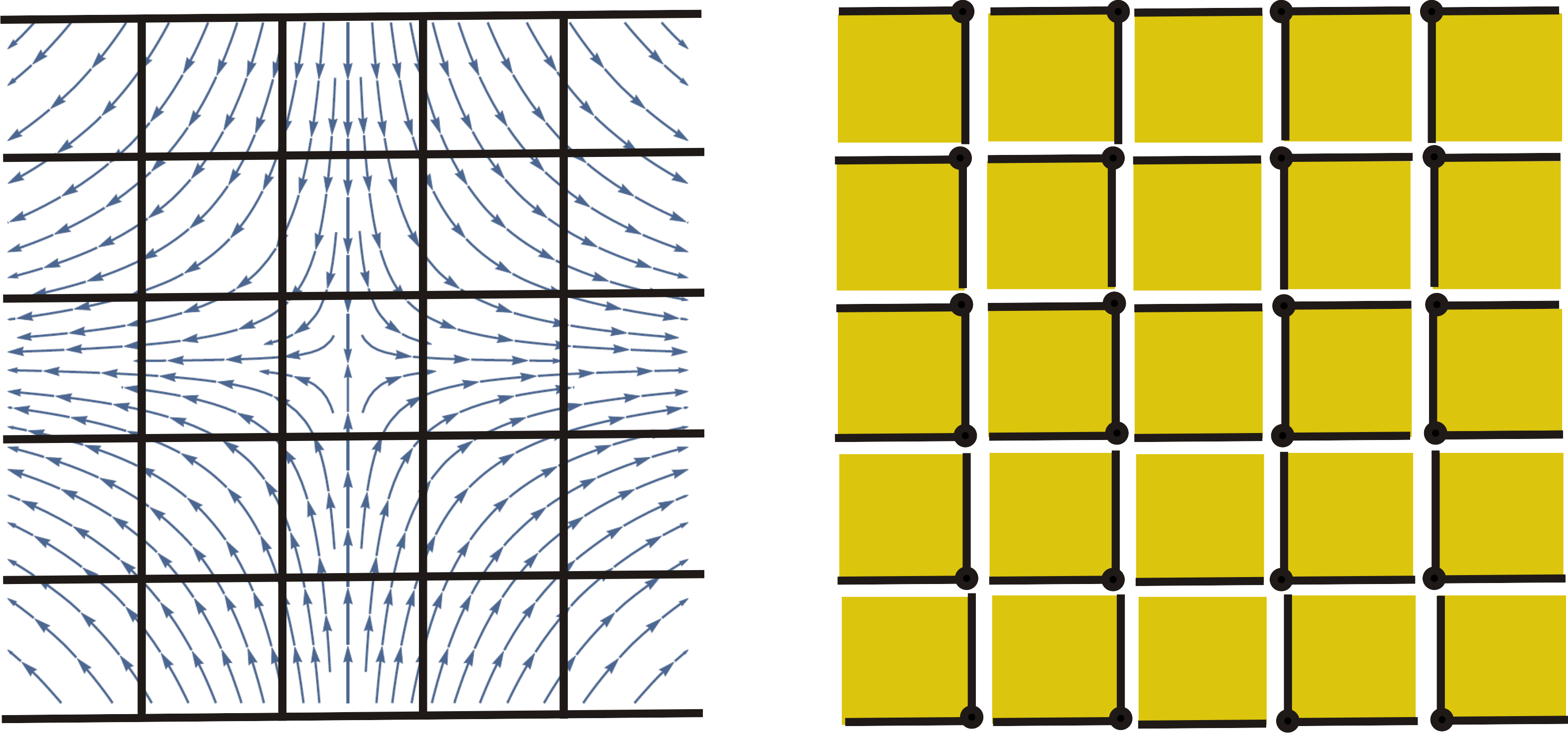}
  \caption{Left: A fragment of cellular space $X$ with cellular decomposition~$\cX$ consisting
  of squares as toplexes together with an admissible flow. Right: The associated multivector field~$\cV$.  }
  \label{fig:saddleWith-cmvf}
  \end{center}
\end{figure}

To illustrate Proposition~\ref{prop:cmvf} consider the fragment of cellular space $X$ presented in the left part of
Figure~\ref{fig:saddleWith-cmvf}. Its cellular decomposition consists of squares as toplexes, edges, and vertices. 
A sample admissible flow, representing a saddle, is visualized in the figure by a collection of its trajectories.
The resulting multivector field is shown on the right. All its multivectors are regular except the one in the middle, 
which is critical. Note that the unique critical multivector covers the unique stationary point of the flow. For all 
other multivectors the trajectories pass through them.

In the sequel, in the setting of the combinatorial multivector
field~$\cV_\phi$, by a \emph{regular} toplex we mean a toplex
$\sigma\in\cXtop$ such that the associated
multivector~$V_\sigma$ is regular.

%%%%%%%%%%%%%%%%%%%%%%%%%%%%%%%%%%%%%%%%%%%%%%%%%%%%%%%%%%%%%%%%%
%%%%%%%%%%%%%%%%%%%%%%%%%%%%%%%%%%%%%%%%%%%%%%%%%%%%%%%%%%%%%%%%%
\subsection{Isolating blocks for admissible semiflows via combinatorial dynamics}
\label{subsec54}

We begin with two straightforward propositions concerning the multivector field~$\cV$
associated with an admissible semiflow~$\phi$.
\begin{prop}[Multivector closures in the induced multivector field]
\label{prop:cl-V-sigma}
  Let~$X$ denote a cellular space with cellular decomposition~$\cX$,
  and assume that~$\phi$ is an admissible semiflow with respect to~$\cX$.
  Furthermore, let~$\cV$ denote the associated multivector field
  introduced in Proposition~\ref{prop:cmvf}. Then for every
  $\sigma\in\cXtop$ we have
  \begin{displaymath}
    \Cl V_\sigma = \Cl\sigma .
  \end{displaymath}
\end{prop}
\proof
   Consider a cell $\tau\in\Cl V_\sigma$. Then we have $\tau\subset\rho$ for some
   $\rho\in V_\sigma$ such that $\ift(\rho)=\sigma$. In view of~\eqref{eq:cmvf-1}
   one obtains $\rho\subset\sigma$, which in turn implies $\tau\subset\rho$ and
   gives $\tau\in\Cl\sigma$. This proves the inclusion $\Cl V_\sigma\subset\Cl\sigma$.
   In order to see the opposite inclusion, observe that $\sigma\in V_\sigma$, which
   immediately yields $\Cl\sigma\subset\Cl V_\sigma$.
\qed
\begin{prop}[Representation of $\cV$-compatible sets]
\label{prop:supp-cl-S}
  Let~$X$ denote a cellular space with cellular decomposition~$\cX$,
  and assume that~$\phi$ is an admissible semiflow with respect to~$\cX$.
  Furthermore, let~$\cV$ denote the associated multivector field
  introduced in Proposition~\ref{prop:cmvf}, and assume that $\cS\subset\cX$ is
  $\cV$-compatible. Then we have
  \begin{displaymath}
    |\Cl\cS|=\bigcup_{\sigma\in\cS\cap\cXtop}\sigma.
  \end{displaymath}
\end{prop}
\proof
Since $\cS$ is  $\cV$-compatible, we have 
\begin{displaymath}
 \cS=\bigcup_{V_\sigma\subset\cS} V_\sigma.
\end{displaymath}
Therefore, by Proposition~\ref{prop:cl-V-sigma}
\begin{displaymath}
   \Cl\cS=\bigcup_{V_\sigma\subset\cS} \Cl V_\sigma=\bigcup_{V_\sigma\subset\cS} \Cl\sigma,
\end{displaymath}
and Proposition~\ref{prop:supp-comb-cl-cell} further implies
\begin{displaymath}
   |\Cl\cS|=\bigcup_{V_\sigma\subset\cS} |\Cl \sigma|=\bigcup_{\sigma\in\cS\cap\cXtop}\sigma,
\end{displaymath}
because clearly one has $V_\sigma\subset\cS$ if and only if $\sigma\in\cS\cap\cXtop$.
\qed\medskip

The next proposition shows that every isolated invariant set of a multivector
field~$\cV$ associated with an admissible flow~$\phi$ induces an isolating block
for~$\phi$ with the same Conley index as~$\cV$.
\begin{thm}[Isolating blocks and Conley index via combinatorics]
\label{thm:iso-block-from-cmvf}
  Let~$X$ denote a cellular space with cellular decomposition~$\cX$,
  and assume that~$\phi$ is an admissible flow with respect to~$\cX$.
  Furthermore, let~$\cV$ denote the associated multivector field
  introduced in Proposition~\ref{prop:cmvf}. Finally, assume that~$\cS$
  is an isolated invariant set for~$\cV$ and let $N:=|\Cl \cS|$. Then
  the following holds:
  \begin{itemize}
  \item[(i)]   The exit set of~$N$ is given by $N^-=|\Mo \cS|$.
  \item[(ii)]  The set~$N$ is an isolating block for~$\phi$.
  \item[(iii)] The Conley indices of~$\cS$ and~$\Inv(N,\phi)$ coincide.
  \end{itemize}
\end{thm}
\proof
It follows from Proposition~\ref{prop:supp-cl-S} that~$N$ is compact as a
finite union of compact sets. To prove~(i), consider first an $x\in N^-$.
Then $x\in N$, which implies that $\tau:=\cell{x}\in\Cl\cS$. Hence, there is
a $\rho\in\cS$ such that $\tau\subset\rho$. Let $\sigma\in\cXtop$ be such
that $\rho\in V_\sigma$. Since~$\cS$ is $\cV$-compatible by
Proposition~\ref{prop:comb-iso-inv-sets}(ii) and $\rho\in\cS$, we see that
also $\sigma\in\cS$. Let $\omega:=\ift(x)=\ift(\tau)$. Then one has
$\omega\not\in\cS$, since otherwise $x\not\in N^-$. From the definition of
the multivector~$V_\omega$ we get $\tau\in V_\omega$. Since $\omega\not\in\cS$
and $\omega\in V_\omega$, again by $\cV$-compatibility of~$\cS$ one obtains
$\tau\not\in\cS$. Therefore, $\tau\in\Mo\cS$ and $x\in |\Mo \cS|$. This establishes
$N^-\subset|\Mo \cS|$. To see the opposite inclusion, consider an $x\in |\Mo \cS|$.
Let $\tau:=\cell{x}$. Then $\tau\in\Cl \cS\setminus\cS$. Hence, $\tau\not\in\cXtop$,
and therefore $x\in\Fr(X)$. Let $\omega:=\ift(x)$. Since $\tau\not\in\cS$, and
$\tau,\omega\in V_\omega$, we get $\omega\not\in\cS$. It follows that
$\opc{\omega}\cap N=\emptyset$, which finally implies $x\in N^-$.

We now turn our attention to~(ii). In view of~(i) and
Corollary~\ref{cor:rep-opn-cls} the exit set~$N^-$ is closed. Thus,
Proposition~\ref{prop:strongexit}(b) immediately implies~(ii), as long
as we can decompose the boundary of~$N$ in the form $\bd N = N^+ \cup N^-$.
Assume that this is false. Then there exists a point $x \in \bd N \setminus
(N^+ \cup N^-)$. Let $\sigma:=\ift(x)$, $\tau:=\ipt(x)$, and $\omega:=\cell x$.
We have $\sigma\in\cS$, because otherwise $x\in N^-$. Similarly, one has to have
$\tau\in\cS$, since otherwise $x\in N^+$. Due to the inclusion $x \in \opc{\omega}
\cap \tau \cap \sigma$, we further get $\omega\subset\tau\cap\sigma$. Moreover,
since $\ift(\omega)=\ift(x)=\sigma$, we see that $\omega\subset V_\sigma\subset\sigma$.
Hence, the $\cV$-compatibility of~$\cS$ implies that $\omega\in\cS$, because
$\sigma\in\cS$. The set $U:=|\Star\omega|$ is open in~$X$ according to
Proposition~\ref{prop:star-open}. Thus, the continuity of~$\phi$ shows that
we can find an open neighborhood~$V$ of~$x$ and an $\epsilon>0$ such that
$\phi(V,[-\epsilon,\epsilon])\subset U$. By possibly decreasing the size of
the neighborhood~$V$ further, we can also ensure that the inclusion
$\phi(V,-\epsilon)\subset\opc{\tau}$ is satisfied, in view of
$\ipt(x)=\tau$. Due to $x\in\bd N$, we can select a point $v\in V\setminus N$.
Without loss of generality we can further assume that $\mu:=\cell{v}\in\cXtop$,
since also some neighborhood of~$v$ is disjoint from~$N$. But then we have to have
$\mu\not\in \cS$, because otherwise $v\in N$. According to our construction, one
also has $\phi(v,-\epsilon)\in\opc{\tau}\subset N$. Therefore, there has to be a
point $z\in  \phi(v,(-\epsilon,0))\cap\bd N$ such that $\ift(z)=\mu$. Let
$\nu:=\cell{z}$. Then we have the inclusion $\nu\in V_\mu$, and since $\mu\not\in\cS$,
the $\cV$-compatibility of~$\cS$ implies $\nu\not\in\cS$. In addition, 
Proposition~\ref{prop:supp-cl-S}, in combination with $z\in\bd N\subset N$,
shows that we have $z\in\rho$ for some $\rho\in\cS\cap\cXtop$. It follows
from (CD1) that $\nu\subset\rho$. Since we also have
\begin{displaymath}
  z \in \phi(v,(-\epsilon,0)) \subset
  \phi(V,[-\epsilon,\epsilon]) \subset U = |\Star\omega|,
\end{displaymath}
we immediately obtain $\nu\in\Star\omega$, and therefore $\omega\subset\nu$.
Since $\omega,\rho\in\cS$ and~$\cS$ is locally closed by 
Proposition~\ref{prop:comb-iso-inv-sets}(i), we finally deduce $\nu\in\cS$,
which is a contradiction. Altogether, we have therefore established the
equality $\bd N = N^+ \cup N^-$, as well as~(ii).

In order to complete the proof of the theorem, we still have to verify the
validity of~(iii). By Proposition~\ref{prop:Sigma-iso} we have the isomorphism
\begin{displaymath}
  H(\Sigma_{|\Cl\cS|}): H(|\Cl\cS|,|\Mo\cS|)\to  H(\Cl\cS,\Mo\cS) .
\end{displaymath}
In addition, $H(|\Cl\cS|,|\Mo\cS|) = H(N,N^-)$ is the Conley index of the
isolated invariant set~$\Inv(N,\phi)$ due to Proposition~\ref{prop:Conley-index},
and $H(\Cl\cS,\Mo\cS)$ is the Conley index of~$\cS$ in view of
Proposition~\ref{prop:Conley-index-combimnatorial}. This completes
the proof.
\qed

%%%%%%%%%%%%%%%%%%%%%%%%%%%%%%%%%%%%%%%%%%%%%%%%%%%%%%%%%%%%%%%%%
%%%%%%%%%%%%%%%%%%%%%%%%%%%%%%%%%%%%%%%%%%%%%%%%%%%%%%%%%%%%%%%%%
\subsection{Periodic solutions via combinatorial Poincar\'e sections}
\label{subsec55}

Let~$X$ denote a cellular space with cellular decomposition~$\cX$, and
assume that~$\phi$ is an admissible flow with respect to~$\cX$.
Theorem~\ref{thm:iso-block-from-cmvf} facilitates the search for interesting
isolated invariant sets for~$\phi$ via an algorithmic study of the associated
multivector field~$\cV_\phi$. Note that due to Proposition~\ref{prop:small-iso-blocks}
every toplex $\sigma\in\cXtop$ is a strong isolating block. Hence, if~$\sigma$ is also
regular, then~$\sigma$ is a brick in the sense of Definition~\ref{def:bricks}.
Therefore, it is natural to expect that the assumptions of Theorem~\ref{thm:bricks}
may in some cases be verified algorithmically. We first make an observation on a
correspondence between certain brick paths and certain solutions of~$\cV_\phi$.
This requires a few more concepts.

We say that a proper pair of bricks~$(B,\tilde{B})$ is \emph{sharp}, if there
exists a point $x\in B\cap\tilde{B}$ and a time $t>0$ such that $\phi(x,[0,t])
\subset \tilde{B}$. Moreover, we say that a brick path is \emph{sharp}, if each
two consecutive bricks in the path form a sharp pair. It is not difficult to verify
that assumption~(f) in Theorem~\ref{thm:bricks} may be weakened by considering only
sharp brick paths. This matters in the context of strongly admissible flows, because
we have the following proposition.
\begin{prop}[Sharp brick pairs and combinatorial paths]
\label{prop:sharp-pair}
  Let~$X$ denote a cellular space with cellular decomposition~$\cX$,
  and assume that~$\phi$ is a strongly admissible semiflow with respect to~$\cX$.
  Furthermore, let~$\cV$ denote the associated multivector field
  introduced in Proposition~\ref{prop:cmvf}. Finally, assume that
  $\tau,\sigma\in\cXtop$ are regular toplexes. Then they are bricks
  for~$\phi$. Moreover, if $(\tau,\sigma)$ is a sharp pair of bricks,
  then there is an $\alpha\in\cX$ such that $\tau,\alpha,\sigma$ is a
  path in~$\cV$. The reverse implication holds true if additionally
  the intersection of~$\tau$ and~$\sigma$ is not circular.
\end{prop}
\proof
It follows from Proposition~\ref{prop:small-iso-blocks} and the definition
of strong admissibility that~$\tau$ and~$\sigma$ are bricks. Suppose now
that~$(\tau,\sigma)$ is a sharp pair of bricks. Let $x\in\tau^-\cap\sigma^+$
and choose the time $t>0$ such that $\phi(x,(0,t)) \subset \sigma$.
Let $\rho:=\ift(x)$. Then we have $\opc{\rho}\cap\sigma\neq\emptyset$.
This in turn implies $\rho\subset\sigma$, and since~$\rho$ is also a toplex,
we immediately obtain $\rho=\sigma$. Now consider the cell $\alpha:=\cell x$.
Since $x\in\tau^-$, we get $\alpha\subset\tau^-\subset\tau$ which yields
$\alpha\in\Cl\tau\subset F_\cV(\tau)$. We also have $\sigma = \rho \in
[\alpha]_\cV\subset F_\cV(\alpha)$, and this proves that $\tau,\alpha,\sigma$
is indeed a path in $\cV$.

Assume now that $\tau,\alpha,\sigma$ is a path in~$\cV$, and that the intersection
of~$\tau$ and~$\sigma$ is not circular. Then $(\tau,\sigma)$ is a proper pair of
bricks. To see that it is even a sharp pair, we observe that in view of
$\sigma\in F_\cV(\alpha)$ one has the inclusion $\sigma\in[\alpha]_\cV$,
because~$\sigma$ is a toplex. Therefore, $\alpha\in V_\sigma$. Now select
an $x\in\opc{\alpha}$. Then $\ift(x)=\ift(\alpha)=\sigma$, which proves that
the pair $(\tau,\sigma)$  is sharp. 
\qed\medskip

Proposition~\ref{prop:sharp-pair} has a straightforward extension to sharp
paths of bricks. To formulate it, we need one more final definition. We say
that a solution~$\gamma$ of the multivector field~$\cV = \cV_\phi$ is
\emph{reduced}, if for $i,i+1,i+2\in\dom\gamma$ we have that $\gamma(i+1)
\in \Cl\gamma(i)$ implies $\gamma(i+2) \in [\gamma(i+1)]_\cV$, and
$\gamma(i+1) \in [\gamma(i)]_\cV$ implies $\gamma(i+2) \in \Cl\gamma(i+1)$.
Then the following result is an immediate consequence of the one above.
\begin{prop}[Correspondence between brick paths and combinatorial solutions]
\label{prop:sharp-path}
In the situation of Proposition~\ref{prop:sharp-pair},
assume that $\cA\subset\cX$ is a $\cV$-compatible collection, which contains
only regular multivectors and is free of circular pairs. Then there is a
one-to-one correspondence between sharp paths of bricks in~$|\cA|$ and
combinatorial solutions of~$\cV$ in~$\cA$ in the following sense. For every
sharp path of bricks in~$|\cA|$ there exists a reduced solution of~$\cV$
in~$\cA$ such that the sequence of toplexes along the solution is exactly
the path of bricks, and the sequence of toplexes along a reduced solution
of~$\cV$ in~$\cA$ is a sharp path of bricks in~$|\cA|$.
\qed
\end{prop}
We now turn our attention to establishing the assumptions of the main
Theorem~\ref{thm:bricks}. Clearly, there are many ways in which the brick
decomposition~$\cA$ of~$A$ can be decomposed into subcollections~$\cA_i$.
Note, however, that at the very least one has to make sure that none
of the collections~$\cA_i$ contain a closed brick path. Hence, by
Proposition~\ref{prop:sharp-path}, it cannot contain a periodic solution
of~$\cV$. Thus, any choice of brick decomposition considered in
Theorem~\ref{thm:bricks} should start with the selection of a family~$\cA_0$
which breaks all periodic solutions of~$\cV$. From the point of view of pure
graph theory this is the feedback vertex set problem~\cite{karp:72a} in the
directed graph whose vertices are cells in~$\cX$ and there is a directed edge
from $\tau\in\cX$ to $\sigma\in\cX$ whenever $\sigma\in F_\cV(x)$. Although
the problem of constructing a minimal feedback vertex set is NP-complete in
general~\cite{karp:72a}, in our setting it is not too difficult, because one
can take as a candidate for~$\cA$ the minimal locally closed, $\cV$-compatible superset of
the support of a nontrivial singular cohomology cocycle in~$\Cl\cS$. To
construct such a cocycle it suffices to find it in a small retract of~$\Cl\cS$.
Typically, both finding a small retract and transferring the cocycle in the
retract back to the original space may be achieved in linear time.
This suggests the following procedure which, if successful, leads to an
appropriate coarsening of the brick decomposition of the isolating neighborhood
considered in Theorem~\ref{thm:bricks}.

Assume that~$\cA$ is a connected isolated invariant set for a combinatorial
multivector field~$\cV$ associated with a strongly admissible flow~$\phi$.
Furthermore, consider a non-empty, locally closed, and $\cV$-compatible
subset $\cP\subset\cA$. Set $\cH:=\Cl(\cA\cap\Mo\cP)$, let $\cR:=\Cl\cA\setminus\cH$,
and define $\cPb:=\Cl_\cR\cP$. Suppose further that
\begin{align}
\label{eq:P-sect-1}
&\text{$\cH$ is non-empty and $\cR$ is connected,} \\
\label{eq:P-sect-2}
&\text{Every maximal solution in~$\Cl\cA$ reaches $\cH\cup\Mo\cA$, and}\\ 
\label{eq:bd-cpb}
&\cH\cap\Cl\Bd_\cR\cPb=\emptyset.
\end{align}
We then say that~$\cP$ is a \emph{combinatorial Poincar\'e section} for~$\cA$. 

For a pair of cells $\sigma,\tau\in\cR$, denote by
$d(\sigma,\tau)$ the length of the shortest  fence
in~$\cR$ joining~$\tau$ and~$\sigma$. 
Such a fence exists, because of assumption  \eqref{eq:P-sect-1}.
For a cell $\sigma\in\cR$ we further
set $d(\sigma,\cPb):=\min\setof{d(\sigma,\tau)\mid \tau\in \cPb}$, while for a
path~$\gamma$ in~$\cR$ we define
\begin{displaymath}
  L(\gamma):=\max\setof{d(\sigma,\cPb)\mid \sigma\in\im\gamma} \; ,
\end{displaymath}
and for a cell $\tau\in\cR$ we set
\begin{displaymath}
  L(\tau):=\max\setof{L(\gamma)\mid \gamma
    \text{ is a path in $\cR $ originating in $\tau$}}.
\end{displaymath}
Then $L:\cR \to\NN$ is a function which decreases along solutions of~$\cV$
in~$\cR$, i.e., the function~$L$ is a Lyapunov function. More precisely, we
have the following straightforward proposition.
\begin{prop}[$L$ is a Lyapunov function]
\label{prop:L}
For any pair of cells $\sigma,\tau\in\cR$ with $\sigma\in F_\cV(\tau)$
we have $L(\sigma) \leq L(\tau)$.
\qed
\end{prop}
Since $\cH\neq\emptyset$ by assumption \eqref{eq:P-sect-1}, also
$\cA\cap\Mo\cP\neq\emptyset$. Therefore, since both $\cP$  and $\cA$ are
$\cV$-compatible, we see that
\begin{displaymath}
  \bar{n}:=\min\setof{L(\tau)\mid \tau\in\Opn\cH\setminus\Cl\cP}
\end{displaymath}
is well-defined. Clearly, since $L(\tau)>0$ for all $\tau\not\in\cPb$,
we also have $\bar{n}>0$.

Using the Lyapunov function~$L$, we can now recursively define a finite
sequence~$(\cL_k)_{i=0}^{\kmax}$ of subsets of~$\cR$ as follows. First,
let $\cL_0:=\emptyset$, as well as $\cL_1:=\cPb$. Now consider $1 < k < \kmax$,
and assume that  the collection~$\cL_{k-1}\subset\cR$ has already been
defined and satisfies $\cL_{k-1}\neq\cR$. 
Then
\begin{displaymath}
  \cB_k := \setof{\sigma\in\cR\setminus\cL_{k-1} \mid
    \cL_{k-1}\cap \Cl\sigma\neq\emptyset}
\end{displaymath}
is non-empty, and $n_k := \max\setof{L(\sigma)\mid \sigma\in\cB_k}$ is well-defined.
Now, if we have $n_k\geq\bar{n}$, then we define $\cL_k:=\cR$, let $\kmax:=k$, and
stop the recursion. Otherwise $L^{-1}([0,n_k])\neq\cR$ and we continue recursion
by setting
\begin{displaymath}
  \cL_k:=L^{-1}([0,n_k]).
\end{displaymath}
A straightforward argument based on Proposition~\ref{prop:L} can be
used to prove the following proposition.
\begin{prop}[Properties of the sets~$\cL_k$]
\label{prop:cL}
  For every $0 \le k < \kmax$  the above-defined set~$\cL_k$ is closed in~$\cR$
  and $\cV$-compatible.
\qed
\end{prop}
Clearly, we have the inclusions $\cB_k\subset \cL_k\setminus \cL_{k-1}$ and
$\cL_{k-1} \subset \cL_k$. From the definitions of~$\bar{n}$ and~$\kmax$
one further obtains
\begin{equation} \label{eq:kmax-1}
  \cL_{\kmax-1}\cap\Opn\cH\subset\Cl\cP.
\end{equation}
For $i=0,1,\ldots,\kmax-1$ we now define the sets
\begin{displaymath}
  \cA_i := \cL_{i+1} \setminus \Int_\cR\cL_i,
\end{displaymath}
and we refer to the sequence~$(\cA_i)_{i=0}^{\kmax-1}$ as the \emph{sequence of shifts}
of~$\cP$. These shifts are clearly non-empty and they are ordered in the following way.
\begin{prop}[Ordering of the shifts of~$\cP$]
\label{prop:j=i+1}
  In the above situation, assume that the inequality $\cA_i \cap \cA_j \neq
  \emptyset$ holds for some $0\leq i<j<\kmax$. Then one has to have $j=i+1$.
\end{prop}
\proof
First observe that 
\begin{displaymath}
  \cA_i \cap \cA_j =
  (\cL_{i+1}\setminus \Int_\cR\cL_i) \cap
    (\cL_{j+1}\setminus \Int_\cR\cL_j) =
  \cL_{i+1} \setminus \Int_\cR\cL_j .
\end{displaymath}
Now consider a cell $\tau \in \cA_i \cap \cA_j$. Since $\tau\not\in\Int_\cR\cL_j$,
we can select a $\sigma\in\cR \cap\Opn\tau\setminus\cL_j$. In view of the inequality
$i+1\leq j$, we have $\cL_{i+1}\subset \cL_j$. Hence, we also have $\sigma\not\in \cL_{i+1}$.
Together with the inclusion $\tau \in \cL_{i+1}\cap \Cl\sigma $, one further obtains
$\sigma\in\cB_{i+2} \subset \cL_{i+2}$. Finally, since $\sigma\not\in\cL_j$, this yields
$j < i+2$, i.e., we have $j = i+1$.
\qed\medskip

We also have the following two propositions.
\begin{prop}[Further properties of the~$\cL_k$]
\label{prop:Lk-int}
For arbitrary $0 < k \leq \kmax$ we have the inclusion $\cL_{k-1}\subset\Int_\cR\cL_k$.
\end{prop}
\proof
Assume to the contrary that there exists a cell $\tau\in\cL_{k-1}\setminus\Int_\cR\cL_k$.
Then there also has to be a cell $\rho\in\Opn_\cR\tau\setminus\cL_k\subset
\Opn_\cR\tau\setminus\cL_{k-1}$. It follows that $\rho\in\cB_k\subset\cL_k$,
which is a contradiction.
\qed

\begin{prop}[Properties of the~$\cA_i$]
\label{prop:extreme-i}
If $0 < i < \kmax-1$, then $\cH\cap\Cl\cA_i=\emptyset$.
\end{prop}
\proof
Assume to the contrary that there exists a cell $\tau\in\cH\cap\Cl\cA_i$.
Choose a $\sigma\in\cA_i$ such that $\tau\in\Cl\sigma$. Then
$\sigma\in\Opn\tau\subset\Opn\cH$ and
\begin{displaymath}
  \sigma\in\cL_{i+1}\setminus\Int_\cR\cL_i\subset \cR\cap\cL_{\kmax-1}\setminus\Int_\cR\cL_1.
\end{displaymath}
Hence, from \eqref{eq:kmax-1} we get $\sigma\in\cR\cap\Cl\cP=\cPb$.
Since $\cL_1=\cPb$ and $\sigma\not\in \Int_\cR\cL_1$, one can easily see that
$\Opn_\cR\sigma\not\subset\cPb$. Hence, $\sigma\in\Bd_\cR\cPb$ and
$\tau\in\cH\cap\Cl\Bd_\cR\cPb$, which contradicts \eqref{eq:bd-cpb}.
\qed
\begin{prop}[Further properties of the~$\cA_i$]
\label{prop:Cl-A-intersection}
   For arbitrary $i=0,1,\ldots,\kmax-1$ we have to have
   $\Cl\cA_i\cap\Cl\cA_{i+1}\neq\emptyset$, where~$i+1$ is taken
   modulo~$\kmax$, i.e., we assume $\cA_{\kmax}=\cA_0$.
\end{prop}
\proof
Assume first that the inequality $i<\kmax-1$ is satisfied, and select
a cell $\sigma\in\cB_{i+2}$. Then $\sigma\not\in\cL_{i+1}$, but we do have
$\cL_{i+1}\cap\Cl\sigma\neq\emptyset$. Hence, we can choose a
$\tau\in\cL_{i+1}\cap\Cl\sigma$. We have $\tau\in\Cl\sigma\subset\cB_{i+2}\subset\Cl\cA_{i+1}$,
yet one cannot have $\tau\in\Int_\cR\cL_i$, because then one would obtain
$\sigma\in\Opn_\cR\tau\subset\cL_i\subset\cL_{i+1}$, which is a contradiction.
But this in turn implies that we have $\tau\in \cL_{i+1}\setminus\Int_\cR\cL_i=
\cA_i\subset\Cl\cA_i$, which yields $\Cl\cA_i\cap\Cl\cA_{i+1}\neq\emptyset$
for $i<\kmax-1$.

Consider now the remaining case $i=\kmax-1$. Then we have both
$\cA_{i+1}=\cA_0=\cL_1=\cPb$ and $\cA_{i}=\cA_{\kmax-1}=\cR\setminus\Int_\cR\cL_{\kmax-1}$.
If one selects a cell $\sigma\in\Opn\cH\setminus\Cl\cP$, then
$\sigma\not\in\cL_{\kmax-1}$, because clearly $n_{\kmax-1}<\bar{n}$. It follows
that $\sigma\in\cR\setminus\Int_\cR\cL_{\kmax-1}=\cA_{\kmax-1}$. Since we also
have $\sigma\in\Opn\cH$, we can choose a $\tau\in\cH$ such that $\sigma \in\Opn\tau$.
Then $\tau\in\Cl\sigma\subset \Cl\cA_{\kmax-1}$, but also $\tau\in\cH\subset\Cl\cP=\Cl\cPb=\Cl\cA_0$.
Therefore, one obtains $\Cl\cA_0\cap \Cl\cA_{\kmax-1}\neq\emptyset$, which gives the
conclusion also for $i=\kmax-1$.
\qed\medskip

After these preparations we can now finally establish the following link between
the combinatorial setting and Theorem~\ref{thm:bricks}.
\begin{thm}[Periodic orbits via combinatorics]
\label{thm:poviacomb}
Assume that~$X$ is a cellular space with cellular decomposition~$\cX$, and that~$\phi$
is a strongly admissible flow on~$X$ with no cyclic intersection of toplexes. Let~$\cA$
be a connected isolated invariant set for the combinatorial multivector field~$\cV$
associated with~$\phi$ such that~$\cA$ contains only regular multivectors, and such that
for either $r = 0$ or $r = 1$ we have
\begin{displaymath}
 \dim H_{2n + r}(\Cl\cA,\Mo\cA) = \dim H_{2n+1 + r}(\Cl\cA,\Mo\cA)
 \quad\mbox{ for all }\quad n \in \ZZ,
\end{displaymath}
where not all of these homology groups are trivial. Assume further that
$\cP \subset \cA$ is a combinatorial Poincar\'e section for~$\cA$, and that the
associated sequence of shifts of~$\cP$ is given by~$(\cA_i)_{i=0}^{\kmax-1}$
with $\kmax \geq 3$. Then the collection $\cAtop := \cA\cap\cXtop$
is a brick decomposition of the isolating block $A := |\Cl\cA|$ for~$\phi$,
with a coarsening~$\overline{\cA}:=\{\cAtop_i\}_{i=0}^{\kmax-1}$, where
we define $\cAtop_i := \cA_{\kmax-1-i}\cap\cXtop$. Moreover, all the
assumptions of Theorem~\ref{thm:bricks} are satisfied. In consequence,
the set~$A$ contains a non-trivial periodic orbit of~$\phi$.
\end{thm}
\proof
We get from Theorem~\ref{thm:iso-block-from-cmvf}(i-ii) that~$A$ is indeed an
isolating block for~$\phi$. As for the assumptions of Theorem~\ref{thm:bricks},
note that~(a) follows from Proposition~\ref{prop:anrs}, while~(b) is a
consequence of Theorem~\ref{thm:iso-block-from-cmvf}(iii). According to
Propositions~\ref{prop:sharp-pair} and~\ref{prop:supp-cl-S} the set~$\cAtop$
is a brick decomposition of~$A$ with coarsening~$\overline{\cA}$. Then assumption~(c)
follows from the fact that we suppose $\kmax \geq 3$. To verify~(e), we define
$A_i := \bigcup_{\sigma \in \cAtop_i} \sigma =|\Cl\cA_i|$, and assume that the inequality $A_i \cap A_j
\neq \emptyset$ is satisfied for some $0 \leq i < j < \kmax$. Let $x\in A_i\cap A_j$
and consider the cell $\tau:=\cell{x}$. If we have $\tau\in\cR$, then
$\tau\in\Cl_\cR\cA_i\cap\Cl_\cR\cA_j=\cA_i\cap\cA_j$, because~$\cA_i$ and~$\cA_j$
are closed in~$\cR$. Therefore, we get from Proposition~\ref{prop:j=i+1} that $j=i+1$.
If $\tau\not\in\cR$, then we have the inclusion $\tau\in\cH$, because clearly
$\tau\in\Cl\cA$. By Proposition~\ref{prop:extreme-i} this is only possible if
$i=0$ and $j=\kmax-1$. Hence, in both of these cases~(e) follows. Assumption~(d)
follows easily from Proposition~\ref{prop:Cl-A-intersection} and  the definition
of the sets~$\cA_i$.

Finally, to verify property~(f) consider a maximal brick path $(\sigma_0,\sigma_1,\ldots)$.
As we mentioned before Proposition~\ref{prop:sharp-pair}, we may assume that the path
is sharp. Hence, by Proposition~\ref{prop:sharp-path} there is a reduced solution~$\gamma$
of~$\cV$ such that the toplexes along this path are precisely the toplexes
$(\sigma_0,\sigma_1,\ldots)$, and clearly we may assume that $\gamma(0)=\sigma_0\in\cAtop$.
Then there exists an index $k\in\{0,1,\ldots, \kmax-1\}$ such that $\sigma_0\in\cAtop_k\subset\cA_k\subset\cL_{k+1}$. 
Obviously, the solution~$\gamma$ is also a maximal solution. Therefore, by~\eqref{eq:P-sect-2}
we may choose a~$j$ such that~$\gamma(j)\in\cH\cup\Mo\cA$. If $\gamma(j)\in\Mo\cA$, then
clearly the corresponding toplex $\sigma$ satisfies $\sigma^-\subset A^-$. Hence, we
now assume $\gamma(j)\in\cH$. Without loss of generality, we may further suppose
that~$j$ is chosen as the first index with this property. Then one has
$\gamma([0,j-1])\subset\cR$, and since $\gamma(0)=\sigma_0\in\cL_{k+1}$, 
we also get from Proposition~\ref{prop:L} the inclusion $\gamma([0,j-1])\subset
\cL_{k+1}$. Now let~$i$ denote the first index such that $\gamma(i)\not\in\cA_k$.
Such an index exists and satisfies $i\leq j$, because both the inclusion
$\gamma(j)\in\cH$ and the inclusion $\cH\cap\cA_k\subset\cH\cap\cR=\emptyset$
are satisfied. It suffices to verify that $\gamma(i)\in\cA_{k-1}$ holds.
If $k=0$, then we have $i=j$, as well as $\gamma(i)\in\cH\subset\cA_{\kmax-1}$.
Hence, we can now assume $k>0$. In this case one has $\gamma(i)\in F_\cV(\gamma(i-1))=
\Cl \gamma(i-1)\cup [\gamma(i-1)]_\cV$. We cannot have $\gamma(i)\in  \Cl \gamma(i-1)$,
because $\gamma(i-1)\in\cA_k$ and then $\gamma(i)\in\Cl_\cR\cA_k=\cA_k$, a contradiction.
Therefore, $\gamma(i)\in [\gamma(i-1)]_\cV$. Note that $\gamma(i-1)\in\cA_k=\cL_{k+1}
\setminus\Int_\cR\cL_k$. In particular, $\gamma(i-1)\in\cL_{k+1}$ and by Proposition~\ref{prop:cL}
also $\gamma(i)\in\cL_{k+1}$. But, $\gamma(i)\not\in\cA_{k}$. Therefore, $\gamma(i)\in
\Int_\cR\cL_k\subset\cL_k$. We will show that $\gamma(i)\in\cA_{k-1}$. Indeed, if
$\gamma(i)\not\in\cA_{k-1}$, then $\gamma(i)\in\Int_\cR\cL_{k-1}$, because
$\gamma(i)\in\cL_k$. Then, by $\cV$-compatibility of $\cL_{k-1}$ we get
$\gamma(i-1)\in\cL_k$ and by Proposition~\ref{prop:Lk-int} we obtain
$\gamma(i-1)\in\Int_\cR\cL_k$. Therefore, $\gamma(i-1)\not\in\cA_k$, which is a
contradiction and proves that $\gamma(i)\in\cA_{k-1}$.
\qed
%
% References
%

%\addcontentsline{toc}{section}{References}
\footnotesize
%
%   Add in the bbl-file the command \parskip 0pt.
%
\bibliography{wanner1a,wanner1b,wanner2a,wanner2b,wanner2c,extra}
\bibliographystyle{abbrv}
\end{document}